%Paper: alg-geom/9410022
%From: demailly@fourier.grenet.fr (Jean Pierre Demailly)
%Date: Fri, 21 Oct 1994 15:17:10 +0100 (MET)
%Date (revised): Tue, 8 Nov 1994 13:53:30 +0100 (MET)
%Date (revised): Tue, 28 Feb 1995 18:04:47 +0100 (MET)

% CIME Course, July 1994, "Transcendental Methods in Algebraic Geometry"
% $L^2$ Vanishing Theorems for Positive Line Bundles and Adjunction Theory
% J.P. Demailly, Universit\'e de Grenoble I, Saint Martin d'H\`eres, France

%%% Macros

\font \tenbf                  = cmb10

\font \tbfontss               = cmbx5  scaled\magstep1
\font \sixbf                  = cmbx6

\font \tbfonts                = cmbx7  scaled\magstep1
\font \ninebf                 = cmbx9
\font \bxf                    = cmbx10
\font \tbfontt                = cmbx10 scaled\magstep1
\font \tafontt                = cmbx10 scaled\magstep2
\font \tasys                  = cmex10 scaled\magstep1

\font \sixi                   = cmmi6
\font \ninei                  = cmmi9
\font \kleinhalbcurs          = cmmib10 scaled 833
\font \tamss                  = cmmib10
\font \tams                   = cmmib10 scaled\magstep1

\font \sixrm                  = cmr6
\font \ninerm                 = cmr9

\font \ninesl                 = cmsl9
\font \tbst                   = cmsy10 scaled\magstep1
\font \tbsss                  = cmsy5  scaled\magstep1

\font \sixsy                  = cmsy6
\font \tbss                   = cmsy7  scaled\magstep1

\font \ninesy                 = cmsy9

\font \nineit                 = cmti9
\font \ninett                 = cmtt9

\catcode`\@=11
\def\hexnumber@#1{\ifnum#1<10 \number#1\else
 \ifnum#1=10 A\else\ifnum#1=11 B\else\ifnum#1=12 C\else
 \ifnum#1=13 D\else\ifnum#1=14 E\else\ifnum#1=15 F\fi\fi\fi\fi\fi\fi\fi}
\def\mathhexbox@#1#2#3{\text{$\m@th\mathchar"#1#2#3$}}
\def\text{\relaxnext@\ifmmode\let\next\text@\else\let\next\text@@\fi\next}
\def\text@@#1{\leavevmode\hbox{#1}}
\font\tenmsx=msxm10
\font\sevenmsx=msxm7
\font\fivemsx=msxm5
\font\tenmsy=msym10
\font\sevenmsy=msym7
\font\fivemsy=msym5
\font\teneuf=eufm10
\font\seveneuf=eufm7
\font\fiveeuf=eufm5
\newfam\msxfam
\newfam\msyfam
\newfam\euffam
\textfont\msxfam=\tenmsx
\scriptfont\msxfam=\sevenmsx
\scriptscriptfont\msxfam=\fivemsx
\textfont\msyfam=\tenmsy
\scriptfont\msyfam=\sevenmsy
\scriptscriptfont\msyfam=\fivemsy
\textfont\euffam=\teneuf
\scriptfont\euffam=\seveneuf
\scriptscriptfont\euffam=\fiveeuf
\def\msx@{\hexnumber@\msxfam}
\def\msy@{\hexnumber@\msyfam}
\mathchardef\boxtimes="2\msx@02
\mathchardef\square="0\msx@03
\mathchardef\twoheadrightarrow="3\msx@10
\mathchardef\twoheadleftarrow="3\msx@11
\mathchardef\upharpoonright="3\msx@16
\let\restriction\upharpoonright
\mathchardef\smallsmile="3\msx@60
\mathchardef\smallfrown="3\msx@61
\mathchardef\complement="0\msx@7B
\mathchardef\ltimes="2\msy@6E
\mathchardef\rtimes="2\msy@6F
\mathchardef\subsetneq="3\msy@28
\mathchardef\supsetneq="3\msy@29
\mathchardef\shortmid="3\msy@70
\mathchardef\shortparallel="3\msy@71
\mathchardef\ssm="2\msy@72
\def\Bbb{\relaxnext@\ifmmode\let\next\Bbb@\else
 \def\next{\errmessage{Use \string\Bbb\space only in math mode}}\fi\next}
\def\Bbb@#1{{\Bbb@@{#1}}}
\def\Bbb@@#1{\fam\msyfam#1}
\def\goth{\relaxnext@\ifmmode\let\next\goth@\else
 \def\next{\errmessage{Use \string\goth\space only in math mode}}\fi\next}
\def\goth@#1{{\goth@@{#1}}}
\def\goth@@#1{\fam\euffam#1}
\def\relaxnext@{\let\next\relax}
\catcode`\@=12

\magnification=\magstep1
\vsize=24.2true cm
\hsize=15.3true cm
\parskip=2pt plus 2pt minus 0pt
\hfuzz=2pt
\tolerance=500
\abovedisplayskip=3 mm plus6pt minus 4pt
\belowdisplayskip=3 mm plus6pt minus 4pt
\abovedisplayshortskip=0mm plus6pt minus 2pt
\belowdisplayshortskip=2 mm plus4pt minus 4pt
\predisplaypenalty=0
\frenchspacing
\newdimen\oldparindent\oldparindent=20pt
\def\newline{\hfil\break}
\skewchar\ninei='177 \skewchar\sixi='177
\skewchar\ninesy='60 \skewchar\sixsy='60
\hyphenchar\ninett=-1

  \mathchardef\Gamma="0100
  \mathchardef\Delta="0101
  \mathchardef\Theta="0102
  \mathchardef\Lambda="0103
  \mathchardef\Xi="0104
  \mathchardef\Pi="0105
  \mathchardef\Sigma="0106
  \mathchardef\Upsilon="0107
  \mathchardef\Phi="0108
  \mathchardef\Psi="0109
  \mathchardef\Omega="010A

\def\bbbr{{\Bbb R}}
\def\bbbn{{\Bbb N}}

\def\bbbc{{\Bbb C}}
\def\bbbz{{\Bbb Z}}
\def\bbbp{{\Bbb P}}
\def\bbbq{{\Bbb Q}}
\def\bbbone{\hbox{\bf 1}}

\def\begpet{\vskip6pt\bgroup\eightpoint}
\def\endpet{\vskip6pt\egroup}

\def\vec#1{\hbox{\textfont1=\bxf\textfont0=\bxf$#1$}}

\def\eightpoint{\def\rm{\fam0\ninerm}%
\textfont0=\ninerm \scriptfont0=\sixrm \scriptscriptfont0=\fiverm
 \textfont1=\ninei \scriptfont1=\sixi \scriptscriptfont1=\fivei
 \textfont2=\ninesy \scriptfont2=\sixsy \scriptscriptfont2=\fivesy
 \def\it{\fam\itfam\nineit}%
 \textfont\itfam=\nineit
 \def\sl{\fam\slfam\ninesl}%
 \textfont\slfam=\ninesl
 \def\bf{\fam\bffam\ninebf}%
 \textfont\bffam=\ninebf \scriptfont\bffam=\sixbf
 \scriptscriptfont\bffam=\fivebf
 \normalbaselineskip=9pt
 \setbox\strutbox=\hbox{\vrule height7pt depth2pt width0pt}%
 \normalbaselines\rm
\def\vec##1{\setbox0=\hbox{$##1$}\hbox{\hbox
to0pt{\copy0\hss}\kern0.45pt\box0}}}%

\headline={\eightpoint\def\newline{ }\def\fonote#1{}\ifodd\pageno
\hfil\botmark\unskip\kern1.4true cm\llap{\folio}\else\leftheadline\fi}
\def\leftheadline{\rlap{\folio}\kern1.4true cm Chaptertitle\hfil}
\mark{Paragraphtitle}
\nopagenumbers
\output={\if N\header\headline={\hfil}\fi\plainoutput
\global\let\header=Y}

\let\lasttitle=N
\let\header=N
 \def\titlea#1#2{%
 {\baselineskip=18pt
  \lineskip=18pt
  \tafontt
  \noindent#1#2}}

\def\leftrunningtitle#1{\gdef\leftheadline{\rlap{\folio}\kern1.4true cm
#1\hfil\ignorespaces}}

\def\checkleftspace#1{\dimen0=\pagetotal
    \ifdim\dimen0<\pagegoal
       \dimen0=\ht0\advance\dimen0 by\dp0\advance\dimen0 by
       #1\normalbaselineskip
       \advance\dimen0 by\pagetotal
       \advance\dimen0 by-\pageshrink
       \ifdim\dimen0>\pagegoal
          \vfill\eject
       \fi
    \fi}

 \def\titleb#1#2{%
     \if N\lasttitle\else\vskip-29pt
     \fi
     \vskip25pt plus 4pt minus4pt
     \bgroup
 \textfont0=\tbfontt \scriptfont0=\tbfonts \scriptscriptfont0=\tbfontss
 \textfont1=\tams \scriptfont1=\tamss \scriptscriptfont1=\kleinhalbcurs
 \textfont2=\tbst \scriptfont2=\tbss \scriptscriptfont2=\tbsss
 \textfont3=\tasys \scriptfont3=\tenex \scriptscriptfont3=\tenex
     \baselineskip=16pt
     \lineskip=16pt
     \raggedright
     \pretolerance=10000
     \tbfontt
     \setbox0=\vbox{\vskip25pt
     \def\fonote##1{}%
     \noindent
     \ignorespaces#1\ #2
     \vskip15pt}%
     \checkleftspace 5
     \noindent\ignorespaces #1\ #2
     \vskip15true pt plus4pt minus4pt\egroup
     \setbox0=\hbox{\eightpoint\def\newline{ }\def\fonote##1{}%
     \kern1.4true cm#1\ #2}\ifdim\wd0>\hsize
     \message{Your TITLEB exceeds the headline, please use a
     short form with TITLEBRUNNING}\mark{Title of section
     suppressed due to excessive length}%
     \else{\def\newline{ }\def\fonote##1{}\mark{#1\ #2}}\fi
     \nobreak
     \global\let\lasttitle=B%
     \parindent=0pt
     \everypar={\global\parindent=\oldparindent
     \global\let\lasttitle=N\global\everypar={}}%
     \ignorespaces}

 \def\titlec#1#2{%
     \if N\lasttitle\else\vskip-3pt\vskip-\baselineskip
     \fi
     \vskip17pt plus 4pt minus 4pt
     \bgroup
     \bxf
     \raggedright
     \pretolerance=10000
     \setbox0=\vbox{\vskip 17pt
     \def\fonote##1{}%
     \noindent
     \ignorespaces#1\ #2
     \vskip10pt}%
     \checkleftspace 4
     \noindent
     \ignorespaces #1\ #2
     \vskip10pt plus4pt minus4pt\egroup
     \nobreak
     \global\let\lasttitle=C%
     \parindent=0pt
     \everypar={\global\parindent=\oldparindent
     \global\let\lasttitle=N\global\everypar={}}%
     \ignorespaces}

 \def \titled#1{
     \if N\lasttitle\else\vskip-3pt\vskip-\baselineskip
     \fi
     \vskip7true pt plus 4pt minus 4pt
     \bgroup
     \bf
     \noindent
     \ignorespaces #1\unskip.\ \egroup
     \ignorespaces}

\def\footnoterule{\kern-3pt\hrule width 2true cm\kern2.6pt}

\newcount\footcount \footcount=0
\def\advftncnt{\advance\footcount by1\global\footcount=\footcount}

\def\fonote#1{\advftncnt$^{\the\footcount}$\begingroup\eightpoint
       \def\textindent##1{\hangindent0.5\oldparindent\noindent\hbox
       to0.5\oldparindent{##1\hss}\ignorespaces}%
\vfootnote{$^{\the\footcount}$}{#1}\endgroup}

\def\item#1{\par\parindent=\oldparindent\noindent
\hangindent\parindent\hbox to\parindent{#1\hss}\ignorespaces}

\def\qed{\ifmmode\square\else{\unskip\nobreak\hfil
\penalty50\hskip1em\null\nobreak\hfil\hbox{$\square$}
\parfillskip=0pt\finalhyphendemerits=0\endgraf}\fi}

\def\ii{{\rm i}}
\def\bu{{\scriptstyle\bullet}}
\def\ld{,\ldots,}
\let\lra=\longrightarrow

\let\wt=\widetilde
\let\ovl=\overline
\def\ovlp{{\overline\partial}}
\def\soul{{\raise-4.6pt\hbox{$-$}\kern-7pt}}
\def\ci{C^\infty}
\def\compact{{\subset\!\subset}}

\def\gge{\raise1.5pt\hbox{$\scriptscriptstyle\ge$}}
\def\dasharrow{\mathrel{\hbox{
  \kern1pt \vrule height2.45pt depth-2.15pt width2pt}
  \kern1pt{\vrule height2.45pt depth-2.15pt width2pt}
  \kern1pt{\vrule height2.45pt depth-2.15pt width2pt}
  \kern1pt{\vrule height2.45pt depth-2.15pt width1.7pt\kern-1.7pt}
  {\raise1.4pt\hbox{$\scriptscriptstyle\succ$}}}}
\def\lraww{\mathrel{\rlap{$\longrightarrow$}\kern-1pt\longrightarrow}}

\def\buildo#1\over#2{\mathrel{\mathop{\null#2}\limits^{#1}}}
\def\buildu#1\under#2{\mathrel{\mathop{\null#2}\limits_{#1}}}

\def\Ker{\mathop{\rm Ker}\nolimits}

\def\Tr{\mathop{\rm Tr}\nolimits}
\def\deg{\mathop{\rm deg}\nolimits}

\def\Im{\mathop{\rm Im}\nolimits}

\def\codim{\mathop{\rm codim}\nolimits}
\def\Supp{\mathop{\rm Supp}\nolimits}

\def\Hom{\mathop{\rm Hom}\nolimits}

\def\Herm{\mathop{\rm Herm}\nolimits}
\def\rank{\mathop{\rm rank}\nolimits}

\def\Psh{\mathop{\rm Psh}\nolimits}
\def\Pic{\mathop{\rm Pic}\nolimits}
\def\Psh{\mathop{\rm Psh}\nolimits}

\def\pr{\mathop{\rm pr}\nolimits}
\def\ord{\mathop{\rm ord}\nolimits}

\def\DR{{\rm DR}}
\def\FS{{\rm FS}}
\def\NS{{\rm NS}}

\def\amp{{\rm amp}}
\def\nef{{\rm nef}}
\def\eff{{\rm eff}}
\def\psef{{\rm psef}}
\def\ac{{\rm ac}}
\def\loc{{\rm loc}}
\def\reg{{\rm reg}}
\def\sing{{\rm sing}}

\def\gm{{\goth m}}

\def\cHom{{\cal H}{\it om}}

\catcode`@=11
\def\cmalign#1{\null\,\vcenter{\normalbaselines\m@th
    \ialign{\hfil$##$&&$##$\hfil\crcr
      \mathstrut\crcr\noalign{\kern-\baselineskip}
      #1\crcr\mathstrut\crcr\noalign{\kern-\baselineskip}}}\,}
\catcode`@=12

\long\def\begstat#1 #2\endstat
{\removelastskip\vskip\baselineskip\noindent{\bf#1.}
{\it\ignorespaces#2}\vskip\baselineskip}

\def\nlni{%
   \par
   \if N\lasttitle
      \ifvmode \removelastskip \fi
      \bigskip
   \fi
   \noindent}
\def\newenvironment#1#2#3#4#5#6%
{\expandafter\gdef\csname beg#1\endcsname##1{#2%
{#3\csname #1\endcsname\if!##1!.\else\ \ignorespaces
##1\unskip\fi\kern0.4em}\bgroup
#4%
\ignorespaces
}%
\expandafter\gdef\csname #1\endcsname{#6}%
\expandafter\gdef\csname end#1\endcsname{\endgraf\egroup#5}%
}
\newenvironment{proof}{\nlni}{\it}{\rm}{\vskip\baselineskip}{Proof}

\def\cA{{\cal A}}

\def\cC{{\cal C}}
\def\cD{{\cal D}}
\def\cE{{\cal E}}

\def\cH{{\cal H}}
\def\cI{{\cal I}}
\def\cJ{{\cal J}}

\def\cL{{\cal L}}

\def\cO{{\cal O}}

\def\ptstep{{\vrule height 0.4pt depth 0pt width 0.4pt}}
\def\verl{{\vrule width 0.4pt height 12pt depth 6pt}}
\def\table#1{\vbox{\baselineskip=0pt\lineskip=0pt
    \halign{$\displaystyle{}##\hfil{}$\hfil&&
    $\displaystyle{}##\hfil{}$\hfil\crcr#1\crcr}}}

\catcode`\@=11
\font\eightrm=cmr8
\font\eighti=cmmi8
\font\eightsy=cmsy8
\font\eightbf=cmbx8

\font\eightit=cmti8
\font\eightsl=cmsl8
\font\sixrm=cmr6
\font\sixi=cmmi6
\font\sixsy=cmsy6
\font\sixbf=cmbx6

\font\smallcap=cmcsc10

\def\eightpoint{%
  \textfont0=\eightrm \scriptfont0=\sixrm \scriptscriptfont0=\fiverm
  \def\rm{\fam\z@\eightrm}%
  \textfont1=\eighti \scriptfont1=\sixi \scriptscriptfont1=\fivei
  \def\oldstyle{\fam\@ne\eighti}%
  \textfont2=\eightsy \scriptfont2=\sixsy \scriptscriptfont2=\fivesy
  \textfont\itfam=\eightit
  \def\it{\fam\itfam\eightit}%
  \textfont\slfam=\eightsl
  \def\sl{\fam\slfam\eightsl}%
  \textfont\bffam=\eightbf \scriptfont\bffam=\sixbf
  \scriptscriptfont\bffam=\fivebf
  \def\bf{\fam\bffam\eightbf}%
  \abovedisplayskip=9pt plus 2pt minus 6pt
  \abovedisplayshortskip=0pt plus 2pt
  \belowdisplayskip=9pt plus 2pt minus 6pt
  \belowdisplayshortskip=5pt plus 2pt minus 3pt
  \smallskipamount=2pt plus 1pt minus 1pt
  \medskipamount=4pt plus 2pt minus 1pt
  \bigskipamount=9pt plus 3pt minus 3pt
  \normalbaselineskip=9pt
  \setbox\strutbox=\hbox{\vrule height7pt depth2pt width0pt}%
  \let\bigf@ntpc=\eightrm \let\smallf@ntpc=\sixrm
  \normalbaselines\rm}

\def\tenpoint{%
  \textfont0=\tenrm \scriptfont0=\sevenrm \scriptscriptfont0=\fiverm
  \def\rm{\fam\z@\tenrm}%
  \textfont1=\teni \scriptfont1=\seveni \scriptscriptfont1=\fivei
  \def\oldstyle{\fam\@ne\teni}%
  \textfont2=\tensy \scriptfont2=\sevensy \scriptscriptfont2=\fivesy
  \textfont\itfam=\tenit
  \def\it{\fam\itfam\tenit}%
  \textfont\slfam=\tensl
  \def\sl{\fam\slfam\tensl}%
  \textfont\bffam=\tenbf \scriptfont\bffam=\sevenbf
  \scriptscriptfont\bffam=\fivebf
  \def\bf{\fam\bffam\tenbf}%
  \textfont\ttfam=\tentt
  \def\tt{\fam\ttfam\tentt}%
  \abovedisplayskip=6pt plus 2pt minus 6pt
  \abovedisplayshortskip=0pt plus 3pt
  \belowdisplayskip=6pt plus 2pt minus 6pt
  \belowdisplayshortskip=7pt plus 3pt minus 4pt
  \smallskipamount=3pt plus 1pt minus 1pt
  \medskipamount=6pt plus 2pt minus 2pt
  \bigskipamount=12pt plus 4pt minus 4pt
  \normalbaselineskip=12pt
  \setbox\strutbox=\hbox{\vrule height8.5pt depth3.5pt width0pt}%
  \let\bigf@ntpc=\tenrm \let\smallf@ntpc=\sevenrm
  \normalbaselines\rm}
\catcode`\@=12

\def\pointdash{\discretionary{.}{}{.\kern.35em---\kern.7em}}
\newif\ifdigit
\def\digit{\digittrue}
\digit
\newdimen\laenge
\def\letter#1|{\setbox3=\hbox{#1}\laenge=\wd3\advance\laenge by 3mm
\digitfalse}

\def\article#1|#2|#3|#4|#5|#6|#7|%
    {{\ifdigit\leftskip=7mm\noindent
     \hangindent=2mm\hangafter=1
\llap{[#1]\hskip1.35em}{\smallcap #2}\pointdash {\it #3}, {\rm #4},
\nobreak{\bf #5} ({\oldstyle #6}), \nobreak #7.\par\else\noindent
\advance\laenge by 4mm \hangindent=\laenge\advance\laenge by -4mm\hangafter=1
\rlap{[#1]}\hskip\laenge{\smallcap #2}\pointdash
{\it #3}, #4, {\bf #5} ({\oldstyle #6}), #7.\par\fi}}

\def\bookref#1|#2|#3|#4|#5|%
    {{\ifdigit\leftskip=7mm\noindent
    \hangindent=2mm\hangafter=1
\llap{[#1]\hskip1.35em}{\smallcap #2}\pointdash{\it #3}, #4, {\oldstyle
#5}.\par
\else\noindent
\advance\laenge by 4mm \hangindent=\laenge\advance\laenge by -4mm
\hangafter=1
\rlap{[#1]}\hskip\laenge{\smallcap #2}\pointdash
{\it #3}, #4, {\oldstyle #5}.\par\fi}}

\def\miscell#1|#2|#3|#4|%
    {{\ifdigit\leftskip=7mm\noindent
    \hangindent=2mm\hangafter=1
     \llap{[#1]\hskip1.35em}{\smallcap #2}\pointdash{\it #3}, {\rm #4}.\par
\else\noindent
\advance\laenge by 4mm \hangindent=\laenge\advance\laenge by -4mm
\hangafter=1
\rlap{[#1]}\hskip\laenge{\smallcap #2}\pointdash{\it #3}, {\rm #4}.\par\fi}}

%%% Main File

\def\twolines{\newline\phantom{$~$}\hfill\newline}

\titlea{}{CIME Session
\twolines Transcendental Methods in Algebraic Geometry
\twolines Cetraro, Italy, July 1994
\twolines\phantom{$~$}\hfill\newline
L\kern-3.5pt\raise7pt\hbox{\tenbf 2} Vanishing Theorems for Positive
\twolines Line Bundles and Adjunction Theory}
\leftrunningtitle{J.-P. Demailly,~~Positive Line Bundles and Adjunction Theory}

\vskip20pt\noindent
\hskip2cm{\tbfontt Lectures by Jean-Pierre Demailly}
\vskip3pt\noindent
\hskip2cm{\tbfontt Universit\'e de Grenoble I, Institut Fourier}

\titlec{}{\phantom{$~$}\newline\tbfontt Contents}
0. Introduction
   \dotfill p.~1\break
1. Preliminary Material
   \dotfill p.~4\break
2. Lelong Numbers and Intersection Theory
   \dotfill p.~13\break
3. Holomorphic Vector Bundles, Connections and Curvature
   \dotfill p.~22\break
4. Bochner Technique and Vanishing Theorems
   \dotfill p.~26\break
5. $L^2$ Estimates and Existence Theorems
   \dotfill p.~31\break
6. Numerically Effective Line Bundles
   \dotfill p.~40\break
7. Seshadri Constants and the Fujita Conjecture
  \dotfill p.~47\break
8. Algebraic Approach of the Fujita Conjecture
  \dotfill p.~53\break
9. Regularization of Currents and Self-intersection Inequalities
  \dotfill p.~62\break
10. Use of Monge-Amp\`ere Equations
  \dotfill p.~68\break
11. Numerical Criteria for Very Ample Line Bundles
  \dotfill p.~74\break
12. Holomorphic Morse Inequalities
  \dotfill p.~82\break
13. Effective Version of Matsusaka's Big Theorem
  \dotfill p.~86\break
References
  \dotfill p.~91

\titleb{0.}{Introduction}
Transcendental methods of algebraic geometry have been extensively
studied since a very long time, starting with the work of Abel, Jacobi
and Riemann in the nineteenth century. More recently, in the period
1940-1970, the work of Hodge, Hirzebruch, Kodaira, Atiyah revealed
still deeper relations between complex analysis, topology, PDE theory
and algebraic geometry. In the last ten years, gauge theory has proved
to be a very efficient tool for the study of many important questions:
moduli spaces, stable sheaves, non abelian Hodge theory, low
dimensional topology $\ldots$

Our main purpose here is to describe a few analytic tools which are useful
to study questions such as linear series and vanishing theorems for algebraic
vector bundles. One of the early success of analytic methods in this context
is Kodaira's use of the Bochner technique in relation with the theory of
harmonic forms, during the decade 1950-60. The idea is to represent
cohomology classes by harmonic forms and to prove vanishing theorems by
means of suitable a priori curvature estimates. The prototype of such
results is the Akizuki-Kodaira-Nakano theorem (1954): if $X$ is a
nonsingular projective algebraic variety and $L$ is a holomorphic line
bundle on $X$ with positive curvature, then $H^q(X,\Omega^p_X\otimes L)=0$
for $p+q>\dim X$ (throughout the paper we set $\Omega^p_X=\Lambda^p
T^\star_X$ and $K_X=\Lambda^nT^\star_X$, $n=\dim X$, viewing these objects
either as holomorphic bundles or as locally free $\cO_X$-modules).
It is only much later that an algebraic proof of this result has been
proposed by Deligne-Illusie, via characteristic $p$ methods, in 1986.

A refinement of the Bochner technique leads to H\"ormander's $L^2$ estimates
(1965), which are best expressed in the geometric form given by
Andreotti-Vesentini [AV65]. Not only vanishing theorems are proved, but more
precise information of a quantitative nature are obtained about
solutions of $\ovlp$-equations. More explicitly, suppose that the bundle
$L$ is equipped with a hermitian metric of weight $e^{-2\varphi}$, where
$\varphi$ is a (locally defined) plurisubharmonic function; then explicit
bounds on the $L^2$ norm $\int_X|f|^2e^{-2\varphi}$ of solutions is obtained.
The result is still more useful if the plurisubharmonic weight $\varphi$
is allowed to have singularities. Following Nadel [Nad89], one defines
the {\it multiplier ideal sheaf} $\cI(\varphi)$ to be the sheaf of germs of
holomorphic functions $f$ such that $|f|^2e^{-2\varphi}$ is locally summable.
Then $\cI(\varphi)$ is a coherent algebraic sheaf over $X$ and
$H^q(X,K_X\otimes L\otimes\cI(\varphi))=0$ for all $q\ge 1$ if the
curvature of $L$ is positive (as a current). This important result
can be seen as a generalization of the Kawamata-Viehweg vanishing theorem
([Kaw82], [Vie82]), which is one of the cornerstones of higher dimensional
algebraic geometry (especially of Mori's minimal model program).

In the dictionary between analytic geometry and algebraic geometry, the
ideal $\cI(\varphi)$ plays a very important role, since it directly converts
an analytic object into an algebraic one, and, simutaneously, takes care
of the singularities in a very efficient way. Another analytic tool used to
deal with singularities is the theory of positive currents introduced by
Lelong [Lel57]. Currents can be seen as generalizations of algebraic cycles,
and many classical results of intersection theory still apply to currents.
The concept of Lelong number of a current is the analytic analogue of the
concept of multiplicity of a germ of algebraic variety. Intersections
of cycles correspond to wedge products of currents (whenever these products
are defined). A convenient measure of local positivity of a holomorphic line
can be defined in this context: the {\it Seshadri constant} of a line bundle
at a point is the largest possible Lelong number for a singular metric of
positive curvature assuming an isolated singularity at the given point
(see [Dem90]). Seshadri constants can also be given equivalent purely
algebraic definitions. We refer to Ein-Lazarsfeld [EL92] and
Ein-K\"uchle-Lazarsfeld [EKL94] for very interesting new results
concerning Seshadri constants.

One of our main motivations has been the study of the following conjecture
of Fujita: if $L$ is an ample (i.e.\ positive) line bundle on a projective
$n$-dimensional algebraic variety $X$, then $K_K+(n+2)L$ is very ample. A
major result obtained by Reider [Rei88] is a proof of the Fujita conjecture
in the case of surfaces (the case of curves is easy).
Reider's approach is based on Bogomolov's inequality for stable vector
bundles and the results obtained are almost optimal. Unfortunately, it
seems difficult to extend Reider's original method to higher dimensions.
In the analytic approch, which works for arbitrary dimensions, one tries
to construct a suitable (singular) hermitian metric on $L$ such that the
the ideal $\cI(\varphi)$ has a given $0$-dimensional subscheme of $X$ as
its zero variety. As we showed in [Dem93b], this can be done essentially
by solving a complex Monge-Amp\`ere equation
$$(\ii d'd''\varphi)^n=
\hbox{\rm linear combination of Dirac measures},$$
via the Aubin-Calabi-Yau theorem ([Aub78], [Yau78]). The solution $\varphi$
then assumes logarithmic poles and the difficulty is to force the
singularity to be an isolated pole; this is the point where intersection
theory of currents is useful. In this way, we can prove e.g.\ that
$2K_X+L$ is very ample under suitable numerical conditions for~$L$.
Alternative algebraic techniques have been developed recently by Koll\'ar
[Kol92], Ein-Lazarsfeld [EL93], Fujita [Fuj93] and [Siu94a, $\,$b]. The
basic idea is to apply the Kawamata-Viehweg vanishing theorem, and to
use the Riemann-Roch formula instead of the Monge-Amp\`ere equation.
The proofs proceed with careful inductions on dimension,
together with an analysis of the base locus of the linear systems
involved. Although the results obtained in low dimensions are slightly
more precise than with the analytic method, it is still not clear
whether the range of applicability of the methods are exactly the same.
Because it fits well with our approach, we have included here a simple
algebraic method due to Y.T.~Siu [Siu94a], showing that $2K_X+mL$
is very ample for $m\ge 2+{3n+1\choose n}$.

Our final concern in these notes is a proof of the effective Matsusaka big
theorem obtained by [Siu93]. Siu's result is the existence of an effective
value $m_0$ depending only on the intersection numbers $L^n$ and $L^{n-1}\cdot
K_X$, such that $mL$ is very ample for $m\ge m_0$. The basic idea is to
combine results on the very ampleness of $2K_X+mL$ together with
the theory of holomorphic Morse inequalities ([Dem85b]). The Morse
inequalities are used to construct sections of $m'L-K_X$ for $m'$ large.
Again this step can be made algebraic (following suggestions by F.~Catanese
and R.~Lazarsfeld), but the analytic formulation apparently has a wider
range of applicability.

These notes are essentially written with the idea of serving as an
analytic toolbox for algebraic geometry. Although efficient algebraic
techniques exist, our feeling is that the analytic techniques are very
flexible and offer a large variety of guidelines for more algebraic
questions (including applications to number theory which are not
discussed here). We made a special effort to use as little prerequisites
and to be as self-contained as possible; hence the rather long
preliminary sections dealing with basic facts of complex differential
geometry. The reader wishing to have a presentation of the algebraic
approach to vanishing theorems and linear series is referred to the
excellent notes written by R.~Lazarsfeld [Laz93]. In the last
years, there has been a continuous and fruitful interplay between the
algebraic and analytic viewpoints on these questions, and I have
greatly benefitted from observations and ideas contained in the works of
J.~Koll\'ar, L.~Ein, R.~Lazarsfeld and Y.T.~Siu. I would like to thank
them for their interest in my work and for their encouragements.
\vfill\eject

\titleb{1.}{Preliminary Material}
\titlec{1.A.}{Dolbeault Cohomology and Sheaf Cohomology}
Let $X$ be a $\bbbc$-analytic manifold of dimension $n$. We denote by
$\Lambda^{p,q}T^\star_X$ the bundle of differential forms of bidegree
$(p,q)$ on $X$, i.e., differential forms which can be written as
$$u=\sum_{|I|=p,\,|J|=q}u_{I,J}dz_I\wedge d\ovl z_J.$$
Here $(z_1\ld z_n)$ denote arbitrary local holomorphic coordinates,
\hbox{$I=(i_1\ld i_p)$}, \hbox{$J=(j_1\ld j_q)$} are multiindices
(increasing sequences of integers in the range $[1\ld n]$, of lengths
$|I|=p$, $|J|=q$), and
$$dz_I:=dz_{i_1}\wedge\ldots\wedge dz_{i_p},\qquad
d\ovl z_J:=d\ovl z_{j_1}\wedge\ldots\wedge d\ovl z_{j_q}.$$
Let $\cE^{p,q}$ be the sheaf of germs of complex valued differential
$(p,q)$-forms with $C^\infty$ coefficients. Recall that the exterior
derivative $d$ splits as $d=d'+d''$ where
$$\eqalign{
d'u&=\sum_{|I|=p,\,|J|=q,1\le k\le n}{\partial u_{I,J}
\over\partial z_k}dz_k\wedge dz_I\wedge d\ovl z_J,\cr
d''u&=\sum_{|I|=p,\,|J|=q,1\le k\le n}{\partial u_{I,J}
\over\partial\ovl z_k}d\ovl z_k\wedge dz_I\wedge d\ovl z_J\cr}$$
are of type $(p+1,q)$, $(p,q+1)$ respectively. The well-known
Dolbeault-Grothendieck lemma asserts that any $d''$-closed
form of type $(p,q)$ with $q>0$ is locally $d''$-exact (this is the
analogue for $d''$ of the usual Poincar\'e lemma for~$d$, see e.g.\
H\"ormander 1966). In other words, the complex of sheaves
$(\cE^{p,\bu},d'')$ is exact in degree $q>0$; in degree $q=0$, $\Ker d''$
is the sheaf $\Omega^p_X$ of germs of holomorphic forms of degree $p$
on~$X$.

More generally, if $F$ is a holomorphic vector bundle
of rank $r$ over~$X$, there is a natural $d''$ operator acting
on the space $C^\infty(X,\Lambda^{p,q}T^\star_X\otimes F)$ of smooth
$(p,q)$-forms with values in~$F$; if $s=\sum_{1\le\lambda\le r}s_\lambda
e_\lambda$ is a $(p,q)$-form expressed in terms of a local holomorphic
frame of $F$, we simply define $d''s:=\sum d''s_\lambda\otimes e_\lambda$,
observing that the holomorphic transition matrices involved in changes
of holomorphic frames do not affect the computation of~$d''$.
It is then clear that the Dolbeault-Grothendieck lemma still
holds for $F$-valued forms. For every integer $p=0,1\ld n$, the
{\it Dolbeault Cohomology} groups $H^{p,q}(X,F)$ are defined to be the
cohomology groups of the complex of global $(p,q)$ forms (graded by $q$):
$$H^{p,q}(X,F)=H^q\big(C^\infty(X,\Lambda^{p,\bu}T^\star_X\otimes F)\big).
\leqno(1.1)$$
Now, let us recall the following fundamental result from sheaf theory
(De Rham-Weil isomorphism theorem): let $(\cL^{\bu},d)$ be a resolution
of a sheaf $\cA$ by acyclic sheaves, i.e.\ a complex of sheaves
$(\cL^\bu,\delta)$ such that there is an exact sequence of sheaves
$$0\lra\cA\buildo j\over\lra\cL^0\buildo \delta^0\over\lra\cL^1\lra\cdots
\lra\cL^q\buildo \delta^q\over\lra\cL^{q+1}\lra\cdots~,$$
and $H^s(X,\cL^q)=0$ for all $q\ge 0$ and $s\ge 1$. Then there
is a fonctorial isomorphism
$$H^q\big(\Gamma(X,\cL^\bu)\big)\lra H^q(X,\cA).\leqno(1.2)$$
We apply this to the following situation: let $\cE(F)^{p,q}$ be the
sheaf of germs of $C^\infty$ sections of $\Lambda^{p,q}T^\star_X\otimes
F$, Then $(\cE(F)^{p,\bu},d'')$ is a resolution of the locally free
$\cO_X$-module $\Omega^p_X\otimes\cO(F)$ (Dolbeault-Grothendieck
lemma), and the sheaves $\cE(F)^{p,q}$ are acyclic as modules over the
soft sheaf of rings~$\cC^\infty$. Hence by (1.2) we get

\begstat{(1.3) Dolbeault Isomorphism Theorem {\rm (1953)}} For every
holomorphic vector bundle $F$ on~$X$, there is a canonical isomorphism
$$H^{p,q}(X,F)\simeq H^q(X,\Omega^p_X\otimes\cO(F)).\eqno\square$$
\endstat

If $X$ is projective algebraic and $F$ is an algebraic vector bundle,
Serre's GAGA theorem [Ser56] shows that the algebraic sheaf cohomology
group $H^q(X,\Omega^p_X\otimes\cO(F))$ computed with algebraic sections
over Zariski open sets is actually isomorphic to the analytic
cohomology group. These results are the most basic tools to attack
algebraic problems via analytic methods. Another important tool
is the theory of plurisubharmonic functions and positive
currents originated by K.~Oka and P.~Lelong in the decades 1940-1960.

\titlec{1.B.}{Plurisubharmonic Functions}
Plurisubharmonic functions have been introduced independently
by Lelong and Oka in the study of holomorphic convexity. We refer to
[Lel67, 69] for more details.

\begstat{(1.4) Definition} A function $u:\Omega\lra[-\infty,+\infty[$ defined
on an open subset $\Omega\subset\bbbc^n$ is said to be plurisubharmonic
$($psh for short$)$ if
\smallskip
\item{\rm a)} $u$ is upper semicontinuous~$;$
\smallskip
\item{\rm b)} for every complex line $L\subset\bbbc^n$,
$u_{\restriction\Omega\cap L}$ is subharmonic on $\Omega\cap L$, that is,
for all $a\in\Omega$ and $\xi\in\bbbc^n$ with $|\xi|<d(a,\complement\Omega)$,
the function $u$ satisfies the mean value inequality
$$u(a)\le{1\over 2\pi}\int_0^{2\pi}u(a+e^{\ii\theta}\,\xi)\,d\theta.$$
\vskip0pt\noindent
The set of psh functions on $\Omega$ is denoted by
$\Psh(\Omega)$.
\endstat

We list below the most basic properties of psh functions.
They all follow easily from the definition.

\begstat{(1.5) Basic properties} \rm
\smallskip
\item{\rm a)} Every function $u\in\Psh(\Omega)$ is subharmonic, namely it
satisfies the mean value inequality on euclidean balls or spheres:
$$u(a)\le{1\over\pi^nr^{2n}/n!}\int_{B(a,r)}u(z)\,d\lambda(z)$$
for every $a\in\Omega$ and $r<d(a,\complement\Omega$. Either $u\equiv-\infty$
or $u\in L^1_\loc$ on every connected component of~$\Omega$.
\smallskip
\item{\rm b)} For any decreasing sequence of psh
functions $u_k\in\Psh(\Omega)$, the limit $u=\lim u_k$ is psh
on $\Omega$.
\smallskip
\item{\rm c)} Let $u\in\Psh(\Omega)$ be such that $u\not\equiv-\infty$
on every connected component of~$\Omega$. If $(\rho_\varepsilon)$ is a
family of smoothing kernels, then $u\star\rho_\varepsilon$ is $\ci$ and
psh on
$$\Omega_\varepsilon=\big\{x\in\Omega\,;\,d(x,\complement\Omega)>\varepsilon
\big\},$$
the family $(u\star\rho_\varepsilon)$ is increasing in $\varepsilon$ and
$\lim_{\varepsilon\to 0}u\star\rho_\varepsilon=u$.
\smallskip
\item{\rm d)} Let $u_1\ld u_p\in\Psh(\Omega)$ and
$\chi:\bbbr^p\lra\bbbr$ be a convex function such that $\chi(t_1\ld t_p)$
is increasing in each $t_j$. Then $\chi(u_1\ld u_p)$ is psh
on $\Omega$. In particular~ $u_1+\cdots+u_p$, $\max\{u_1\ld u_p\}$,
$\log(e^{u_1}+\cdots+e^{u_p})$ are psh on~$\Omega$.\qed
\smallskip
\endstat

\begstat{(1.6) Lemma} A function $u\in C^2(\Omega,\bbbr)$ is psh on
$\Omega$ if and only if the hermitian form $Hu(a)(\xi)=\sum_{1\le j,k\le n}
\partial^2 u/\partial z_j\partial\ovl z_k(a)\,\xi_j\ovl \xi_k$
is semipositive at every point $a\in\Omega$.
\endstat

\begproof{} This is an easy consequence of the following standard formula
$${1\over 2\pi}\int_0^{2\pi}u(a+e^{\ii\theta}\,\xi)\,d\theta - u(a)=
{2\over\pi}\int_0^1{dt\over t}\int_{|\zeta|<t}Hu(a+\zeta\xi)(\xi)\,
d\lambda(\zeta),$$
where $d\lambda$ is the Lebesgue measure on $\bbbc$.
Lemma~1.6 is a strong evidence that plurisubharmonicity is the natural
complex analogue of linear convexity.\qed
\endproof

For non smooth functions, a similar characterization of plurisubharmonicity
can be obtained by means of a regularization process.

\begstat{(1.7) Theorem} If $u\in\Psh(\Omega)$, $u\not\equiv-\infty$ on every
connected component of $\Omega$, then for all $\xi\in\bbbc^n$
$$Hu(\xi)=\sum_{1\le j,k\le n}{\partial^2 u\over\partial z_j\partial\ovl z_k}
\,\xi_j\ovl\xi_k\in\cD'(\Omega)$$
is a positive measure. Conversely, if $v\in\cD'(\Omega)$ is such that
$Hv(\xi)$ is a positive measure for every $\xi\in\bbbc^n$, there exists a
unique function $u\in\Psh(\Omega)$ which is locally integrable on $\Omega$
and such that $v$ is the distribution associated to $u$.\qed
\endstat

In order to get a better geometric insight of this notion, we assume more
generally that $u$ is a function on a complex $n$-dimensional manifold
$X$. If $\Phi:X\to Y$ is a holomorphic mapping and if $v\in C^2(Y,\bbbr)$, we
have $\hbox{$d'd''(v\circ \Phi)$}=\Phi^\star d'd''v$, hence
$$H(v\circ \Phi)(a,\xi)=Hv\big(\Phi(a),\Phi'(a).\xi\big).$$
In particular $Hu$, viewed as a hermitian form on $T_X$, does not depend
on the choice of coordinates $(z_1\ld z_n)$. Therefore, the notion of
psh function makes sense on any complex manifold. More generally,
we have

\begstat{(1.8) Proposition} If $\Phi:X\lra Y$ is a holomorphic map and
$v\in\Psh(Y)$, then $v\circ \Phi\in\Psh(X)$.\qed
\endstat

\begstat{(1.9) Example} \rm It is a standard fact that $\log|z|$ is psh
(i.e.\ subharmonic) on $\bbbc$. Thus $\log|f|\in\Psh(X)$ for every
holomorphic function $f\in H^0(X,\cO_X)$. More generally
$$\log\big(|f_1|^{\alpha_1}+\cdots+|f_q|^{\alpha_q}\big)\in\Psh(X)$$
for every $f_j\in H^0(X,\cO_X)$ and $\alpha_j\ge 0$ (apply Property~1.5$\,$d
with $u_j=\alpha_j\,\log|f_j|$). We will be especially interested in the
singularities obtained at points of the zero variety $f_1=\ldots=f_q=0$,
when the $\alpha_j$ are rational numbers.\qed
\endstat

\begstat{(1.10) Definition} A psh function $u\in\Psh(X)$ will be said to
have analytic singularities if $u$ can be written locally as
$$u={\alpha\over 2}\log\big(|f_1|^2+\cdots+|f_N|^2\big)+v,$$
where $\alpha\in\bbbr_+$, $v$ is a locally bounded function and the $f_j$
are holomorphic functions. If $X$ is algebraic, we say that $u$ has
algebraic singularities if $u$ can be written as above on sufficiently
small Zariski open sets, with $\alpha\in\bbbq_+$ and $f_j$ algebraic.
\endstat

We then introduce the ideal $\cJ=\cJ(u/\alpha)$ of germs of holomorphic
functions $h$ such that $|h|\le Ce^{u/\alpha}$ for some constant~$C$, i.e.
$$|h|\le C\big(|f_1|+\cdots+|f_N|\big).$$
This is a globally defined ideal sheaf on $X$, locally equal to the
integral closure $\ovl\cI$ of the ideal sheaf $\cI=(f_1\ld f_N)$, thus $\cJ$
is coherent on~$X$. If $(g_1\ld g_{N'})$ are local generators of $\cJ$,
we still have
$$u={\alpha\over 2}\log\big(|g_1|^2+\cdots+|g_{N'}|^2\big)+O(1).$$
If $X$ is projective algebraic and $u$ has analytic singularities with
$\alpha\in\bbbq_+$, then $u$ automatically has algebraic singularities.
{}From an algebraic point of view, the singularities of $u$ are in 1:1
correspondence with the ``algebraic data'' $(\cJ,\alpha)$. Later
on, we will see another important method for associating an ideal sheaf
to a psh function.

\begstat{(1.11) Exercise} \rm Show that the above definition of the integral
closure of an ideal $\cI$ is equivalent to the following more algebraic
definition: $\ovl\cI$ consists of all germs $h$ satisfying
an integral equation
$$h^d+a_1h^{d-1}+\ldots+a_{d-1}h+a_d=0,\qquad a_k\in\cI^k.$$
{\it Hint}. One inclusion is clear. To prove the other inclusion, consider
the normalization of the blow-up of $X$ along the (non necessarily reduced)
zero variety $V(\cI)$.\qed
\endstat

\titlec{1.C.}{Positive Currents}
The reader can consult [Fed69] for a more thorough treatment of
current theory. Let us first recall a few basic definitions. A {\it
current} of degree $q$ on an oriented differentiable manifold $M$ is
simply a differential $q$-form $\Theta$ with distribution coefficients.
The space of currents of degree $q$ over $M$ will be denoted by
$\cD^{\prime q}(M)$. Alternatively, a current of degree $q$ can be seen
as an element $\Theta$ in the dual space $\cD'_p(M):=\big(\cD^p(M)\big)'$ of
the space $\cD^p(M)$ of smooth differential forms of degree $p=\dim M-q$ with
compact support; the duality pairing is given by
$$\langle \Theta,\alpha\rangle=\int_M \Theta\wedge\alpha,~~~~
\alpha\in\cD^p(M).\leqno(1.12)$$
A basic example is the {\it current of integration} $[S]$ over a compact
oriented submanifold $S$ of $M\,$:
$$\langle[S],\alpha\rangle=\int_S\alpha,~~~~{\rm deg}\alpha=p=\dim_\bbbr S.
\leqno(1.13)$$
Then $[S]$ is a current with measure coefficients, and
Stokes' formula shows that $d[S]=(-1)^{q-1}[\partial S]$, in particular
$d[S]=0$ if $S$ has no boundary. Because of this example, the integer $p$
is said to be the dimension of~$\Theta$ when $\Theta\in\cD_p'(M)$. The
current $\Theta$ is said to be {\it closed} if $d\Theta=0$.

On a complex manifold $X$, we have similar notions of bidegree and
bidimension; as in the real case, we denote by
$$\cD^{\prime p,q}(X)=\cD'_{n-p,n-q}(X),\qquad n=\dim X,$$
the space of currents of bidegree $(p,q)$ and bidimension $(n-p,n-q)$ on~$X$.
According to [Lel57], a current $\Theta$
of bidimension $(p,p)$ is said to be {\it $($weakly$)$ positive} if
for every choice of smooth $(1,0)$-forms $\alpha_1\ld\alpha_p$ on~$X$
the distribution
$$\Theta\wedge\ii\alpha_1\wedge\ovl\alpha_1\wedge\ldots\wedge
\ii\alpha_p\wedge\ovl\alpha_p~~~~\hbox{\rm is a positive measure.}
\leqno(1.14)$$

\begstat{(1.15) Exercise} \rm If $\Theta$ is positive, show that the
coefficients $\Theta_{I,J}$ of $\Theta$ are complex measures, and that,
up to constants, they are dominated by the trace measure
$$\sigma_\Theta=\Theta\wedge{1\over p!}\beta^p=2^{-p}\sum\Theta_{I,I},\qquad
\beta={\ii\over 2}d'd''|z|^2={\ii\over 2}\sum_{1\le j\le n}dz_j\wedge
d\ovl z_j,$$
which is a positive measure.\newline
{\it Hint}. Observe that $\sum\Theta_{I,I}$ is invariant by unitary changes
of coordinates and that the $(p,p)$-forms
$\ii\alpha_1\wedge\ovl\alpha_1\wedge\ldots\wedge\ii\alpha_p\wedge\ovl\alpha_p$
generate $\Lambda^{p,p}T^\star_{\bbbc^n}$ as a $\bbbc$-vector space.\qed
\endstat

A current $\Theta=\ii\sum_{1\le j,k\le n}\Theta_{jk}dz_j\wedge dz_k$
of bidegree $(1,1)$ is easily seen to be positive if and only if the
complex measure $\sum\lambda_j\ovl\lambda_k \Theta_{jk}$ is a positive
measure for every $n$-tuple $(\lambda_1\ld\lambda_n)\in\bbbc^n$.

\begstat{(1.16) Example} \rm If $u$ is a (not identically $-\infty$)
psh function on $X$, we can associate with $u$ a (closed) positive
current $\Theta=\ii\partial\ovl\partial u$ of bidegree~$(1,1)$. Conversely,
every closed positive current of bidegree $(1,1)$ can be written under
this form on any open subset $\Omega\subset X$ such that
$H^2_{DR}(\Omega,\bbbr)=H^1(\Omega,\cO)=0$, e.g.\ on small coordinate
balls (exercise to the reader).\qed
\endstat

It is not difficult to show that a product $\Theta_1\wedge\ldots\wedge
\Theta_q$ of positive currents of bidegree $(1,1)$ is positive whenever
the product is well defined (this is certainly the case if all
$\Theta_j$ but one at most are smooth; much finer conditions will be
discussed in Section~2).

We now discuss another very important example of closed positive current.
In fact, with every closed analytic set $A\subset X$ of pure dimension $p$
is associated a current of integration
$$\langle [A],\alpha\rangle=\int_{A_\reg}\alpha,~~~~\alpha\in\cD^{p,p}(X),
\leqno(1.17)$$
obtained by integrating over the regular points of~$A$. In order to show
that (1.17) is a correct definition of a current on~$X$, one must show that
$A_\reg$ has locally finite area in a neighborhood of $A_\sing$.
This result, due to [Lel57] is shown as follows. Suppose that $0$
is a singular point of $A$. By the local parametrization theorem for
analytic sets, there is a linear change of coordinates on $\bbbc^n$
such that all projections
$$\pi_I:(z_1\ld z_n)\mapsto (z_{i_1}\ld z_{i_p})$$
define a finite ramified covering of the intersection $A\cap\Delta$ with
a small polydisk $\Delta$ in $\bbbc^n$ onto a small
polydisk $\Delta_I$ in~$\bbbc^p$. Let $n_I$ be the sheet number. Then the
\hbox{$p$-dimensional} area of $A\cap\Delta$ is bounded above by the sum
of the areas of its projections counted with multiplicities, i.e.
$$\hbox{\rm Area}(A\cap\Delta)\le\sum n_I\hbox{\rm Vol}(\Delta_I).$$
The fact that $[A]$ is positive is also easy. In fact
$$\ii\alpha_1\wedge\ovl\alpha_1\wedge\ldots\wedge
\ii\alpha_p\wedge\ovl\alpha_p=|\det(\alpha_{jk})|^2\,
\ii w_1\wedge\ovl w_1\wedge\ldots\wedge
\ii w_p\wedge\ovl w_p$$
if $\alpha_j=\sum\alpha_{jk}dw_k$ in terms of local coordinates
$(w_1\ld w_p)$ on $A_\reg$. This shows that all such forms are $\ge 0$
in the canonical orientation defined by
$\ii w_1\wedge\ovl w_1\wedge\ldots\wedge\ii w_p\wedge\ovl w_p$.
More importantly, Lelong [Lel57] has shown that $[A]$ is $d$-closed in~$X$,
even at points of $A_\sing$. This last result can be seen today as a
consequence of the Skoda-El Mir extension theorem. For this we need the
following definition: a {\it complete pluripolar} set is a set $E$ such
that there is an open covering $(\Omega_j)$ of $X$ and psh functions $u_j$ on
$\Omega_j$ with $E\cap\Omega_j=u_j^{-1}(-\infty)$. Any (closed) analytic
set is of course complete pluripolar (take $u_j$ as in Example~1.9).

\begstat{(1.18) Theorem {\rm(Skoda [Sko81], El Mir [EM84], Sibony [Sib85])}}
Let $E$ be a closed complete pluripolar set in $X$,
and let $\Theta$ be a closed positive current on $X\ssm E$ such that the
coefficients $\Theta_{I,J}$ of $\Theta$ are measures with locally finite
mass near~$E$. Then the trivial extension $\wt \Theta$ obtained by extending
the measures $\Theta_{I,J}$ by $0$ on $E$ is still closed on~$X$.
\endstat

Lelong's result $d[A]=0$ is obtained by applying the Skoda-El Mir
theorem to $\Theta=[A_\reg]$ on $X\ssm A_\sing$.

\begproof{of Theorem 1.18.} The statement is local on $X$, so we may
work on a small open set $\Omega$ such that $E\cap\Omega=v^{-1}(-\infty)$,
$v\in\Psh(\Omega)$. Let $\chi:\bbbr\to\bbbr$ be a convex increasing function
such that $\chi(t)=0$ for $t\le -1$ and $\chi(0)=1$. By shrinking $\Omega$
and putting $v_k=\chi(k^{-1}v\star\rho_{\varepsilon_k})$ with
$\varepsilon_k\to 0$ fast, we get a sequence of functions
$v_k\in\Psh(\Omega)\cap\ci(\Omega)$ such that $0\le v_k\le 1$, $v_k=0$ in
a neighborhood of $E\cap\Omega$ and $\lim v_k(x)=1$ at every point of
$\Omega\ssm E$. Let $\theta\in\ci([0,1])$ be a
function such that $\theta=0$ on $[0,1/3]$, $\theta=1$ on $[2/3,1]$ and
$0\le\theta\le 1$. Then $\theta\circ v_k=0$ near $E\cap\Omega$ and
$\theta\circ v_k\to 1$ on $\Omega\ssm E$. Therefore
$\tilde\Theta=\lim_{k\to+\infty}(\theta\circ v_k)\Theta$ and
$$d'\tilde\Theta=\lim_{k\to+\infty}\Theta\wedge d'(\theta\circ v_k)$$
in the weak topology of currents. It is therefore sufficient to verify that
\hbox{$\Theta\wedge d'(\theta\circ v_k)$} converges weakly to $0$
(note that $d''\tilde \Theta$ is conjugate to $d'\tilde \Theta$,
thus $d''\tilde \Theta$ will also vanish).

Assume first that $\Theta\in\cD^{\prime n-1,n-1}(X)$. Then
$\Theta\wedge d'(\theta\circ v_k)\in\cD^{\prime n,n-1}(\Omega)$, and we have
to show that
$$\langle \Theta\wedge d'(\theta\circ v_k),\ovl\alpha\rangle=\langle\Theta,
\theta'(v_k)d'v_k\wedge\ovl\alpha\rangle
\buildu{k\to+\infty}\under\lra 0,~~~~\forall\alpha\in\cD^{1,0}(\Omega).$$
As $\gamma\mapsto\langle\Theta,\ii\gamma\wedge\ovl\gamma\rangle$ is a
non-negative hermitian form on $\cD^{1,0}(\Omega)$, the Cauchy-Schwarz
inequality yields
$$\big|\langle \Theta,\ii\beta\wedge\ovl\gamma\rangle\big|^2\le
\langle \Theta,\ii\beta\wedge\ovl\beta\rangle~
\langle \Theta,\ii\gamma\wedge\ovl\gamma\rangle,~~~~
\forall\beta,\gamma\in\cD^{1,0}(\Omega).$$
Let $\psi\in\cD(\Omega)$, $0\le\psi\le 1$, be equal to $1$ in a neighborhood
of ${\rm Supp}\,\alpha$. We find
$$\big|\langle \Theta,\theta'(v_k)d'v_k\wedge\ovl\alpha\rangle\big|^2\le
\langle \Theta,\psi \ii d'v_k\wedge d''v_k\rangle~
\langle \Theta,\theta'(v_k)^2 \ii\alpha\wedge\ovl\alpha\rangle.$$
By hypothesis $\int_{\Omega\ssm E}\Theta\wedge \ii\alpha\wedge\ovl\alpha<
+\infty$ and $\theta'(v_k)$ converges everywhere to $0$ on $\Omega$, thus
$\langle \Theta,\theta'(v_k)^2 \ii\alpha\wedge\ovl\alpha\rangle$ converges
to $0$ by Lebesgue's dominated convergence theorem. On the other hand
$$\eqalign{
&\ii d'd''v_k^2=2v_k\,\ii d'd''v_k+2\ii d'v_k\wedge d''v_k\ge 2\ii d'v_k
\wedge d''v_k,\cr
&2\langle\Theta,\psi \ii d'v_k\wedge d''v_k\rangle\le\langle\Theta,\psi
\ii d'd''v_k^2\rangle.\cr}$$
As $\psi\in\cD(\Omega)$, $v_k=0$ near $E$ and $d\Theta=0$ on $\Omega\ssm E$,
an integration by parts yields
$$\langle \Theta,\psi \ii d'd''v_k^2\rangle=\langle \Theta,v_k^2
\ii d'd''\psi\rangle\le C\int_{\Omega\ssm E}\|\Theta\|<+\infty$$
where $C$ is a bound for the coefficients of $\ii d'd''\psi$. Thus
$\langle \Theta,\psi \ii d'v_k\wedge d''v_k\rangle$ is bounded,
and the proof is complete when $\Theta\in\cD^{\prime n-1,n-1}$.

In the general case $\Theta\in\cD^{\prime p,p}$, $p<n$, we simply apply the
result already proved to all positive currents $\Theta\wedge\gamma\in
\cD^{\prime n-1,n-1}$ where
$\gamma=\ii\gamma_1\wedge\ovl\gamma_1\wedge\ldots\wedge\ii\gamma_{n-p-1,}
\wedge\ovl\gamma_{n-p-1}$ runs over a basis of forms of
$\Lambda^{n-p-1,n-p-1}T^\star_\Omega$ with constant coefficients (Lemma~1.4).
Then we get $d(\tilde \Theta\wedge\gamma)=d\tilde \Theta\wedge\gamma=0$ for
all such $\gamma$, hence $d\tilde \Theta=0$.\qed
\endproof

\begstat{(1.19) Corollary} Let $\Theta$ be a closed positive current on
$X$ and let $E$ be a complete pluripolar set. Then $\bbbone_E \Theta$
and $\bbbone_{X\ssm E}\Theta$ are closed positive currents. In fact,
\hbox{$\wt\Theta=\bbbone_{X\ssm E}\Theta$} is the trivial extension of
$\Theta_{\restriction X\ssm E}$ to~$X$, and
$\bbbone_E\Theta=\Theta-\wt\Theta$.\qed
\endstat

As mentioned above, any current $\Theta=\ii d'd''u$ associated with a psh
function $u$ is a closed positive $(1,1)$-current. In the special
case $u=\log|f|$ where $f\in H^0(X,\cO_X)$ is a non zero holomorphic
function, we have the important

\vbox{\begstat{(1.20) Lelong-Poincar\'e equation}
Let $f\in H^0(X,\cO_X)$ be a non zero holomorphic function,
$Z_f=\sum m_jZ_j$, $m_j\in\bbbn$, the zero divisor of~$f$ and
$[Z_f]=\sum m_j[Z_j]$ the associated current of integration. Then
$${\ii\over\pi}\partial\ovl\partial\log|f|=[Z_f].$$
\endstat\removelastskip}

\begproof{(sketch).} It is clear that $\ii d'd''\log|f|=0$ in a neighborhood
of every point $x\notin\Supp(Z_f)=\bigcup Z_j$, so it is enough to check
the equation in a neighborhood of every point of $\Supp(Z_f)$. Let $A$ be
the set of singular points of $\Supp(Z_f)$, i.e.\ the union of the
pairwise intersections $Z_j\cap Z_k$ and of the singular loci
$Z_{j,\sing}$; we thus have $\dim A\le n-2$. In a neighborhood of any
point $x\in\Supp(Z_f)\ssm A$ there are local coordinates $(z_1\ld z_n)$
such that $f(z)=z_1^{m_j}$ where $m_j$ is the multiplicity of $f$ along
the component $Z_j$ which contains~$x$ and $z_1=0$ is an equation for
$Z_j$ near~$x$. Hence
$${\ii\over\pi}d'd''\log|f|=m_j{\ii\over\pi}d'd''\log|z_1|=m_j[Z_j]$$
in a neighborhood of $x$, as desired (the identity comes from the
standard formula ${\ii\over\pi}d'd''\log|z|={}$ Dirac measure $\delta_0$
in $\bbbc$). This shows that the equation holds on $X\ssm A$. Hence
the difference ${\ii\over\pi}d'd''\log|f|-[Z_f]$ is a closed current of
degree $2$ with measure coefficients, whose support is contained in~$A$.
By Exercise~1.21, this current must be $0$, for $A$ has too small
dimension to carry its support ($A$ is stratified by submanifolds of real
codimension $\ge 4$).\qed
\endproof

\begstat{(1.21) Exercise} \rm Let $\Theta$ be a current of degree $q$
on a real manifold~$M$, such that both $\Theta$ and $d\Theta$ have
measure coefficients (``normal current''). Suppose that $\Supp \Theta$
is contained in a real submanifold $A$ with $\codim_\bbbr A>q$. Show
that $\Theta=0$.\newline
{\it Hint:} Let $m=\dim_\bbbr M$ and let $(x_1\ld x_m)$ be a coordinate
system in a neighborhood $\Omega$ of a point $a\in A$ such that
$A\cap\Omega=\{x_1=\ldots=x_k=0\}$, $k>q$. Observe that
$x_j\Theta=x_jd\Theta=0$ for $1\le j\le k$, thanks to the hypothesis on
supports and on the normality of~$\Theta$, hence $dx_j\wedge
\Theta=d(x_j\Theta)-x_jd\Theta=0$, $1\le j\le k$. Infer from this that
all coefficients in $\Theta=\sum_{|I|=q}\Theta_Idx_I$ vanish.\qed
\endstat

We now recall a few basic facts of slicing theory (the reader will
profitably consult [Fed69] and [Siu74] for further developments).
Let \hbox{$\sigma:M\to M'$} be a submersion of smooth differentiable
manifolds and let $\Theta$ be a {\it locally flat} current on~$M$,
that is, a current which can be written locally as $\Theta=U+dV$
where $U$, $V$ have $L^1_\loc$ coefficients. It is a standard fact
(see Federer) that every current $\Theta$ such that both $\Theta$ and
$d\Theta$ have measure coefficients is locally flat; in particular,
closed positive currents are locally flat. Then, for almost every
$x'\in M'$, there is a well defined {\it slice} $\Theta_{x'}$,
which is the current on the fiber $\sigma^{-1}(x')$ defined by
$$\Theta_{x'}=U_{\restriction\sigma^{-1}(x')}+
  dV_{\restriction\sigma^{-1}(x')}.$$
The restrictions of $U$, $V$ to the fibers exist for almost all $x'$ by
the Fubini theorem. The slices $\Theta_{x'}$ are currents on the fibers with
the same degree as~$\Theta$ (thus of dimension $\dim \Theta-\dim
\hbox{(fibers)}$). Of course, every slice $\Theta_{x'}$ coincides with
the usual restriction of $\Theta$ to the fiber if $\Theta$ has smooth
coefficients. By using a regularization $\Theta_\varepsilon=\Theta\star
\rho_\varepsilon$, it is easy to show that the slices of a
closed positive current are again closed and positive: in fact
$U_{\varepsilon,x'}$ and $V_{\varepsilon,x'}$ converge to $U_{x'}$ and
$V_{x'}$ in $L^1_{\rm loc}(\sigma^{-1}(x'))$, thus
$\Theta_{\varepsilon,x'}$ converges weakly to $\Theta_{x'}$ for almost
every~$x'$. Now, the basic slicing formula is
$$\int_M \Theta\wedge\alpha\wedge\sigma^\star\beta=
\int_{x'\in M'}\Big(\int_{x''\in\sigma^{-1}(x')}\Theta_{x'}(x'')\wedge
\alpha_{\restriction\sigma^{-1}(x')}(x'')\Big)\beta(x')\leqno(1.22)$$
for every smooth form $\alpha$ on $M$ and $\beta$ on $M'$, such that
$\alpha$ has compact support and $\deg\alpha=\dim M-dim M'-\deg \Theta$,
$\deg\beta=\dim M'$. This is an easy consequence of the usual Fubini
theorem applied to $U$ and $V$ in the decomposition $\Theta=U+dV$, if we
identify locally $\sigma$ with a projection map
$M=M'\times M''\to M'$, $x=(x',x'')\mapsto x'$, and use a partitition
of unity on the support of $\alpha$.

To conclude this section, we discuss De Rham and Dolbeault cohomology theory
in the context of currents. A basic observation is that the Poincar\'e and
Dolbeault-Grothendieck lemma still hold for currents. Namely, if
$(\cD^{\prime q},d)$ and $(\cD'(F)^{p,q},d'')$ denote the complex of sheaves
of degree $q$ currents (resp.\ of $(p,q)$-currents with values in a
holomorphic vector bundle~$F$), we still have De Rham and Dolbeault sheaf
resolutions
$$0\to\bbbr\to\cD^{\prime \bu},\qquad 0\to\Omega^p_X\otimes\cO(F)\to
\cD'(F)^{p,\bu}.$$
Hence we get canonical isomorphisms
$$\leqalignno{
H^q_\DR(M,\bbbr)&=H^q\big((\Gamma(M,\cD^{\prime\bu}),d)\big),&(1.23)\cr
H^{p,q}(X,F)&=H^q\big((\Gamma(X,\cD'(F)^{p,\bu}),d'')\big).&\cr}$$
In other words, we can attach a cohomology class
$\{\Theta\}\in H^q_\DR(M,\bbbr)$ to any closed current $\Theta$ of degree~$q$,
resp.\ a cohomology class $\{\Theta\}\in H^{p,q}(X,F)$ to any $d''$-closed
current of bidegree~$(p,q)$. Replacing if necessary every current by a
smooth representative in the same cohomology class, we see that there is a
well defined cup product given by the wedge product of differential forms
$$\eqalign{
H^{q_1}(M,\bbbr)\times\ldots\times H^{q_m}(M,\bbbr)&\lra
H^{q_1+\ldots+q_m}(M,\bbbr),\cr
(\{\Theta_1\}\ld\{\Theta_1\})&\longmapsto\{\Theta_1\}\wedge\ldots\wedge
\{\Theta_m\}.\cr}$$
In particular, if $M$ is a compact oriented variety and
$q_1+\ldots+q_m=\dim M$, there is a well defined intersection number
$$\{\Theta_1\}\cdot\{\Theta_2\}\cdot\,\cdots\,\cdot
\{\Theta_m\}=\int_M\{\Theta_1\}\wedge\ldots\wedge\{\Theta_m\}.$$
However, as we will see in the next section, the pointwise product
\hbox{$\Theta_1\wedge\ldots\wedge\Theta_m$} need not exist in general.
\vfill\eject

\titleb{2.}{Lelong Numbers and Intersection Theory}
\titlec{2.A.}{Multiplication of Currents and Monge-Amp\`ere Operators}
Let $X$ be a $n$-dimensional complex manifold. We set
$$d^c={1\over2\ii\pi}(d'-d'').$$
It follows in particular that $d^c$ is a real operator, i.e. $\ovl{d^cu}=
d^c\ovl u$, and that $dd^c={\ii\over\pi}d'd''$. Although not quite standard,
the $1/2\ii\pi$ normalization is very convenient for many purposes, since
we may then forget the factor $\pi$ or $2\pi$ almost everywhere (e.g.\ in
the Lelong-Poincar\'e equation (1.20)).

Let $u$ be a psh function and let $\Theta$ be a closed
positive current on~$X$. Our desire is to define the wedge product
$dd^cu\wedge\Theta$ even when neither $u$ nor $\Theta$ are smooth.
In general, this product does not make sense because $dd^cu$ and $\Theta$
have measure coefficients and measures cannot be multiplied; see
Kiselman [Kis84] for interesting counterexamples. Even in the algebraic
setting considered here, multiplication of currents is not always possible:
suppose e.g.\ that $\Theta=[D]$ is the exceptional divisor of a blow-up
in a surface; then $D\cdot D=-1$ cannot be the cohomology class of a
closed positive current $[D]^2$. Assume however that
$u$ is a {\it locally bounded} psh function. Then the
current $u\Theta$ is well defined since $u$ is a locally bounded Borel
function and $\Theta$ has measure coefficients. According to
Bedford-Taylor [BT82] we define
$$dd^cu\wedge\Theta=dd^c(u\Theta)$$
where $dd^c(~~)$ is taken in the sense of distribution theory.

\begstat{(2.1) Proposition} If $u$ is a locally bounded psh function, the
wedge product $dd^cu\wedge\Theta$ is again a closed positive current.
\endstat

\begproof{} The result is local. Use a convolution $u_\nu=u\star\rho_{1/\nu}$
to get a decreasing sequence of smooth psh functions converging to~$u$.
Then write
$$dd^c(u\Theta)=\lim_{\nu\to+\infty}dd^c(u_\nu\Theta)=
dd^cu_\nu\wedge\Theta$$
as a weak limit of closed positive currents. Observe that $u_\nu\Theta$
converges weakly to $u\Theta$ by Lebesgue's monotone convergence theorem.\qed
\endproof

More generally, if $u_1\ld u_m$ are locally bounded
psh functions, we can define
$$dd^cu_1\wedge\ldots\wedge dd^cu_m\wedge\Theta=
dd^c\big(u_1dd^cu_2\wedge\ldots\wedge dd^cu_m\wedge\Theta\big)$$
by induction on~$m$. (Chern, Levine and Nirenberg, 1969) noticed
the following useful inequality. Define the {\it mass} of a current
$\Theta$ on a compact set $K$ to be
$$||\Theta||_K=\int_K\sum_{I,J}|\Theta_{I,J}|$$
whenever $K$ is contained in a coordinate patch and
$\Theta=\sum\Theta_{I,J}dz_I\wedge d\ovl z_J$. Up to seminorm equivalence,
this does not depend on the choice of coordinates. If $K$ is not contained
in a coordinate patch, we use a partition of unity to define a suitable
seminorm $||\Theta||_K$. If $\Theta\ge 0$, Exercise~1.15 shows that the mass
is controlled by the trace measure, i.e.\
$||\Theta||_K\le C\int_K\Theta\wedge\beta^p$.

\begstat{(2.2) Chern-Levine-Nirenberg inequality}
For all compact subsets $K,L$ of $X$ with $L\subset K^\circ$, there exists
a constant $C_{K,L}\ge0$ such that
$$||dd^cu_1\wedge\ldots\wedge dd^cu_m\wedge\Theta||_L\le C_{K,L}~
||u_1||_{L^\infty(K)}\ldots||u_m||_{L^\infty(K)}\,||\Theta||_K$$
\endstat

\begproof{} By induction, it is sufficient to prove the result for
$m=1$ and $u_1=u$. There is a covering of $L$ by a family of open balls
$B'_j\compact B_j\subset K$ contained in coordinate patches of~$X$. Let
$(p,p)$ be the bidimension of $\Theta$, let $\beta={\ii\over 2}d'd''|z|^2$,
and let $\chi\in\cD(B_j)$ be equal to $1$ on $\smash{\ovl B'_j}$. Then
$$||dd^cu\wedge\Theta||_{L\cap\ovl B'_j}\le C\int_{\ovl B'_j}dd^cu\wedge
\Theta\wedge\beta^{p-1}\le C\int_{B_j}\chi\,dd^c u\wedge\Theta\wedge
\beta^{p-1}.$$
As $\Theta$ and $\beta$ are closed, an integration by parts yields
$$||dd^cu\wedge\Theta||_{L\cap\ovl B'_j}\le C\int_{B_j}u\,\Theta\wedge
dd^c\chi\wedge\beta^{p-1}\le C'||u||_{L^\infty(K)}||\Theta||_K$$
where $C'$ is equal to $C$ multiplied by a bound for the coefficients
of the smooth form $dd^c\chi\wedge\beta^{p-1}$.\qed
\endproof

Various examples (cf.\ [Kis84]) show however that products of
$(1,1)$-currents $dd^cu_j$ cannot be defined in a reasonable way for
arbitrary psh functions~$u_j$. However, functions $u_j$
with $-\infty$ poles can be admitted if the polar sets are sufficiently
small.

\begstat{(2.3) Proposition} Let $u$ be a psh function
on $X$, and let $\Theta$ be a closed positive current
of bidimension $(p,p)$. Suppose that $u$ is locally bounded on $X\ssm A$,
where $A$ is an analytic subset of $X$ of dimension $<p$ at each point.
Then $dd^cu\wedge\Theta$ can be defined in such a way that
$dd^cu\wedge\Theta=\lim_{\nu\to+\infty}dd^cu_\nu\wedge\Theta$
in the weak topology of currents, for any decreasing sequence
$(u_\nu)_{\nu\ge 0}$ of psh functions converging to~$u$.
\endstat

\begproof{} When $u$ is locally bounded everywhere, we have
$\,\lim u_\nu\,\Theta=u\,\Theta$ by the monotone convergence
theorem and the result follows from the continuity of $dd^c$ with
respect to the weak topology.

First assume that $A$ is discrete. Since our results are local, we
may suppose that $X$ is a ball $B(0,R)\subset\bbbc^n$ and that $A=\{0\}$.
For every $s\le 0$, the function $u^{\gge s}=\max(u,s)$
is locally bounded on $X$, so the product $\Theta\wedge
dd^cu^{\gge s}$ is well defined. For $|s|$ large, the function
$\smash{u^{\gge s}}$ differs from $u$ only in a small neighborhood
of the origin, at which $u$ may have a $-\infty$ pole. Let $\gamma$
be a $(p-1,p-1)$-form with constant coefficients and set
$s(r)=\liminf_{|z|\to r-0}u(z)$. By Stokes' formula, we see that
the integral
$$I(s):=\int_{B(0,r)}dd^cu^{\gge s}\wedge\Theta\wedge\gamma\leqno(2.4)$$
does not depend on $s$ when $s<s(r)$, for the difference $I(s)-I(s')$
of two such integrals involves the $dd^c$ of a current
$(u^{\gge s}-u^{\gge s'})\wedge\Theta\wedge\gamma$ with compact support
in $B(0,r)$. Taking $\gamma=(dd^c|z|^2)^{p-1}$, we see that the current
$dd^cu\wedge\Theta$ has finite mass on $B(0,r)\ssm\{0\}$
and we can define $\langle\bbbone_{\{0\}}(dd^cu\wedge\Theta),
\gamma\rangle$ to be the limit of the integrals $(2.4)$ as $r$ tends to
zero and $s<s(r)$. In this case, the weak convergence statement is easily
deduced from the locally bounded case discussed above.

In the case where $0<\dim A<p$, we use a slicing technique to reduce the
situation to the discrete case. Set $q=p-1$. There are linear
coordinates $(z_1\ld z_n)$ centered at any point of $A$,
such that $0$ is an isolated point of $A\cap\big(\{0\}\times\bbbc^{n-q}\big)$.
Then there are small balls $B'=B(0,r')$ in $\bbbc^q$, $B''=B(0,r'')$ in
$\bbbc^{n-q}$ such that $A\cap(\smash{\ovl B}'\times\partial B'')=
\emptyset$, and the projection map
$$\pi:\bbbc^n\to\bbbc^q,~~~~z=(z_1\ld z_n)\mapsto z'=(z_1\ld z_q)$$
defines a finite proper mapping $A\cap(B'\times B'')\to B'$.
These properties are preserved if we slightly change the direction of
projection. Take sufficiently many projections $\pi_m$ associated
to coordinate systems $(z^m_1\ld z^m_n)$, \hbox{$1\le m\le N$}, in such
a way that the family of $(q,q)$-forms
$$i\,dz^m_1\wedge d\ovl z^m_1\wedge\ldots\wedge i\,dz^m_q\wedge d\ovl z^m_q$$
defines a basis of the space of $(q,q)$-forms. Expressing any
compactly supported smooth $(q,q)$-form in such a basis, we see
that we need only define
$$\leqalignno{
&\int_{B'\times B''}dd^cu\wedge\Theta\wedge f(z',z'')\,
\ii\,dz_1\wedge d\ovl z_1\wedge\ldots\wedge\ii\,dz_q\wedge d\ovl z_q=&(2.5)\cr
&\int_{B'}\Big\{\int_{B''}f(z',\bu)\,
dd^cu(z',\bu)\wedge\Theta(z',\bu)\Big\}
\ii\,dz_1\wedge d\ovl z_1\wedge\ldots\wedge\ii\,dz_q\wedge d\ovl z_q\cr}$$
where $f$ is a test function with compact support in $B'\times B''$,
and $\Theta(z',\bu)$ denotes the slice of $\Theta$ on the fiber
$\{z'\}\times B''$ of the projection \hbox{$\pi:\bbbc^n\to\bbbc^q$}.
Each integral $\int_{B''}$ in the right hand side of (2.5) makes sense
since the slices $(\{z'\}\times B'')\cap A$ are discrete.
Moreover, the double integral $\int_{B'}\int_{B''}$ is convergent.
Indeed, observe that $u$ is bounded on any compact cylinder
$$K_{\delta,\varepsilon}=\smash{\ovl B}\big((1-\delta)r'\big)\times
\Big(\smash{\ovl B}(r'')\ssm\smash{\ovl B}\big((1-\varepsilon)r''\big)\Big)$$
disjoint from $A$. Take $\varepsilon\ll\delta\ll 1$ so small that
$${\rm Supp}\,f\subset\smash{\ovl B}\big((1-\delta)r'\big)\times
\smash{\ovl B}\big((1-\varepsilon)r''\big).$$
For all $z'\in\smash{\ovl B}((1-\delta)r')$, the proof of the
Chern-Levine-Nirenberg inequality 2.2 with a cut-off function $\chi(z'')$
equal to $1$ on $B((1-\varepsilon)r'')$ and with support in
$B((1-\varepsilon/2)r'')$ shows that
$$\eqalign{
&\int_{B((1-\varepsilon)r'')}dd^cu(z',\bu)\wedge\Theta(z',\bu)\cr
&\qquad\quad{}\le C_\varepsilon ||u||_{L^\infty(K_{\delta,\varepsilon})}
\int_{z''\in B((1-\varepsilon/2)r'')}
\Theta(z',z'')\wedge dd^c|z''|^2.\cr}$$
This implies that the double integral is convergent. Now replace $u$
everywhere by $u_\nu$ and observe that $\lim_{\nu\to+\infty}\int_{B''}$
is the expected integral for every $z'$ such that $\Theta(z',\bu)$ exists
(apply the discrete case already proven). Moreover, the
Chern-Levine-Nirenberg inequality yields uniform bounds for all
functions $u_\nu$, hence Lebesgue's dominated convergence theorem can
be applied to $\int_{B'}$. We conclude from this that the sequence of
integrals (2.5) converges when \hbox{$u_\nu\downarrow u$}, as
expected.\qed
\endproof

\begstat{(2.6) Remark} \rm In the above proof, the fact that $A$ is an
analytic set does not play an essential role. The main point is just
that the slices $(\{z'\}\times B'')\cap A$ consist of isolated points
for generic choices of coordinates $(z',z'')$. In fact, the proof even
works if the slices are totally discontinuous, in particular if they are
of zero Hausdorff measure~$\cH_1$. It follows that Proposition~2.3 still
holds whenever $A$ is a closed set such that $\cH_{2p-1}(A)=0$.\qed
\endstat

\titlec{2.B.}{Lelong Numbers}
The concept of Lelong number is an analytic analogue of the algebraic
notion of multiplicity. It is a very useful technique to extend
results in the intersection theory of algebraic cycles to currents.
Lelong numbers have been introduced for the first time by Lelong in [Lel57].
See also [Lel69], [Siu74], [Dem82a, 85a, 87] for further developments.

Let us first recall a few definitions.
Let $\Theta$ be a closed positive current of bidimension $(p,p)$
on a coordinate open set $\Omega\subset\bbbc^n$ of a complex manifold~$X$.
The Lelong number of $\Theta$ at a point $x\in\Omega$ is defined to be
the limit
$$\nu(\Theta,x)=\lim_{r\to 0{\scriptscriptstyle+}}\nu(\Theta,x,r),\qquad
\hbox{where}~~
\nu(\Theta,x,r)={\sigma_\Theta(B(x,r))\over\pi^p r^{2p}/p!}$$
measures the ratio of the area of $\Theta$ in the ball $B(x,r)$ to the
area of the ball of radius $r$ in $\bbbc^p$. As $\sigma_\Theta=\Theta\wedge
{1\over p!}(\pi dd^c|z|^2)^p$ by 1.15, we also get
$$\nu(\Theta,x,r)={1\over r^{2p}}\int_{B(x,r)}\Theta(z)\wedge
(dd^c|z|^2)^p.\leqno(2.7)$$
The main results concerning Lelong numbers are summarized in the following
theorems, due respectively to Lelong, Thie and Siu.

\begstat{(2.8) Theorem {\rm([Lel57])}}
\smallskip
\item{\rm a)}  For every positive current $\Theta$, the
ratio $\nu(\Theta,x,r)$ is a nonnegative increasing function of $r$, in
particular the limit $\nu(\Theta,x)$ as $r\to0{\scriptstyle+}$
always exists.
\smallskip
\item{\rm b)} If $\Theta=dd^cu$ is the bidegree $(1,1)$-current associated
with a psh function~$u$, then
$$\nu(\Theta,x)=\sup\big\{\gamma\ge 0\,;\,u(z)\le\gamma\log|z-x|+O(1)
{}~~\hbox{at $x$}\big\}.$$
In particular, if $u=\log|f|$ with $f\in H^0(X,\cO_X)$ and
$\Theta=dd^cu=[Z_f]$, we have
$$\nu([Z_f],x)=\ord_x(f)=\max\{m\in\bbbn\,;\,D^\alpha f(x)=0,\,|\alpha|<m\}.$$
\vskip0pt
\endstat

\begstat{(2.9) Theorem {\rm([Thi67])}} In the case where $\Theta$
is a current of integration $[A]$ over an analytic subvariety $A$, the
Lelong number $\nu([A],x)$ coincides with the multiplicity of $A$ at $x$
$($defined e.g.\ as the sheet number in the ramified covering obtained by
taking a generic linear projection of the germ $(A,x)$ onto a
$p$-dimensional linear subspace through $x$ in any coordinate
patch~$\Omega)$.
\endstat

\begstat{(2.10) Theorem {\rm([Siu74])}} Let $\Theta$ be a closed positive
current of bidimension $(p,p)$ on the complex manifold~$X$.
\smallskip
\item{\rm a)} The Lelong number $\nu(\Theta,x)$ is invariant by holomorphic
changes of local coordinates.
\smallskip
\item{\rm b)} For every $c>0$, the set
$E_c(\Theta)=\big\{x\in X\,;\,\nu(\Theta,x)\ge c\big\}$
is a closed analytic subset of~$X$ of dimension $\le p$.\vskip0pt
\endstat

The most important result is 2.10~b), which is a deep application of
H\"ormander $L^2$ estimates (see Section~5). The earlier proofs of all
other results were rather intricate in spite of their rather simple nature.
We reproduce below a sketch of elementary arguments based on the use of
a more general and more flexible notion of Lelong number introduced in
[Dem87]. Let $\varphi$ be a continuous psh function
with an isolated $-\infty$ pole at $x$, e.g.\ a function of the form
$\varphi(z)=\log\sum_{1\le j\le N}|g_j(z)|^{\gamma_j}$, $\gamma_j>0$,
where $(g_1\ld g_N)$ is an ideal of germs of holomorphic functions in
$\cO_x$ with $g^{-1}(0)=\{x\}$. The {\it generalized
Lelong number} $\nu(\Theta,\varphi)$ of $\Theta$ with respect to the
weight $\varphi$ is simply defined to be the mass of the measure
$\Theta\wedge(dd^c\varphi)^p$ carried by the point $x$
(the measure $\Theta\wedge(dd^c\varphi)^p$ is always well defined thanks
to Prop.~2.3). This number can also
be seen as the limit $\nu(\Theta,\varphi)=\lim_{t\to-\infty}
\nu(\Theta,\varphi,t)$, where
$$\nu(\Theta,\varphi,t)=\int_{\varphi(z)<t}\Theta\wedge
(dd^c\varphi)^p.\leqno(2.11)$$
The relation with our earlier definition of Lelong numbers (as well as part
a) of Theorem~2.8) comes from the identity
$$\nu(\Theta,x,r)=\nu(\Theta,\varphi,\log r),~~~~
\varphi(z)=\log|z-x|,\leqno(2.12)$$
in particular $\nu(\Theta,x)=\nu(\Theta,\log|\bu-x|)$. This equality is
in turn a consequence of the following general formula, applied to
$\chi(t)=e^{2t}$ and $t=\log r\,$:
$$\int_{\varphi(z)<t}\Theta\wedge(dd^c\chi\circ\varphi)^p=
\chi'(t-0)^p\int_{\varphi(z)<t}\Theta\wedge(dd^c\varphi)^p,
\leqno(2.13)$$
where $\chi$ is an arbitrary convex increasing function. To prove the
formula, we use a regularization and thus suppose that $\Theta$,
$\varphi$ and $\chi$ are smooth, and that $t$ is a non critical value of
$\varphi$. Then Stokes' formula shows that the integrals on the left and
right hand side of (2.13) are equal respectively to
$$\int_{\varphi(z)=t}\Theta\wedge(dd^c\chi\circ\varphi\big)^{p-1}
\wedge d^c(\chi\circ\varphi),~~~~
\int_{\varphi(z)=t}\Theta\wedge\big(dd^c\varphi\big)^{p-1}
\wedge d^c\varphi,$$
and the differential form of bidegree $(p-1,p)$ appearing in the integrand
of the first integral is equal to
$(\chi'\circ\varphi)^p\,(dd^c\varphi)^{p-1}\wedge d^c\varphi$.
The expected formula follows. Part b) of Theorem 2.8 is a consequence of
the Jensen-Lelong formula, whose proof is left as an exercise to the reader.

\begstat{(2.14) Jensen-Lelong formula} Let $u$ be any psh function
on~$X$. Then $u$ is integrable with respect to the measure
$\mu_r=(dd^c\varphi)^{n-1}\wedge d^c\varphi$ supported by the
pseudo-sphere $\{\varphi(z)=r\}$ and
$$\mu_r(u)=\int_{\{\varphi<r\}}u(dd^c\varphi)^n+
\int_{-\infty}^r\nu(dd^cu,\varphi,t)\,dt.\eqno\square$$
\endstat

In our case, we set $\varphi(z)=\log|z-x|$. Then $(dd^c\varphi)^n=\delta_x$
and $\mu_r$ is just the unitary invariant mean value measure on the
sphere $S(x,e^r)$. For $r<r_0$, Formula 2.14 implies
$$\mu_r(u)-\mu_{r_0}(u)=\int_{r_0}^r\nu(dd^cu,x,t)\sim
(r-r_0)\nu(dd^cu,x)\qquad\hbox{as $r\to-\infty$}.$$
{}From this, using the Harnack inequality for subharmonic functions, we get
$$\liminf_{z\to x}{u(z)\over\log|z-x|}=
\lim_{r\to-\infty}{\mu_r(u)\over r}=\nu(dd^cu,x).$$
These equalities imply statement 2.8~b).

Next, we show that the Lelong numbers $\nu(T,\varphi)$ only
depend on the asymptotic behaviour of $\varphi$ near the polar set
$\varphi^{-1}(-\infty)$. In a precise way:

\begstat{(2.15) Comparison theorem} Let $\Theta$ be a closed positive
current on~$X$, and let $\varphi,\psi:X\to[-\infty,+\infty[$
be continuous psh functions with isolated poles at some point
$x\in X$. Assume that
$$\ell:=\limsup_{z\to x}{\psi(z)\over\varphi(z)}<+\infty.$$
Then $\nu(\Theta,\psi)\le \ell^p\nu(\Theta,\varphi)$, and the equality holds
if $\ell=\lim\psi/\varphi$.
\endstat

\begproof{} (2.12) shows that
$\nu(\Theta,\lambda\varphi)=\lambda^p\nu(\Theta,\varphi)$
for every positive constant $\lambda$. It is thus sufficient to verify the
inequality $\nu(\Theta,\psi)\le\nu(\Theta,\varphi)$ under the hypothesis
$\limsup\psi/\varphi<1$. For any $c>0$, consider the psh function
$$u_c=\max(\psi-c,\varphi).$$
Fix $r\ll 0$. For $c>0$ large enough, we have $u_c=\varphi$ on a
neighborhood of $\varphi^{-1}(r)$ and Stokes' formula gives
$$\nu(\Theta,\varphi,r)=\nu(\Theta,u_c,r)\ge\nu(\Theta,u_c).$$
The hypothesis $\limsup\psi/\varphi<1$ implies on the other hand that
there exists $t_0<0$ such that $u_c=\psi-c$ on $\{u_c<t_0\}$. We thus get
$$\nu(\Theta,u_c)=\nu(\Theta,\psi-c)=\nu(\Theta,\psi),$$
hence $\nu(\Theta,\psi)\le\nu(\Theta,\varphi)$. The equality case is obtained
by reversing the roles of $\varphi$ and $\psi$ and observing that
$\lim\varphi/\psi=1/l$.\qed
\endproof

Part a) of Theorem~2.10 follows immediately from 2.15 by considering
the weights $\varphi(z)=\log|\tau(z)-\tau(x)|$, $\psi(z)=\log
|\tau'(z)-\tau'(x)|$ associated to coordinates systems
$\tau(z)=(z_1\ld z_n)$, $\tau'(z)=(z'_1\ld z'_n)$ in a neighborhood
of~$x$. Another application is a direct simple proof of Thie's Theorem
2.9 when $\Theta=[A]$ is the current of integration over an
analytic set $A\subset X$ of pure dimension~$p$. For this, we have to
observe that Theorem 2.15 still holds provided that $x$ is an isolated
point in $\Supp(\Theta)\cap\varphi^{-1}(-\infty)$ and
$\Supp(\Theta)\cap\psi^{-1}(-\infty)$ (even though $x$ is not isolated in
$\varphi^{-1}(-\infty)$ or $\psi^{-1}(-\infty)$), under the weaker assumption
that $\limsup_{\Supp(\Theta)\ni z\to x}\psi(z)/\varphi(z)=\ell$. The reason
for this is that all integrals involve currents supported on $\Supp(\Theta)$.
Now, by a generic choice of local coordinates \hbox{$z'=(z_1\ld z_p)$}
and \hbox{$z''=(z_{p+1}\ld z_n)$} on $(X,x)$, the germ $(A,x)$ is contained
in a cone \hbox{$|z''|\le C|z'|$}. If
$B'\subset\bbbc^p$ is a ball of center $0$ and radius $r'$ small,
and \hbox{$B''\subset\bbbc^{n-p}$} is the ball of center $0$ and
radius~$r''=Cr'$, the projection
      $$\pr:A\cap(B'\times B'')\longrightarrow B'$$
is a ramified covering with finite sheet number $m$.
When $z\in A$ tends to $x=0$, the functions
$$\varphi(z)=\log|z|=\log(|z'|^2+|z''|^2)^{1/2},~~~~\psi(z)=\log|z'|.$$
satisfy $\lim_{z\to x}\psi(z)/\varphi(z)=1$. Hence Theorem~2.15 implies
      $$\nu([A],x)=\nu([A],\varphi)=\nu([A],\psi).$$
Now, Formula 2.13 with $\chi(t)=e^{2t}$ yields
$$\eqalign{\nu([A],\psi,\log t)
&=t^{-2p}\int_{\{\psi<\log t\}}[A]\wedge\Big({1\over 2}dd^ce^{2\psi}\Big)^p\cr
&=t^{-2p}\int_{A\cap\{|z'|<t\}}\Big({1\over 2}\pr^\star dd^c|z'|^2\Big)^p\cr
&=m\,t^{-2p}\int_{\bbbc^p\cap\{|z'|<t\}}\Big({1\over 2}dd^c|z'|^2\Big)^p=m,
\cr}$$
hence $\nu([A],\psi)=m$.  Here, we have used the fact that pr is an \'etale
covering with $m$ sheets over the complement of the ramification locus
$S\subset B'$, and the fact that $S$ is of zero Lebesgue measure in~$B'$.

\begstat{(2.16) Proposition} Under the assumptions of Proposition 2.3,
we have
$$\nu(dd^cu\wedge\Theta,x)\ge\nu(u,x)\,\nu(\Theta,x)$$
at every point $x\in X$.
\endstat

\begproof{} Assume that $X=B(0,r)$ and $x=0$. By definition
$$\nu(dd^c u\wedge\Theta,x)=
\lim_{r\to 0}\int_{|z|\le r}dd^cu\wedge\Theta\wedge(dd^c\log|z|)^{p-1}.$$
Set $\gamma=\nu(u,x)$ and
$$u_\nu(z)=\max\big(u(z),(\gamma-\varepsilon)\log|z|-\nu\big)$$
with $0<\varepsilon<\gamma$ (if $\gamma=0$, there is nothing to prove).
Then $u_\nu$ decreases to $u$ and
$$\int_{|z|\le r}\kern-7pt dd^c u\wedge\Theta\wedge(dd^c\log|z|)^{p-1}\ge
\limsup_{\nu\to+\infty}\int_{|z|\le r}\kern-7pt dd^cu_\nu\wedge
\Theta\wedge(dd^c\log|z|)^{p-1}$$
by the weak convergence of $dd^cu_\nu\wedge\Theta$; here
$(dd^c\log|z|)^{p-1}$ is not
smooth on $\ovl B(0,r)$, but the integrals remain unchanged if we replace
$\log|z|$ by $\chi(\log|z|/r)$ with a smooth convex function $\chi$ such
that $\chi(t)=t$ for $t\ge -1$ and $\chi(t)=0$ for $t\le -2$. Now, we have
$u(z)\le\gamma\log|z|+C$ near $0$, so $u_\nu(z)$ coincides with
$(\gamma-\varepsilon)\log|z|-\nu$ on a small ball $B(0,r_\nu)\subset
B(0,r)$ and we infer
$$\eqalign{
\int_{|z|\le r}dd^c u_\nu\wedge\Theta\wedge
(dd^c\log|z|)^{p-1}&\ge(\gamma-\varepsilon)
\int_{|z|\le r_\nu}\Theta\wedge(dd^c\log|z|)^p\cr
&\ge(\gamma-\varepsilon)\nu(\Theta,x).\cr}$$
As $r\in{}]0,R[$ and $\varepsilon\in{}]0,\gamma[$ were arbitrary, the
desired inequality follows.\qed
\endproof

We will later need an important decomposition formula of [Siu74].
We start with the following lemma.

\begstat{(2.17) Lemma} If $\Theta$ is a closed
positive current of bidimension $(p,p)$ and $Z$ is an irreducible
analytic set in $X$, we set
$$m_Z=\inf\{x\in Z\,;\,\nu(\Theta,x)\}.$$
\item{\rm a)} There is a countable family of proper analytic subsets
$(Z'_j)$ of $Z$ such that $\nu(\Theta,x)=m_Z$ for all
$x\in Z\ssm\bigcup Z'_j$. We say that $m_Z$ is the generic Lelong number
of $\Theta$ along $Z$.
\smallskip
\item{\rm b)} If $\dim Z=p$, then $\Theta\ge m_Z[Z]$ and
$\bbbone_Z\Theta=m_Z[Z]$.\vskip0pt
\endstat

\begproof{} a) By definition of $m_Z$ and $E_c(\Theta)$, we
have $\nu(\Theta,x)\ge m_Z$ for every $x\in Z$ and
$$\nu(\Theta,x)=m_Z~~~~\hbox{\rm on}~~
Z\ssm\bigcup_{c\in\bbbq,\,c>m_Z}Z\cap E_c(\Theta).$$
However, for $c>m_Z$, the intersection $Z\cap E_c(\Theta)$ is a proper
analytic subset of~$A$.

\noindent b) Left as an exercise to the reader. It is enough to prove
that $\Theta\ge m_Z[Z_\reg]$ at regular points of $Z$, so one may assume
that $Z$ is a $p$-dimensional linear subspace in~$\bbbc^n$. Show that
the measure $(\Theta-m_Z[Z])\wedge(dd^c|z|^2)^p$ has nonnegative mass
on every ball $|z-a|<r$ with center $a\in Z$. Conclude by using
arbitrary affine changes of coordinates that $\Theta-m_Z[Z]\ge 0$.\qed
\endproof

\begstat{(2.18) Decomposition formula {\rm([Siu74])}} Let $\Theta$
be a closed positive current of bidimension~$(p,p)$. Then
$\Theta$ can be written as a convergent series of
closed positive currents
$$\Theta=\sum_{k=1}^{+\infty}\lambda_k\,[Z_k]+R,$$
where $[Z_k]$ is a current of integration over an irreducible
analytic set of dimension~$p$, and $R$ is a residual current
with the property that $\dim E_c(R)<p$ for every $c>0$.
This decomposition is locally and globally unique: the sets $Z_k$
are precisely the $p$-dimensional components occurring in the
sublevel sets $E_c(\Theta)$, and $\lambda_k=\min_{x\in Z_k}\nu(\Theta,x)$ is
the generic Lelong number of $\Theta$ along $Z_k$.
\endstat

\begproof{of uniqueness.} If $\Theta$ has such a decomposition, the
$p$-dimensional components of $E_c(\Theta)$ are $(Z_j)_{\lambda_j\ge c}$,
for $\nu(\Theta,x)=\sum\lambda_j\nu([Z_j],x)+\nu(R,x)$ is non zero
only on $\bigcup Z_j\cup\bigcup E_c(R)$, and is equal to
$\lambda_j$ generically on $Z_j$ $\big($more precisely,
$\nu(\Theta,x)=\lambda_j$ at every regular
point of $Z_j$ which does not belong to any intersection $Z_j\cup Z_k$,
$k\ne j$ or to $\bigcup E_c(R)\big)$. In particular $Z_j$ and
$\lambda_j$ are unique.
\endproof

\begproof{of existence.} Let $(Z_j)_{j\ge 1}$ be the countable collection of
$p$-dimensional components occurring in one of the sets $E_c(\Theta)$,
$c\in\bbbq_+^\star$, and let $\lambda_j>0$ be the generic Lelong
number of $\Theta$ along $Z_j$. Then Lemma~2.17 shows by induction on $N$
that $R_N=\Theta-\sum_{1\le j\le N}\lambda_j[Z_j]$ is positive. As $R_N$
is a decreasing sequence, there must be a
limit $R=\lim_{N\to+\infty}R_N$ in the weak topology. Thus we have
the asserted decomposition. By construction, $R$ has zero generic
Lelong number along $Z_j$, so $\dim E_c(R)<p$ for every $c>0$.\qed
\endproof

It is very important to note that some components of lower dimension
can actually occur in $E_c(R)$, but they cannot be subtracted because
$R$ has bidimension $(p,p)$. A typical case is the case of a
bidimension \hbox{$(n-1,n-1)$} current $\Theta=dd^cu$ with
$u=\log(|f_j|^{\gamma_1}+\ldots|f_N|^{\gamma_N})$ and
$f_j\in H^0(X,\cO_X)$. In general $\bigcup E_c(\Theta)=\bigcap f_j^{-1}(0)$
has dimension~$<n-1$.

\begstat{Corollary 2.19} Let $\Theta_j=dd^c u_j$, $1\le j\le p$,
be closed positive \hbox{$(1,1)$-currents}
on a complex manifold $X$. Suppose that there are analytic sets
$A_2\supset\ldots\supset A_p$ in $X$ with
$\codim A_j\ge j$ at every point such that each $u_j$, $j\ge 2$,
is locally bounded on $X\ssm A_j$. Let $\{A_{p,k}\}_{k\ge 1}$ be the
irreducible components of $A_p$ of codimension $p$ exactly and let
$\nu_{j,k}=\min_{x\in A_{p,k}}\nu(\Theta_j,x)$ be the generic Lelong number
of $\Theta_j$ along $A_{p,k}$. Then $\Theta_1\wedge\ldots\wedge \Theta_p$
is well-defined and
$$\Theta_1\wedge\ldots\wedge \Theta_p\ge\sum_{k=1}^{+\infty}\nu_{1,k}\ldots
\nu_{p,k}\,[A_{p,k}].$$
\endstat

\begproof{} By induction on $p$, Prop.~2.3 shows that
$\Theta_1\wedge\ldots\wedge \Theta_p$ is well defined. Moreover, Prop.~2.16
implies
$$\nu(\Theta_1\wedge\ldots\wedge \Theta_p,x)\ge
\nu(\Theta_1,x)\ldots\nu(\Theta_p,x)\ge\nu_{1,k}\ldots\nu_{p,k}$$
at every point $x\in A_{p,k}$. The desired inequality is then a consequence
of Siu's decomposition theorem.\qed
\endproof

\titleb{3.}{Hermitian Vector Bundles, Connections and Curvature}
The goal of this section is to recall the most basic definitions of
hemitian differential geometry related to the concepts of connection,
curvature and first Chern class of a line bundle.

Let $F$ be a complex vector bundle of rank $r$ over a smooth
differentiable manifold~$M$. A {\it connection} $D$ on $F$ is a
linear differential operator of order $1$
$$D:C^\infty(M,\Lambda^qT^\star_M\otimes F) \to C^\infty(M,
\Lambda^{q+1}T^\star_M\otimes F)$$
such that
$$D(f\wedge u) = df\wedge u + (-1)^{\deg\,f}f\wedge Du\leqno(3.1)$$
for all forms $f\in C^\infty(M,\Lambda^pT^\star_M)$, $u\in C^\infty(X,
\Lambda^qT^\star_M\otimes F)$. On an open set $\Omega \subset M$ where $F$
admits a trivialization $\theta:F_{|\Omega }\buildo\simeq\over\lra\Omega
\times \bbbc^r$, a connection $D$ can be written
$$Du\simeq_\theta du+\Gamma \wedge u$$
where $\Gamma \in C^\infty(\Omega ,\Lambda^1T^\star_M\otimes\Hom
(\bbbc^r,\bbbc^r))$ is an arbitrary matrix of $1$-forms and $d$ acts
componentwise. It is then easy to check that
$$D^2u\simeq_\theta (d\Gamma +\Gamma \wedge \Gamma )\wedge u\quad
\hbox{on~}\Omega.$$
Since $D^2$ is a globally defined operator, there is a global $2$-form
$$\Theta(D) \in C^\infty(M,\Lambda^2T^\star_M\otimes\Hom(F,F))\leqno(3.2)$$
such that $D^2u=\Theta(D)\wedge u$ for every form $u$ with values in~$F$.

Assume now that $F$ is endowed with a $C^\infty$ hermitian metric along
the fibers and that the isomorphism $F_{|\Omega }\simeq \Omega \times
\bbbc^r$ is given by a $C^\infty$ frame $(e_\lambda)$. We then have
a canonical sesquilinear pairing
$$\leqalignno{
C^\infty(M,\Lambda^pT^\star_M\otimes F)\times
C^\infty(M,\Lambda^qT^\star_M\otimes F)
&\longrightarrow C^\infty(M,\Lambda^{p+q}T^\star_M\otimes\bbbc)&(3.3)\cr
(u,v) &\longmapsto \{u,v\}\cr}$$
given by
$$\{u,v\}=\sum_{\lambda ,\mu }u_\lambda \wedge \overline v_\mu \langle
e_\lambda,e_\mu \rangle,\qquad
u=\sum u_\lambda \otimes e_\lambda,\quad
v=\sum v_\mu \otimes e_\mu.$$
The connection $D$ is said to be {\it hermitian} if it satisfies the
additional property
$$d\{u,v\}=\{Du,v\}+(-1)^{{\rm deg~}u}\{u,Dv\}.$$
Assuming that $(e_\lambda )$ is orthonormal, one easily checks that $D$
is hermitian if and only if $\Gamma^\star=-\Gamma$. In this case
$\Theta(D)^\star=-\Theta(D)$, thus
$$\ii\Theta(D) \in C^\infty(M,\Lambda^2T^\star_M\otimes\Herm(F,F)).$$
\medskip

\begstat{(3.4) Special case} \rm For a bundle $F$ of rank 1, the
connection form $\Gamma$ of a hermitian connection $D$ can be seen as
a $1$-form with purely imaginary coefficients $\Gamma=\ii A$ ($A$ real).
Then we have $\Theta(D)=d\Gamma=\ii dA$. In particular $\ii\Theta(F)$
is a closed 2-form. The {\it First Chern class} of $F$
is defined to be the cohomology class
$$c_1(F)_\bbbr=\Big\{{\ii\over 2\pi}\Theta(D)\Big\}\in H^2_\DR(M,\bbbr).$$
The cohomology class is actually independent of the connection, since
any other connection $D_1$ differs by a global 1-form, $D_1u=Du+B\wedge u$,
so that $\Theta(D_1)=\Theta(D)+dB$. It is well-known that $c_1(F)_\bbbr$
is the image in $H^2(M,\bbbr)$ of an integral class $c_1(F)\in
H^2(M,\bbbz)\,$; by using the exponential exact sequence
$$0\to\bbbz\to\cE\to\cE^\star\to 0,$$
$c_1(F)$ can be defined in \v Cech cohomology theory as the image
by the coboundary map $H^1(M,\cE^\star)\to H^2(M,\bbbz)$ of the cocycle
$\{g_{jk}\}\in H^1(M,\cE^\star)$ defining $F\,$; see e.g.\ [GH78] for
details.\qed
\endstat

We now concentrate ourselves on the complex analytic case. If $M=X$ is a
complex manifold $X$, every connection $D$ on a complex $\ci$ vector
bundle $F$ can be splitted in a unique way as a
sum of a $(1,0)$ and of a $(0,1)$-connection, $D=D'+D''$. In a local
trivialization $\theta$ given by a $C^\infty$ frame, one can write
$$\leqalignno{
D'u &\simeq_\theta d'u + \Gamma '\wedge u,&(3.5')\cr
D''u &\simeq_\theta d''u + \Gamma ''\wedge u,&(3.5'')\cr}$$
with $\Gamma =\Gamma '+\Gamma ''$. The connection is hermitian if and
only if $\Gamma '=-(\Gamma '')^\star$ in any orthonormal frame. Thus there
exists a unique hermitian connection $D$ corresponding to a prescribed
$(0,1)$ part $D''$.

Assume now that the bundle $F$ itself has a {\it holomorphic} structure.
The unique hermitian connection for which $D''$ is the $d''$ operator
defined in \S~1 is called the {\it Chern connection} of $F$. In a local
holomorphic frame $(e_\lambda )$ of $E_{|\Omega }$, the metric is given
by the hermitian matrix $H=(h_{\lambda\mu})$, $h_{\lambda\mu}
=\langle e_\lambda,e_\mu \rangle$. We have
$$\{u,v\}=\sum_{\lambda,\mu}h_{\lambda\mu}u_\lambda \wedge
\ovl v_\mu=u^\dagger\wedge H\ovl v,$$
where $u^\dagger$ is the transposed matrix of $u$, and
easy computations yield
$$\eqalign{d\{u,v\}
&=(du)^\dagger \wedge H\ovl v+(-1)^{\deg u} u^\dagger
\wedge (dH\wedge \ovl v + H\ovl{dv})\cr
&=\big(du +\ovl H^{-1}d'\ovl H \wedge u\big)^\dagger \wedge H\ovl v
+(-1)^{\deg u}u^\dagger\wedge(\ovl{dv+ \ovl H^{-1}d'\ovl H\wedge v)}\cr}$$
using the fact that $dH = d'H+\ovl{d'\ovl H}$ and
$\ovl H^\dagger=H$. Therefore the Chern connection $D$ coincides with the
hermitian connection defined by
$$\left\{\eqalign{
Du&{}\simeq_\theta du+\ovl H^{-1}d'\ovl H \wedge u,\cr
D'&{}\simeq_\theta d'+\ovl H^{-1} d'\ovl H\wedge\bu=
\ovl H^{-1} d'(\ovl H\bu),~~~~D'' = d''.\cr}\right.\leqno(3.6)$$
It is clear from this relations that $D^{\prime 2}=D^{\prime\prime 2}=0$.
Consequently $D^2$ is given by to $D^2=D'D''+D''D'$, and the curvature
tensor $\Theta(D)$ is of type~$(1,1)$. Since $d'd''+d''d'=0$, we get
$$\eqalign{
(D'D''+D''D')u&\simeq_\theta\ovl H^{-1}d'\ovl H\wedge d''u
+d''(\ovl H^{-1}d'\ovl H\wedge u)\cr
&=d''(\ovl H^{-1}d'\ovl H)\wedge u.\cr}$$

\begstat{(3.7) Proposition} The Chern curvature tensor $\Theta(F):=
\Theta(D)$ is such that
$$\ii\,\Theta(F) \in C^\infty(X,\Lambda^{1,1}T^\star_X\otimes\Herm(F,F)).$$
If $\theta: E_{\restriction\Omega}\to\Omega\times\bbbc^r$ is a holomorphic
trivialization and if $H$ is the hermitian matrix representing the
metric along the fibers of $F_{\restriction\Omega}$, then
$$\ii\,\Theta(F)\simeq_\theta
\ii\,d''(\ovl H^{-1} d'\ovl H)\quad\hbox{\rm on}~\Omega.\eqno\square$$
\endstat

Let $(z_1\ld z_n)$ be holomorphic coordinates on $X$ and let
$(e_\lambda)_{1\le\lambda\le r}$ be an orthonormal frame of $F$.
Writing
$$\ii\Theta(F)=\sum_{1\le j,k\le n,\,1\le\lambda,\mu\le r}
c_{jk\lambda\mu}dz_j\wedge dz_k\otimes e_\lambda^\star\otimes e_\mu,$$
we can identify the curvature tensor to a hermitian form
$$\wt\Theta(F)(\xi\otimes v)=\sum_{1\le j,k\le n,\,1\le\lambda,\mu\le r}
c_{jk\lambda\mu}\xi_j\ovl\xi_kv_\lambda\ovl v_\mu\leqno(3.8)$$
on $T_X\otimes F$. This leads in a natural way to positivity concepts,
following definitions introduced by Kodaira [Kod53], Nakano [Nak55] and
Griffiths [Gri66].

\begstat{(3.9) Definition} The hermitian vector bundle $F$ is said to be
\smallskip
\item{\rm a)} positive in the sense of Nakano if
$\wt\Theta(F)(\tau)>0$ for all non zero tensors $\tau=\sum
\tau_{j\lambda}\partial/\partial z_j\otimes e_\lambda\in T_X\otimes F$.
\smallskip
\item{\rm b)} positive in the sense of Griffiths if
$\wt\Theta(F)(\xi\otimes v)>0$ for all non zero decomposable tensors
$\xi\otimes v\in T_X\otimes F\,;$
\smallskip\noindent
Corresponding semipositivity concepts are defined by relaxing the strict
inequalities.
\endstat

\begstat{(3.10) Special case of rank 1 bundles} \rm Assume that
$F$ is a line bundle. The hermitian matrix $H=(h_{11})$ associated
to a trivialization $\theta:F_{\restriction\Omega}\simeq\Omega\times\bbbc$
is simply a positive function which we find convenient to denote by
$e^{-2\varphi}$, \hbox{$\varphi\in C^\infty(\Omega,\bbbr)$}. In this case
the curvature form $\Theta(F)$ can be identified to the $(1,1)$-form
$2d'd''\varphi$, and
$${\ii\over 2\pi}\Theta(F)={\ii\over\pi}d'd''\varphi=dd^c\varphi$$
is a real $(1,1)$-form. Hence $F$ is semipositive (in either Nakano or
Griffiths sense) if and only if $\varphi$ is psh, resp.\ positive if
and only if $\varphi$ is {\it strictly psh}. In this setting, the
Lelong-Poincar\'e equation can be generalized as follows: let
$\sigma\in H^0(X,F)$ be a non zero holomorphic section. Then
$$dd^c\log\|\sigma\|=[Z_\sigma]-{\ii\over 2\pi}\Theta(F).\leqno(3.11)$$
Formula (3.11) is immediate if we write $\|\sigma\|=|\theta(\sigma)|
e^{-\varphi}$ and if we apply (1.20) to the holomorphic function
$f=\theta(\sigma)$. As we shall see later, it is very important for
applications to consider also singular hermitian metrics.
\endstat

\begstat{(3.12) Definition} A singular $($hermitian$)$ metric on a line
bundle $F$ is a metric which is given in any trivialization
$\theta:F_{\restriction\Omega}
\buildo\simeq\over\longrightarrow\Omega\times\bbbc$ by
$$\|\xi\|=|\theta(\xi)|\,e^{-\varphi(x)},~~~~x\in\Omega,~\xi\in F_x$$
where $\varphi\in L^1_\loc(\Omega)$ is an arbitrary function, called the
weight of the metric with respect to the trivialization $\theta$.
\endstat

If $\theta':F_{\restriction\Omega'}\lra\Omega'\times\bbbc$
is another trivialization, $\varphi'$ the associated weight and
$g\in\cO^\star(\Omega\cap\Omega')$ the transition function, then
$\theta'(\xi)=g(x)\,\theta(\xi)$ for $\xi\in F_x$, and so
$\varphi'=\varphi+\log|g|$ on $\Omega\cap\Omega'$. The curvature form
of $F$ is then given formally by the closed $(1,1)$-current
${\ii\over 2\pi}\Theta(F)=dd^c\varphi$ on $\Omega\,$; our assumption
$\varphi\in L^1_\loc(\Omega)$ guarantees that $\Theta(F)$ exists in
the sense of distribution theory. As in the smooth case, ${\ii\over 2\pi}
\Theta(F)$ is globally defined on $X$ and independent of the choice of
trivializations, and its De Rham cohomology class is the image of the first
Chern class $c_1(F)\in H^2(X,\bbbz)$ in $H^2_{DR}(X,\bbbr)$.
Before going further, we discuss two basic examples.

\begstat{(3.13) Example} \rm Let $D=\sum\alpha_jD_j$ be a divisor with
coefficients $\alpha_j\in\bbbz$ and let $F=\cO(D)$ be the associated
invertible sheaf of meromorphic functions $u$ such that
${\rm div}(u)+D\ge 0\,$; the corresponding line bundle can be
equipped with the singular metric defined by $\|u\|=|u|$. If $g_j$ is a
generator of the ideal of $D_j$ on an open set $\Omega\subset X$ then
$\theta(u)=u\prod\smash{g_j^{\alpha_j}}$ defines a trivialization of
$\cO(D)$ over $\Omega$, thus our singular metric is associated to the
weight $\varphi=\sum\alpha_j\log|g_j|$. By the Lelong-Poincar\'e equation,
we find
$${\ii\over 2\pi}\Theta\big(\cO(D)\big)=dd^c\varphi=[D],$$
where $[D]=\sum\alpha_j[D_j]$ denotes the current of integration
over~$D$.\qed
\endstat

\begstat{(3.14) Example} \rm Assume that $\sigma_1\ld \sigma_N$ are non zero
holomorphic sections of $F$. Then we can define a natural (possibly
singular) hermitian metric
on $F^\star$ by
$$\|\xi^\star\|^2=\sum_{1\le j\le n}\big|\xi^\star.\sigma_j(x)\big|^2
{}~~~~\hbox{\rm for}~~\xi^\star\in F^\star_x.$$
The dual metric on $F$ is given by
$$\|\xi\|^2={|\theta(\xi)|^2\over|\theta(\sigma_1(x))|^2+\ldots+
|\theta(\sigma_N(x))|^2}$$
with respect to any trivialization $\theta$. The associated weight function is
thus given by $\varphi(x)=\log\big(\sum_{1\le j\le N}|\theta(\sigma_j(x))|^2
\big){}^{1/2}$. In this case $\varphi$ is a psh function, thus
$\ii\Theta(F)$ is a closed positive current. Let us denote by $\Sigma$
the linear system defined by $\sigma_1\ld \sigma_N$ and by
$B_\Sigma=\bigcap\sigma_j^{-1}(0)$ its base locus. We have a meromorphic
map
$$\Phi_\Sigma:X\ssm B_\Sigma\to\bbbp^{N-1},\qquad
x\mapsto(\sigma_1(x):\sigma_2(x):\ldots:\sigma_N(x)).$$
Then ${\ii\over 2\pi}\Theta(F)$ is equal to the pull-back over $X\ssm
B_\Sigma$ of the Fubini-Study metric $\omega_\FS={\ii\over 2\pi}
\log(|z_1|^2+\ldots+ |z_N|^2)$ of $\bbbp^{N-1}$ by $\Phi_\Sigma$.\qed
\endstat

\begstat{(3.15) Ample and very ample line bundles} A holomorphic line
bundle $F$ over a compact complex manifold $X$ is said to be
\smallskip
\item{\rm a)} very ample if the map $\Phi_{|F|}:X\to\bbbp^{N-1}$ associated
to the complete linear system $|F|=P(H^0(X,F))$ is a regular embedding
$($by this we mean in particular that the base locus is empty, i.e.\
$B_{|F|}=\emptyset)$.
\smallskip
\item{\rm b)} ample if some multiple $mF$, $m>0$, is very ample.
\endstat

Here we use an additive notation for $\Pic(X)=H^1(X,\cO^\star)$, hence
the symbol $mF$ denotes the line bundle $F^{\otimes m}$. By
Example~3.15, every ample line bundle $F$ has a smooth hermitian metric
with positive definite curvature form; indeed, if the linear system
$|mF|$ gives an embedding in projective space, then we get a smooth
hermitian metric on $F^{\otimes m}$, and the $m$-th root yields a
metric on $F$ such that ${\ii\over 2\pi}\Theta(F)={1\over m}
\Phi_{|mF|}^\star\omega_\FS$. Conversely, the Kodaira embedding theorem
[Kod54] tells us that every positive line bundle $F$ is ample (see
Exercise~5.14 for a straightforward analytic proof of the Kodaira
embedding theorem).

\titleb{4.}{Bochner Technique and Vanishing Theorems}
We first recall briefly a few basic facts of Hodge theory. Assume for
the moment that $M$ is a differentiable manifold equipped with a riemannian
metric $g=\sum g_{ij}dx_i\otimes dx_j$. Given a $q$-form $u$ on $M$ with
values in $F$, we consider the global $L^2$ norm
$$\|u\|^2 = \int_M |u(x)|^2dV_g(x)$$
where $|u|$ is the pointwise hermitian norm and $dV_g$ is the
riemannian volume form. The Laplace-Beltrami operator associated to the
connection $D$ is
$$\Delta =DD^\star+D^\star D$$
where
$$D^\star:\ci(M,\Lambda^qT^\star_M\otimes F)\to
\ci(M,\Lambda^{q-1}T^\star_M\otimes F)$$
is the (formal) adjoint of $D$ with respect to the $L^2$ inner product.
Assume that $M$ is {\it compact}. Since
$$\Delta:\ci(M,\Lambda^qT^\star_M\otimes F)\to
\ci(M,\Lambda^qT^\star_M\otimes F)$$
is a self-adjoint elliptic operator in each degree, standard results of
PDE theory show that there is an orthogonal decomposition
$$\ci(M,\Lambda^qT^\star_M\otimes F)=\cH^q(M,F)\oplus\Im\Delta$$
where $\cH^q(M,F)=\Ker\Delta$ is the space of harmonic forms of degree~$q$;
$\cH^q(M,F)$ is a finite dimensional space.
Assume moreover that the connection $D$ is {\it integrable}, i.e.\ that
$D^2=0$. It is then easy to check that there is an orthogonal direct sum
$$\Im\Delta=\Im D\oplus\Im D^\star,$$
indeed $\langle Du,D^\star v\rangle=\langle D^2u,v\rangle=0$ for all $u,v$.
Hence we get an orthogonal decomposition
$$\ci(M,\Lambda^qT^\star_M\otimes F)=\cH^q(M,F)\oplus\Im D\oplus\Im D^\star,$$
and $\Ker\Delta$ is precisely equal to $\cH^q(M,F)\oplus\Im D$. Especially,
the $q$-th cohomology group $\Ker\Delta/\Im\Delta$ is isomorphic to
$\cH^q(M,F)$. All this can be applied for example in the case of the
De Rham groups $H^q_\DR(M,\bbbc)$, taking $F$ to be the trivial bundle
$F=M\times\bbbc$ (notice, however, that a nontrivial bundle $F$
usually does not admit any integrable connection):

\begstat{(4.1) Hodge Fundamental Theorem} If $M$ is a compact riemannian
manifold, there is an isomorphism
$$H^q_\DR(M,\bbbc)\simeq\cH^q(M,\bbbc)$$
from De Rham cohomology groups onto spaces of harmonic forms.\qed
\endstat

A rather important consequence of the Hodge fundamental theorem is a proof
of the {\it Poincar\'e duality theorem}. Assume that the Riemannian manifold
$(M,g)$ is oriented. Then there
is a (conjugate linear) Hodge star operator
$$\star:\Lambda^qT^\star_M\otimes\bbbc\to\Lambda^{m-q}T^\star_M\otimes\bbbc,
\qquad m=\dim_\bbbr M$$
defined by $u\wedge\star v=\langle u,v\rangle dV_g$ for any two complex
valued $q$-forms $u$, $v$. A standard computation shows
that $\star$ commutes with $\Delta$, hence $\star u$ is harmonic if and only
if $u$ is. This implies that the natural pairing
$$H^q_\DR(M,\bbbc)\times H^{m-q}_\DR(M,\bbbc),\qquad
(\{u\},\{v\})\mapsto\int_M u\wedge v\leqno(4.2)$$
is a nondegenerate duality, the dual of a class $\{u\}$ represented by a
harmonic form being $\{\star u\}$.

Let us now suppose that $X$ is a compact complex manifold
equipped with a hermitian metric $\omega=\sum\omega_{jk}dz_j\wedge d\ovl
z_k$. Let $F$ be a holomorphic vector bundle on $X$ equipped with a
hermitian metric, and let $D=D'+D''$ be its Chern curvature form.
All that we said above for the Laplace-Beltrami operator $\Delta$
still applies to the complex Laplace operators
$$\Delta'=D'D^{\prime\star}+D^{\prime\star}D',\qquad
\Delta''=D''D^{\prime\prime\star}+D^{\prime\prime\star}D'',$$
with the great advantage that we always have
$D^{\prime 2}=D^{\prime\prime 2}=0$. Especially, if $X$ is a compact
complex manifold, there are isomorphisms
$$H^{p,q}(X,F)\simeq\cH^{p,q}(X,F)\leqno(4.3)$$
between Dolbeault cohomology groups $H^{p,q}(X,F)$ and spaces
$\cH^{p,q}(X,F)$ of $\Delta''$-harmonic forms of bidegree $(p,q)$ with
values in~$F$. Now, there is a generalized Hodge star operator
$$\star:\Lambda^{p,q}T^\star_X\otimes F\to
\Lambda^{n-p,n-q}T^\star_X\otimes F^\star,\qquad n=\dim_\bbbc X,$$
such that $u\wedge\star v=\langle u,v\rangle dV_g$, when the
for any two $F$-valued $(p,q)$-forms, when the wedge product
$u\wedge\star v$ is combined with the pairing $F\times F^\star\to\bbbc$.
This leads to the {\it Serre duality theorem} [Ser55]: the bilinear pairing
$$H^{p,q}(X,F)\times H^{n-p,n-q}(X,F^\star),\qquad
(\{u\},\{v\})\mapsto\int_X u\wedge v\leqno(4.4)$$
is a nondegenerate duality. Combining this with the Dolbeault isomorphism,
we may restate the result in the form of the duality formula
$$H^q(X,\Omega^p_X\otimes\cO(F))^\star\simeq H^{n-q}(X,
\Omega^{n-p}_X\otimes\cO(F^\star)).\leqno(4.4')$$

We now proceed to explain the basic ideas of the Bochner technique used
to prove vanishing theorems. Great simplifications occur in the computations
if the hermitian metric on $X$ is supposed to be {\it K\"ahler}, i.e.\
if the associated {\it fundamental $(1,1)$-form}
$$\omega=\ii\sum\omega_{jk}dz_j\wedge d\ovl z_k$$
satisfies $d\omega=0$. It can be easily shown that $\omega$ is K\"ahler
if and only if there are holomorphic coordinates $(z_1\ld z_n)$ centered
at any point $x_0\in X$ such that the matrix of coefficients $(\omega_{jk})$
is tangent to identity at order $2$, i.e.\
$$\omega_{jk}(z)=\delta_{jk}+O(|z|^2)\quad\hbox{at $x_0$}.$$
It follows that all order 1 operators $D$, $D'$, $D''$ and their
adjoints $D^\star$, $D^{\prime\star}$, $D^{\prime\prime\star}$ admit at $x_0$
the same expansion as the analogous operators obtained when all hermitian
metrics on $X$ or $F$ are constant. From this, the basic commutation
relations of K\"ahler geometry can be checked. If $A,B$ are differential
operators acting on the algebra $C^\infty(X,\Lambda^{\bu,\bu}T^\star_X\otimes
F)$, their graded commutator (or graded Lie bracket) is defined by
$$[A,B]= AB - (-1)^{ab} BA$$
where $a,b$ are the degrees of $A$ and $B$ respectively. If $C$ is another
endomorphism of degree $c$, the following purely formal {\it Jacobi identity}
holds:
$$(-1)^{ca} \big[A,[B,C]\big] + (-1)^{ab}\big[B,[C,A]\big]
+ (-1)^{bc} \big[C,[A,B]\big] = 0.$$

\begstat{(4.5) Basic commutation relations} Let $(X,\omega)$ be a K\"ahler
manifold and let $L$ be the operators defined by $Lu=\omega\wedge u$ and
$\Lambda=L^\star$. Then
$$\eqalign{
[D^{\prime\prime\star},L]&=\ii D',\cr
[\Lambda,D'']&=-\ii D^{\prime\star},\cr}
\qquad
\eqalign{
[D^{\prime\star},L]&=-\ii D'',\cr
[\Lambda,D']&=\ii D^{\prime\prime\star}.\cr}$$
\endstat

\begproof{(sketch).} The first step is to check the identity
$[d^{\prime\prime\star},L]=\ii d'$ for constant metrics on $X=\bbbc^n$ and
$F=X\times\bbbc$, by a brute force calculation. All three other identities
follow by taking conjugates or adjoints. The case of variable metrics follows
by looking at Taylor expansions up to order~1.\qed
\endproof

\begstat{(4.6) Bochner-Kodaira-Nakano identity} If $(X,\omega)$ is
K\"ahler, the complex Laplace operators $\Delta'$ and $\Delta''$ acting
on $F$-valued forms satisfy the identity
$$\Delta''=\Delta'+[\ii\Theta(F),\Lambda].$$
\endstat

\begproof{} The last equality in (4.5) yields $D^{\prime\prime\star}
=-\ii[\Lambda,D']$, hence
$$\Delta''=[D'',\delta'']=-\ii[D'',\big[\Lambda,D']\big].$$
By the Jacobi identity we get
$$\big[D'',[\Lambda,D']\big]=\big[\Lambda,[D',D'']]+\big[D',[D'',\Lambda]\big]
=[\Lambda,\Theta(F)]+\ii[D',D^{\prime\star}],$$
taking into account that $[D',D'']=D^2=\Theta(F)$. The formula follows.\qed
\endproof

Assume that $X$ is compact and that $u\in\ci(X,\Lambda^{p,q}T^\star
X\otimes F)$ is an arbitrary $(p,q)$-form. An integration by parts yields
$$\langle\Delta'u,u\rangle=\|D'u\|^2+\|D^{\prime\star}u\|^2\ge0$$
and similarly for $\Delta''$, hence we get the basic a priori inequality
$$\|D''u\|^2+\|D^{\prime\prime\star}u\|^2\ge\int_X
\langle[\ii\Theta(F),\Lambda]u,u\rangle dV_\omega.\leqno(4.7)$$
This inequality is known as the {\it Bochner-Kodaira-Nakano} inequality
(see [Boc48], [Kod53], [Nak55]). When $u$ is $\Delta''$-harmonic, we get
$$\int_X\big(\langle[\ii\Theta(F),\Lambda]u,u\rangle+
\langle T_\omega u,u\rangle\big)dV\le 0.$$
If the hermitian operator $[\ii\Theta(F),\Lambda]$ acting on
$\Lambda^{p,q}T^\star_X\otimes F$ is positive on each fiber, we infer
that $u$ must be zero, hence
$$H^{p,q}(X,F)=\cH^{p,q}(X,F)=0$$
by Hodge theory. The main point is thus to compute the curvature form
$\Theta(F)$ and find sufficient conditions under which the operator
$[\ii\Theta(F),\Lambda]$ is positive definite. Elementary (but somewhat
tedious) calculations yield the following formulae: if the curvature of
$F$ is written as in (3.8) and
$u=\sum u_{J,K,\lambda}dz_I\wedge d\ovl z_J\otimes e_\lambda$,
$|J|=p$, $|K|=q$, $1\le\lambda\le r$ is a $(p,q)$-form with values in~$F$,
then
$$\leqalignno{\langle[\ii\Theta(F),\Lambda]u,u\rangle=
&\phantom{{}+{}}\sum_{j,k,\lambda,\mu,J,S}{\,}
c_{jk\lambda\mu}\,u_{J,jS,\lambda}~\ovl{u_{J,kS,\mu}}&(4.8)\cr
&{}+\sum_{j,k,\lambda,\mu,R,K}
c_{jk\lambda\mu}\,u_{kR,K,\lambda}~\ovl{u_{jR,K,\mu}}\cr
&{}-{\,}\sum_{j,\lambda,\mu,J,K}\phantom{k}
c_{jj\lambda\mu}\,u_{J,K,\lambda}~\ovl{u_{J,K,\mu}},\cr}$$
where the sum is extended to all indices $1\le j,k\le n$, $1\le\lambda,
\mu\le r$ and multiindices $|R|=p-1$, $|S|=q-1$ (here the notation
$u_{JK\lambda}$ is extended to non necessarily increasing multiindices by
making it alternate with respect to permutations). It is usually hard to
decide the sign of the curvature term (4.8), except in some special cases.

The easiest case is when $p=n$. Then all terms in the second summation of
(4.8) must have $j=k$ and $R=\{1\ld n\}\ssm\{j\}$, therefore the second
and third summations are equal. It follows that $[\ii\Theta(F),\Lambda]$ is
positive on $(n,q)$-forms under the assumption that $F$ is positive in
the sense of Nakano. In this case $X$ is automatically K\"ahler since
$$\omega=\Tr_F(\ii\Theta(F))=\ii\sum_{j,k,\lambda}c_{jk\lambda\lambda}dz_j
\wedge d\ovl z_k=\ii\Theta(\det F)$$
is a K\"ahler metric.

\begstat{(4.9) Nakano vanishing theorem {\rm (1955)}} Let $X$ be a
compact complex manifold and let $F$ be a Nakano positive vector
bundle on~$X$. Then
$$H^{n,q}(X,F)=H^q(X,K_X\otimes F)=0\qquad\hbox{for every $q\ge 1$.}
\eqno\square$$
\endstat

Another tractable case is the case where $F$ is a line bundle ($r=1$).
Indeed, at each point $x\in X$, we may then choose a coordinate
system which diagonalizes simultaneously the hermitians forms
$\omega(x)$ and $\ii\Theta(F)(x)$, in such a way that
$$\omega(x)=\ii\sum_{1\le j\le n}dz_j\wedge d\ovl z_j,\qquad
\ii\Theta(F)(x)=\ii\sum_{1\le j\le n}\gamma_jdz_j\wedge d\ovl z_j$$
with $\gamma_1\le\ldots\le\gamma_n$. The curvature eigenvalues
$\gamma_j=\gamma_j(x)$ are then uniquely defined and depend continuously
on~$x$. With our previous notation, we have $\gamma_j=c_{jj11}$ and all
other coefficients $c_{jk\lambda\mu}$ are zero. For any $(p,q)$-form
$u=\sum u_{JK}dz_J\wedge d\ovl z_K\otimes e_1$, this gives
$$\leqalignno{\langle[\ii\Theta(F),\Lambda]u,u\rangle
&=\sum_{|J|=p,\,|K|=q}\Big(\sum_{j\in J}\gamma_j+\sum_{j\in K}\gamma_j-
\sum_{1\le j\le n}\gamma_j\Big)|u_{JK}|^2\cr
&\ge(\gamma_1+\ldots+\gamma_q-\gamma_{n-p+1}-\ldots-\gamma_n)|u|^2.
&(4.10)\cr}$$
Assume that $\ii\Theta(F)$ is positive. It is then natural to make the
special choice $\omega=\ii\Theta(F)$ for the K\"ahler metric. Then
$\gamma_j=1$ for $j=1,2\ld n$ and we obtain
$\langle[\ii\Theta(F),\Lambda]u,u\rangle=(p+q-n)|u|^2$.
As a consequence:

\begstat{(4.11) Akizuki-Kodaira-Nakano vanishing theorem {\rm ([AN54])}}
If $F$ is a positive line bundle on a compact complex manifold~$X$, then
$$H^{p,q}(X,F)=H^q(X,\Omega^p_X\otimes F)=0\qquad\hbox{for}~~p+q\ge n+1.
\eqno\square$$
\endstat

More generally, if $F$ is a Griffiths positive (or ample) vector bundle
of rank~$r\ge 1$, Le Potier [LP75] proved that $H^{p,q}(X,F)=0$ for
$p+q\ge n+r$. The proof is not a direct consequence of the Bochner
technique. A rather easy proof has been found by M.~Schneider [Sch74],
using the Leray spectral sequence associated to the projectivized bundle
projection $\bbbp(F)\to X$.

\begstat{(4.12) Exercise} \rm It is important for various applications
to obtain vanishing theorems which are also valid in the case of semipositive
line bundles. The easiest case is the following result of Girbau [Gir76]:
let $(X,\omega)$ be compact K\"ahler; assume that $F$ is a line bundle and
that $\ii\Theta(F)\ge 0$ has at least $n-k$ positive eigenvalues at each
point, for some integer~$k\ge 0$; show that $H^{p,q}(X,F)=0$ for $p+q\ge
n+k+1$.\newline
{\it Hint:}\/ use the K\"ahler metric $\omega_\varepsilon=
\ii\Theta(F)+\varepsilon\omega$ with $\varepsilon>0$ small.
\smallskip
A stronger and more natural ``algebraic version'' of this result has been
obtained by Sommese [Som78]: define $F$ to be $k$-ample if some multiple
$mF$ is such that the canonical map
$$\Phi_{|mF|}:X\ssm B_{|mF|}\to\bbbp^{N-1}$$
has at most $k$-dimensional fibers and $\dim B_{|mF|}\le k$. If
$X$ is projective and $F$ is $k$-ample, show that $H^{p,q}(X,F)=0$ for
$p+q\ge n+k+1$.\newline
{\it Hint:}\/ prove the dual result $H^{p,q}(X,F^{-1})=0$ for $p+q\le n-k-1$
by induction on~$k$. First show that $F$ 0-ample $\Rightarrow$ $F$ positive;
then use hyperplane sections $Y\subset X$ to prove the induction step,
thanks to the exact sequences
$$\eqalignno{
&0\lra\Omega^p_X\otimes F^{-1}\otimes\cO(-Y)\lra\Omega^p_X\otimes F^{-1}
\lra\big(\Omega^p_X\otimes F^{-1}\big)_{\restriction Y}\lra 0,&\cr
&0\lra\Omega^{p-1}_Y\otimes F^{-1}_{\restriction Y}\lra
\big(\Omega^p_X\otimes F^{-1}\big)_{\restriction Y}\lra
\Omega^p_Y\otimes F^{-1}_{\restriction Y}\lra 0.&\square\cr}$$
\endstat

\titleb{5.}{$L^2$ Estimates and Existence Theorems}
The starting point is the following $L^2$ existence theorem, which is
essentially due to H\"ormander [H\"or65, 66], and Andreotti-Vesentini
[AV65]. We will only  outline the main ideas, referring e.g.\ to [Dem82b]
for a detailed exposition of the technical situation considered here.

\begstat{(5.1) Theorem} Let $(X,\omega)$ be a K\"ahler manifold. Here $X$
is not necessarily compact, but we assume that the geodesic distance
$\delta_\omega$ is complete on~$X$. Let $F$ be a hermitian vector bundle of
rank $r$ over $X$, and assume that the curvature operator
$A=A^{p,q}_{F,\omega}=[\ii\Theta(F),\Lambda_\omega]$ is positive
definite everywhere on $\Lambda^{p,q}T^\star_X\otimes F$, $q\ge 1$.
Then for any form $g\in L^2(X,\Lambda^{p,q}T^\star_X\otimes F)$
satisfying $D''g=0$ and $\int_X\langle A^{-1}g,g\rangle
\,dV_\omega<+\infty$, there exists $f\in
L^2(X,\Lambda^{p,q-1}T^\star_X\otimes F)$ such that  $D''f=g$ and
$$\int_X|f|^2\,dV_\omega\le\int_X\langle A^{-1}g,g\rangle\,dV_\omega.$$
\endstat

\begproof{} The assumption that $\delta_\omega$ is complete implies the
existence of cut-off functions $\psi_\nu$ with arbitrarily large compact
support such that $|d\psi_\nu|\le 1$ (take $\psi_\nu$ to be a function of
the distance $x\mapsto \delta_\omega(x_0,x)$, which is an almost
everywhere differentiable 1-Lipschitz function, and regularize if
necessary). From this, it follows that very form $u\in
L^2(X,\Lambda^{p,q}T^\star_X\otimes F)$ such that $D''u\in L^2$ and
$D^{\prime\prime\star}u\in L^2$ in the sense of distribution theory is a
limit of a sequence of smooth forms $u_\nu$ with compact support, in such
a way that $u_\nu\to u$, $D''u_\nu\to D''u$ and $D^{\prime\prime\star}u_\nu
\to D^{\prime\prime\star}u$ in~$L^2$ (just take $u_\nu$ to be a
regularization of $\psi_\nu u$). As a consequence, the basic a priori
inequality (4.7) extends to arbitrary forms $u$ such that
$u,~D''u,D^{\prime\prime\star}u\in L^2$ . Now, consider the Hilbert
space orthogonal decomposition
$$L^2(X,\Lambda^{p,q}T^\star_X\otimes F)=\Ker
D''\oplus({\Ker D''})^\perp,$$
observing that $\Ker D''$ is weakly (hence strongly) closed. Let $v=v_1+v_2$
be the decomposition of a smooth form $v\in\cD^{p,q}(X,F)$ with compact
support according to this decomposition ($v_1$, $v_2$ do not have compact
support in general$\,$!).
Since $(\Ker D'')^\perp\subset\Ker D^{\prime\prime\star}$ by duality
and $g,v_1\in\Ker D''$ by hypothesis, we get
$D^{\prime\prime\star}v_2=0$ and
$$|\langle g,v\rangle|^2=|\langle g,v_1\rangle|^2\le
\int_X\langle A^{-1}g,g\rangle \,dV_\omega
\int_X\langle Av_1,v_1\rangle \,dV_\omega$$
thanks to the Cauchy-Schwarz inequality. The a priori inequality (4.7)
applied to $u=v_1$ yields
$$\int_X\langle Av_1,v_1\rangle \,dV_\omega\le
\|D''v_1\|^2+\|D^{\prime\prime\star}v_1\|^2=\|D^{\prime\prime\star}v_1\|^2=
\|D^{\prime\prime\star}v\|^2.$$
Combining both inequalities, we find
$$|\langle g,v\rangle|^2\le\Big(\int_X\langle A^{-1}g,g\rangle \,dV_\omega\Big)
\|D^{\prime\prime\star}v\|^2$$
for every smooth $(p,q)$-form $v$ with compact support. This shows that we
have a well defined linear form
$$w=D^{\prime\prime\star}v\longmapsto\langle v,g\rangle,\qquad
L^2(X,\Lambda^{p,q-1}T^\star_X\otimes F)
\supset D^{\prime\prime\star}(\cD^{p,q}(F))\lra\bbbc$$
on the range of $D^{\prime\prime\star}$. This linear form
is continuous in $L^2$ norm and has norm $\le C$ with
$$C=\Big(\int_X\langle A^{-1}g,g\rangle \,dV_\omega\Big)^{1/2}.$$
By the Hahn-Banach theorem, there is an element
$f\in L^2(X,\Lambda^{p,q-1}T^\star_X\otimes F)$ with $||f||\le C$, such that
$\langle v,g\rangle=\langle D^{\prime\prime\star}v,f\rangle$ for every $v$,
hence $D''f=g$ in the sense of distributions. The inequality $||f||\le C$
is equivalent to the last estimate in the theorem.\qed
\endproof

The above $L^2$ existence theorem can be applied in the fairly general
context of {\it weakly pseudoconvex} manifolds. By this, we mean a complex
manifold $X$ such that there exists a smooth psh exhaustion function
$\psi$ on $X$ ($\psi$ is said to be an exhaustion if for every $c>0$ the
sublevel set $X_c=\psi^{-1}(c)$ is relatively compact, i.e.\ $\psi(z)$
tends to $+\infty$ when $z$ is taken outside larger and larger compact
subsets of~$X$). In particular, every compact complex manifold $X$ is
weakly pseudoconvex (take $\psi=0$), as well as every Stein
manifold, e.g.\ affine algebraic submanifolds of $\bbbc^N$ (take
$\psi(z)=|z|^2$), open balls $X=B(z_0,r)$ $\big($take $\psi(z)=
1/(r-|z-z_0|^2)\big)$, convex open subsets, etc. Now, a basic observation
is that every weakly pseudoconvex K\"ahler manifold $(X,\omega)$ carries a
{\it complete} K\"ahler metric: let $\psi\ge 0$ be a psh exhaustion function
and set
$$\omega_\varepsilon=\omega+\varepsilon\,\ii d'd''\psi^2=
\omega+2\varepsilon(2\ii\psi d'd''\psi+\ii d'\psi\wedge d''\psi).$$
Then $|d\psi|_{\omega_\varepsilon}\le 1/\varepsilon$ and
$|\psi(x)-\psi(y)|\le\varepsilon^{-1}\delta_{\omega_\varepsilon}(x,y)$.
It follows easily from this estimate that the geodesic balls are relatively
compact, hence $\delta_{\omega_\varepsilon}$ is complete
for every $\varepsilon>0$. Therefore, the $L^2$ existence theorem can be
applied to each K\"ahler metric $\omega_\varepsilon$, and by passing to
the limit it can even be applied to the non necessarily complete
metric~$\omega$. An important special case is the following

\begstat{(5.2) Theorem} Let $(X,\omega)$ be a K\"ahler manifold,
$\dim X=n$. Assume that $X$ is weakly pseudoconvex. Let $F$ be a hermitian
line bundle and let
$$\gamma_1(x)\le\ldots\le\gamma_n(x)$$
be the curvature eigenvalues $($i.e.\ the eigenvalues of $\ii\Theta(F)$
with respect to the metric~$\omega)$ at every point. Assume that the
curvature is positive, i.e.\ $\gamma_1>0$ everywhere. Then for any form
$g\in L^2(X,\Lambda^{n,q}T^\star_X\otimes F)$ satisfying $D''g=0$ and
\hbox{$\int_X\langle(\gamma_1+\ldots+\gamma_q)^{-1}|g|^2\,dV_\omega<+\infty$},
there exists \hbox{$f\in L^2(X,\Lambda^{p,q-1}T^\star_X\otimes F)$} such that
$D''f=g$ and
$$\int_X|f|^2\,dV_\omega\le\int_X
(\gamma_1+\ldots+\gamma_q)^{-1}|g|^2\,dV_\omega.$$
\endstat

\begproof{} Indeed, for $p=n$, Formula~4.10 shows that
$$\langle Au,u\rangle\ge(\gamma_1+\ldots+\gamma_q)|u|^2,$$
hence $\langle A^{-1}u,u\rangle\ge(\gamma_1+\ldots+\gamma_q)^{-1}|u|^2$.\qed
\endproof

An important observation is that the above theorem still applies when the
hermitian metric on $F$ is a singular metric with positive curvature in
the sense of currents. In fact, by standard
regularization techniques (convolution of psh functions by smoothing
kernels), the metric can be made smooth and the solutions obtained by
(5.1) or (5.2) for the smooth metrics have limits satisfying the desired
estimates. Especially, we get the following

\begstat{(5.3) Corollary} Let $(X,\omega)$ be a K\"ahler manifold,
$\dim X=n$. Assume that $X$ is weakly pseudoconvex. Let $F$ be a
holomorphic line bundle equipped with a singular metric whose local
weights are denoted $\varphi\in L^1_\loc$. Suppose that
$$\ii\Theta(F)=2\ii d'd''\varphi\ge\varepsilon\omega$$
for some $\varepsilon>0$. Then for any form $g\in L^2(X,\Lambda^{n,q}
T^\star_X\otimes F)$ satisfying $D''g=0$, there exists
$f\in L^2(X,\Lambda^{p,q-1}T^\star_X\otimes F)$ such that  $D''f=g$ and
$$\int_X|f|^2e^{-2\varphi}\,dV_\omega\le{1\over q\varepsilon}\int_X
|g|^2e^{-2\varphi}\,dV_\omega.\eqno\square$$
\endstat

Here we denoted somewhat incorrectly the metric by $|f|^2e^{-2\varphi}$,
as if the weight $\varphi$ was globally defined on~$X$ (of course, this
is so only if $F$ is globally trivial). We will use this notation anyway,
because it clearly describes the dependence of the $L^2$ norm on the psh
weights.

We now introduce the concept of {\it multiplier ideal sheaf},
following A.~Nadel [Nad89]. The main idea actually goes back to the
fundamental works of Bombieri [Bom70] and H.~Skoda [Sko72a].

\begstat{(5.4) Definition} Let $\varphi$ be a
psh function on an open subset $\Omega\subset X\,$; to $\varphi$ is
associated the ideal subsheaf $\cI(\varphi)\subset\cO_\Omega$ of germs of
holomorphic functions $f\in\cO_{\Omega,x}$ such that $|f|^2e^{-2\varphi}$
is integrable with respect to the Lebesgue measure in some local
coordinates near~$x$.
\endstat

The zero variety $V(\cI(\varphi))$ is thus the set of points in
a neighborhood of which $e^{-2\varphi}$ is non integrable. Of course,
such points occur only if $\varphi$ has logarithmic poles. This is made
precise as follows.

\begstat{(5.5) Definition} A psh function $\varphi$ is said to have a
logarithmic pole of coefficient $\gamma$ at a point $x\in X$ if the Lelong
number
$$\nu(\varphi,x):=\liminf_{z\to x}{\varphi(z)\over\log|z-x|}$$
is non zero and if $\nu(\varphi,x)=\gamma$.
\endstat

\begstat{(5.6) Lemma {\rm (Skoda [Sko72a])}} Let $\varphi$ be a psh function
on an open set $\Omega$ and let $x\in\Omega$.
\smallskip
\item{\rm a)} If $\nu(\varphi,x)<1$, then $e^{-2\varphi}$ is integrable in
a neighborhood of~$x$, in particular $\cI(\varphi)_x=\cO_{\Omega,x}$.
\smallskip
\item{\rm b)} If $\nu(\varphi,x)\ge n+s$ for some integer $s\ge 0$, then
$e^{-2\varphi}\ge C|z-x|^{-2n-2s}$ in a neighborhood of $x$ and
$\cI(\varphi)_x\subset\gm_{\Omega,x}^{s+1}$, where $\gm_{\Omega,x}$ is the
maximal ideal of $\cO_{\Omega,x}$.
\smallskip
\item{\rm c)} The zero variety $V(\cI(\varphi))$ of $\cI(\varphi)$
satisfies
$$E_n(\varphi)\subset V(\cI(\varphi))\subset E_1(\varphi)$$
where $E_c(\varphi)=\{x\in X\,;\,\nu(\varphi,x)\ge c\}$ is the $c$-sublevel
set of Lelong numbers of~$\varphi$.
\vskip0pt
\endstat

\begproof{} a) Set $\Theta=dd^c\varphi$ and $\gamma=\nu(\Theta,x)=
\nu(\varphi,x)$. Let $\chi$ be a cut-off function will support in a small
ball $B(x,r)$, equal to $1$ in $B(x,r/2)$. As $(dd^c\log|z|)^n=\delta_0$,
we get
$$\eqalign{\varphi(z)
&=\int_{B(x,r)}\chi(\zeta)\varphi(\zeta)(dd^c\log|\zeta-z|)^n\cr
&=\int_{B(x,r)}dd^c(\chi(\zeta)\varphi(\zeta))\wedge
\log|\zeta-z|(dd^c\log|\zeta-z|)^{n-1}\cr}$$
for $z\in B(x,r/2)$. Expanding $dd^c(\chi\varphi)$ and observing that
$d\chi=dd^c\chi=0$ on $B(x,r/2)$, we find
$$\varphi(z)=\int_{B(x,r)}\chi(\zeta)\Theta(\zeta)\wedge
\log|\zeta-z|(dd^c\log|\zeta-z|)^{n-1}+\hbox{smooth terms}$$
on $B(x,r/2)$. Fix $r$ so small that
$$\int_{B(x,r)}\chi(\zeta)\Theta(\zeta)\wedge(dd^c\log|\zeta-x|)^{n-1}
\le\nu(\Theta,x,r)<1.$$
By continuity, there exists $\delta,\varepsilon>0$ such that
$$I(z):=\int_{B(x,r)}\chi(\zeta)\Theta(\zeta)\wedge(dd^c\log|\zeta-z|)^{n-1}
\le 1-\delta$$
for all $z\in B(x,\varepsilon)$. Applying Jensen's convexity inequality to
the probability measure
$$d\mu_z(\zeta)=I(z)^{-1}\chi(\zeta)\Theta(\zeta)\wedge
(dd^c\log|\zeta-z|)^{n-1},
$$
we find
$$\eqalign{
-\varphi(z)&=\int_{B(x,r)}I(z)\log|\zeta-z|^{-1}\,d\mu_z(\zeta)+O(1)\quad
\Longrightarrow\cr
e^{-2\varphi(z)}&\le C\int_{B(x,r)}|\zeta-z|^{-2I(z)}\,d\mu_z(\zeta).\cr}$$
As
$$d\mu_z(\zeta)\le C_1|\zeta-z|^{-(2n-2)}\Theta(\zeta)\wedge
(dd^c|\zeta|^2)^{n-1}=C_2|\zeta-z|^{-(2n-2)}d\sigma_\Theta(\zeta),$$
we get
$$e^{-2\varphi(z)}\le C_3\int_{B(x,r)}|\zeta-z|^{-2(1-\delta)-(2n-2)}
d\sigma_\Theta(\zeta),$$
and the Fubini theorem implies that $e^{-2\varphi(z)}$ is integrable on a
neighborhood of~$x$.
\smallskip

\noindent b) If $\nu(\varphi,x)=\gamma$, the convexity properties of
psh functions, namely, the convexity of $\log r\mapsto\sup_{|z-x|=r}
\varphi(z)$ implies that
$$\varphi(z)\le\gamma\log|z-x|/r_0+M,$$
where $M$ is the supremum on $B(x,r_0)$. Hence there exists a constant
$C>0$ such that $e^{-2\varphi(z)}\ge C|z-x|^{-2\gamma}$ in a neighborhood
of $x$. The desired result follows from the identity
$$\int_{B(0,r_0)}{\big|\sum a_\alpha z^\alpha\big|^2\over|z|^{2\gamma}}dV(z)=
{\rm Const}\int_0^{r_0}\Big(\sum |a_\alpha|^2r^{2|\alpha|}\Big)
r^{2n-1-2\gamma}\,dr,$$
which is an easy consequence of Parseval's formula. In fact, if $\gamma$
has integral part $[\gamma]=n+s$, the integral converges if and only if
$a_\alpha=0$ for $|\alpha|\le s$.
\smallskip
\noindent c) is just a simple formal consequence of a) and b).\qed
\endproof

\begstat{(5.7) Proposition {\rm([Nad89])}} For any psh function
$\varphi$ on~$\Omega\subset X$, the sheaf
$\cI(\varphi)$ is a coherent sheaf of ideals over $\Omega$.
\endstat

\begproof{} Since the result is local, we may assume that $\Omega$ is the
unit ball in $\bbbc^n$. Let $E$ be the set of all holomorphic functions $f$
on $\Omega$ such that $\int_\Omega|f|^2e^{-2\varphi}\,d\lambda<+\infty$.
By the strong noetherian property of coherent sheaves, the set $E$
generates a coherent ideal sheaf $\cJ\subset\cO_\Omega$. It is clear that
$\cJ\subset\cI(\varphi)$; in order to prove the equality,
we need only check that $\cJ_x+\cI(\varphi)_x\cap
\gm^{s+1}_{\Omega,x}=\cI(\varphi)_x$ for every integer $s$, in
view of the Krull lemma. Let $f\in\cI(\varphi)_x$ be defined
in a neighborhood $V$ of $x$ and let $\theta$ be a cut-off function
with support in $V$ such that $\theta=1$ in a neighborhood of $x$.
We solve the equation $d''u=g:=d''(\theta f)$ by means of H\"ormander's
$L^2$ estimates~5.3, where $F$ is the trivial line bundle
$\Omega\times\bbbc$ equipped with the strictly psh weight
$$\wt\varphi(z)=\varphi(z)+(n+s)\log|z-x|+|z|^2.$$
We get a solution $u$ such that
$\int_\Omega |u|^2e^{-2\varphi}|z-x|^{-2(n+s)}d\lambda<\infty$, thus
$F=\theta f-u$ is holomorphic, $F\in E$ and
$f_x-F_x=u_x\in\cI(\varphi)_x\cap\gm_{\Omega,x}^{s+1}$.
This proves our contention.\qed
\endproof

The multiplier ideal sheaves satisfy the following basic
fonctoriality property with respect to direct images of sheaves
by modifications.

\begstat{(5.8) Proposition} Let $\mu:X'\to X$ be a modification
of non singular complex manifolds $($i.e.\ a proper generically 1:1
holomorphic map$)$, and let $\varphi$ be a psh function on~$X$. Then
$$\mu_\star\big(\cO(K_{X'})\otimes\cI(\varphi\circ\mu)\big)=
\cO(K_X)\otimes\cI(\varphi).$$
\endstat

\begproof{} Let $n=\dim X=\dim X'$ and let $S\subset X$ be an analytic
set such that $\mu:X'\ssm S'\to X\ssm S$ is a biholomorphism. By
definition of multiplier ideal sheaves, $\cO(K_X)\otimes\cI(\varphi)$ is
just the sheaf of holomorphic $n$-forms $f$ on open sets $U\subset X$
such that
$\ii^{n^2}f\wedge\ovl f\,e^{-2\varphi}\in L^1_\loc(U)$. Since
$\varphi$ is locally bounded from above, we may even consider forms $f$
which are a priori defined only on $U\ssm S$, because $f$ will be
in $L^2_\loc(U)$ and therefore will automatically extend through~$S$.
The change of variable formula yields
$$\int_U\ii^{n^2}f\wedge\ovl f\,e^{-2\varphi}=
\int_{\mu^{-1}(U)}\ii^{n^2}\mu^\star f\wedge\ovl{\mu^\star f}
\,e^{-2\varphi\circ\mu},$$
hence $f\in\Gamma(U,\cO(K_X)\otimes\cI(\varphi))$ iff
$\mu^\star f\in\Gamma(\mu^{-1}(U),\cO(K_{X'})\otimes\cI(\varphi\circ
\mu))$. Proposition~5.8 is proved.\qed
\endproof

\begstat{(5.9) Remark} \rm If $\varphi$ has analytic singularities
(according to Definition~1.10), the computation of $\cI(\varphi)$
can be reduced to a purely algebraic problem.

The first observation is that $\cI(\varphi)$ can be computed easily if
$\varphi$ has the form $\varphi=\sum\alpha_j\log|g_j|$ where
$D_j=g^{-1}_j(0)$ are nonsingular irreducible divisors with normal
crossings. Then $\cI(\varphi)$ is the sheaf of functions $h$ on open
sets $U\subset X$ such that
$$\int_U|h|^2\prod|g_j|^{-2\alpha_j}dV<+\infty.$$
Since locally the $g_j$ can be taken to be coordinate functions from a
local coordinate system $(z_1\ld z_n)$, the condition is that
$h$ is divisible by $\prod g_j^{m_j}$ where $m_j-\alpha_j>-1$ for
each $j$, i.e.\ $m_j\ge\lfloor\alpha_j\rfloor$ (integer part). Hence
$$\cI(\varphi)=\cO(-\lfloor D\rfloor)=\cO(-\sum\lfloor\alpha_j\rfloor D_j)$$
where $\lfloor D\rfloor$ denotes the integral part of the $\bbbq$-divisor
$D=\sum\alpha_jD_j$.

Now, consider the general case of analytic singularities and suppose
that $\varphi\sim{\alpha\over 2}\log\big(|f_1|^2+\cdots+|f_N|^2\big)$
near the poles. By the remarks after Definition~1.10, we may assume
that the $(f_j)$ are generators of the integrally closed ideal
sheaf $\cJ=\cJ(\varphi/\alpha)$, defined as the sheaf of holomorphic
functions $h$ such that $|h|\le C\exp(\varphi/\alpha)$. In this case,
the computation is made as follows (see also L.~Bonavero's work [Bon93],
where similar ideas are used in connection with ``singular'' holomorphic
Morse inequalities).

First, one computes a smooth modification $\mu:\wt X\to X$ of $X$ such
that $\mu^\star\cJ$ is an invertible sheaf $\cO(-D)$ associated with
a normal crossing divisor $D=\sum\lambda_jD_j$, where $(D_j)$ are the
components of the exceptional divisor of $\wt X$ (take the blow-up $X'$
of $X$ with respect to the ideal $\cJ$ so that the pull-back of $\cJ$
to $X'$ becomes an invertible sheaf $\cO(-D')$, then blow up again by
Hironaka [Hir64] to make $X'$ smooth and $D'$ have normal crossings).
Now, we have $K_{\wt X}=\mu^\star K_X+R$ where $R=\sum\rho_j
D_j$ is the zero divisor of the Jacobian function $J_\mu$ of the blow-up
map. By the direct image formula~5.8, we get
$$\cI(\varphi)=\mu_\star\big(\cO(K_{\wt X}-\mu^\star K_X)\otimes
\cI(\varphi\circ\mu)\big)=\mu_\star\big(\cO(R)\otimes
\cI(\varphi\circ\mu)\big).$$
Now, $(f_j\circ\mu)$ are generators of the ideal $\cO(-D)$, hence
$$\varphi\circ\mu\sim\alpha\sum\lambda_j\log|g_j|$$
where $g_j$ are local generators of $\cO(-D_j)$. We are thus reduced to
computing multiplier ideal sheaves in the case where the poles are
given by a $\bbbq$-divisor with normal crossings $\sum\alpha\lambda_jD_j$.
We obtain $\cI(\varphi\circ\mu)=\cO(-\sum\lfloor\alpha\lambda_j\rfloor D_j)$,
hence
$$\cI(\varphi)=\mu_\star\cO_{\wt X}\big(\sum(\rho_j-\lfloor\alpha\lambda_j
\rfloor)D_j\big).\eqno\square$$
\endstat

\begstat{(5.10) Exercise} \rm Compute the multiplier ideal sheaf
$\cI(\varphi)$ associated with
$\varphi=\log(|z_1|^{\alpha_1}+\ldots+|z_p|^{\alpha_p})$
for arbitrary real numbers $\alpha_j>0$.\newline
{\it Hint:} using Parseval's formula and polar coordinates $z_j=r_j
e^{\ii\theta_j}$, show that the problem is equivalent to determining for
which $p$-tuples $(\beta_1\ld\beta_p)\in\bbbn^p$ the integral
$$\int_{[0,1]^p}{r_1^{2\beta_1}\ldots r_p^{2\beta_p}\,r_1dr_1\ldots r_pdr_p
\over r_1^{2\alpha_1}+\ldots+r_p^{2\alpha_p}}
=\int_{[0,1]^p}{t_1^{(\beta_1+1)/\alpha_1}\ldots t_p^{(\beta_p+1)/\alpha_p}
\over t_1+\ldots+t_p}{dt_1\over t_1}\ldots{dt_p\over t_p}$$
is convergent. Conclude from this that $\cI(\varphi)$ is generated by
the monomials $z_1^{\beta_1}\ldots z_p^{\beta_p}$ such that
$\sum(\beta_p+1)/\alpha_p>1$. (This exercise shows that the analytic
definition of $\cI(\varphi)$ is sometimes also quite convenient for
computations).\qed
\endstat

Let $F$ be a line bundle over $X$ with a singular metric $h$ of
curvature current $\Theta_h(F)$. If $\varphi$ is the weight
representing the metric in an open set $\Omega\subset X$, the ideal
sheaf $\cI(\varphi)$ is independent of the choice of the trivialization
and so it is the restriction to $\Omega$ of a global coherent sheaf
$\cI(h)$ on~$X$. We will sometimes still write $\cI(h)=\cI(\varphi)$ by
abuse of notation. In this context, we have the following fundamental
vanishing theorem, which is probably one of the most central results of
analytic and algebraic geometry (as we will see later, it contains the
Kawamata-Viehweg vanishing theorem as a special case).

\begstat{(5.11) Nadel vanishing theorem {\rm([Nad89], [Dem93b])}} Let
$(X,\omega)$ be a K\"ahler weakly pseudoconvex manifold, and let $F$ be a
holomorphic line bundle over $X$ equipped with a singular hermitian metric
$h$ of weight~$\varphi$. Assume that $\ii\Theta_h(F)\ge\varepsilon\omega$
for some continuous positive function $\varepsilon$ on~$X$. Then
$$H^q\big(X,\cO(K_X+F)\otimes\cI(h)\big)=0\qquad
\hbox{for all $q\ge 1$.}$$
\endstat

\begproof{} Let $\cL^q$ be the sheaf of germs of $(n,q)$-forms
$u$ with values in $F$ and with measurable coefficients, such that
both $|u|^2e^{-2\varphi}$ and $|d''u|^2e^{-2\varphi}$ are locally
integrable. The $d''$ operator defines a complex of sheaves
$(\cL^\bu,d'')$ which is a resolution of the sheaf
$\cO(K_X+F)\otimes\cI(\varphi)$: indeed, the kernel of $d''$ in degree
$0$ consists of all germs of holomorphic $n$-forms with values in $F$
which satisfy the integrability condition; hence the coefficient
function lies in $\cI(\varphi)$; the exactness in degree $q\ge 1$
follows from Corollary~5.3 applied on arbitrary small balls.
Each sheaf $\cL^q$ is a ${\cal C}^\infty$-module, so $\cL^\bu$
is a resolution by acyclic sheaves. Let $\psi$ be a smooth psh
exhaustion function on~$X$. Let us apply Corollary~5.3 globally on~$X$,
with the original metric of $F$ multiplied by the factor
$e^{-\chi\circ\psi}$, where $\chi$ is a convex increasing function of
arbitrary fast growth at infinity. This factor can be used to ensure
the convergence of integrals at infinity. By Corollary~5.3, we conclude
that $H^q\big(\Gamma(X,\cL^\bu)\big)=0$ for $q\ge 1$. The theorem
follows.\qed
\endproof

\begstat{(5.12) Corollary} Let $(X,\omega)$, $F$ and $\varphi$ be as in
Theorem~5.11 and let $x_1\ld x_N$ be isolated points in the zero
variety~$V(\cI(\varphi))$. Then there is a surjective map
$$H^0(X,K_X+F)\lraww\bigoplus_{1\le j\le N}\cO(K_X+L)_{x_j}\otimes
\big(\cO_X/\cI(\varphi)\big)_{x_j}.$$
\endstat

\begproof{} Consider the long exact sequence of cohomology associated to the
short exact sequence $0\to\cI(\varphi)\to\cO_X\to\cO_X/\cI(\varphi)\to 0$
twisted by $\cO(K_X+F)$, and apply Theorem~5.11 to obtain the vanishing of
the first $H^1$ group. The asserted surjectivity property follows.\qed
\endproof

\begstat{(5.13) Corollary} Let $(X,\omega)$, $F$ and $\varphi$ be as in
Theorem~5.11 and suppose that the weight function $\varphi$ is such that
$\nu(\varphi,x)\ge n+s$ at some point $x\in X$ which is an isolated point
of $E_1(\varphi)$. Then $H^0(X,K_X+F)$ generates all $s$-jets at~$x$.
\endstat

\begproof{} The assumption is that $\nu(\varphi,y)<1$ for $y$ near $x$,
$y\ne x$. By Skoda's lemma 5.6~b), we conclude that $e^{-2\varphi}$ is
integrable at all such points~$y$, hence $\cI(\varphi)_y=\cO_{X,y}$, whilst
$\cI(\varphi)_x\subset\gm_{X,x}^{s+1}$ by 5.6~a). Corollary 5.13 is thus
a special case of 5.12.\qed
\endproof

The philosophy of these results (which can be seen as generalizations
of the H\"ormander-Bombieri-Skoda theorem [Bom70], [Sko72a,~75]) is that
the problem of constructing holomorphic sections of $K_X+F$ can be
solved by constructing suitable hermitian metrics on $F$ such that the
weight $\varphi$ has isolated poles at given points~$x_j$.

\begstat{(5.14) Exercise} \rm Assume that $X$ is compact and that $L$ is
a positive line bundle on~$X$. Let $\{x_1\ld x_N\}$ be a finite set.
Show that there are constants $a,b\ge 0$ depending only on $L$ and $N$
such that $H^0(X,mL)$ generates jets of any order $s$ at all points $x_j$
for $m\ge as+b$.\newline
{\it Hint:} Apply Corollary~5.12 to $F=-K_X+mL$, with a singular
metric on $L$ of the form $h=h_0e^{-\varepsilon\psi}$, where $h_0$ is
smooth of positive curvature, $\varepsilon>0$ small and
$\psi(z)\sim\log|z-x_j|$ in a neighborhood of $x_j$.\newline
Derive the Kodaira embedding theorem from the above result:
\endstat

\begstat{(5.15) Theorem {\rm(Kodaira)}} If $L$ is a line bundle on a compact
complex manifold, then $L$ is ample if and only if $L$ is positive.\qed
\endstat

\begstat{(5.16) Exercise {\rm (solution of the Levi problem)}} \rm Show
that the following two properties are equivalent.
\smallskip
\item{\rm a)} $X$ is strongly pseudoconvex, i.e.\ $X$ admits a strongly psh
exhaustion function.
\smallskip
\item{\rm b)} $X$ is Stein, i.e.\ the global holomorphic functions
$H^0(X,\cO_X)$ separate points and yield local coordinates at any point,
and $X$ is holomorphically convex $($this means that for any discrete
sequence $z_\nu$ there is a function $f\in H^0(X,\cO_X)$ such that
$|f(z_\nu)|\to\infty)$.\qed
\vskip0pt
\endstat

\begstat{(5.17) Remark} \rm As long as forms of bidegree $(n,q)$ are
considered, the $L^2$ estimates can be extended to complex spaces with
arbitrary singularities. In fact, if $X$ is a complex space and $\varphi$
is a psh weight function on~$X$, we may still define a sheaf $K_X(\varphi)$
on~$X$, such that the sections on an open set $U$ are the holomorphic
$n$-forms $f$ on the regular part $U\cap X_\reg$, satisfying the
integrability condition $\ii^{n^2}f\wedge\ovl f\,e^{-2\varphi}\in
L^1_\loc(U)$. In this setting, the fonctoriality property 5.8 becomes
$$\mu_\star\big(K_{X'}(\varphi\circ\mu)\big)=K_X(\varphi)$$
for arbitrary complex spaces $X$, $X'$ such that $\mu:X'\to X$
is a modification. If $X$ is nonsingular we have $K_X(\varphi)=
\cO(K_X)\otimes\cI(\varphi)$, however, if $X$ is
singular, the symbols $K_X$ and $\cI(\varphi)$ must not be
dissociated. The statement of the Nadel vanishing theorem becomes
$H^q(X,\cO(F)\otimes K_X(\varphi))=0$ for $q\ge 1$, under the same
assumptions ($X$ K\"ahler and weakly pseudoconvex, curvature
$\ge\varepsilon\omega$). The proof can be obtained by restricting
everything to $X_\reg$. Although in general $X_\reg$ is not weakly
pseudoconvex (e.g.\ in case $\codim X_\sing\ge 2$), $X_\reg$ is
always K\"ahler complete (the complement of a proper analytic subset in
a K\"ahler weakly pseudoconvex space is complete K\"ahler, see
e.g.\ [Dem82a]). As a consequence, Nadel's vanishing theorem is
essentially insensitive to the presence of singularities.\qed
\endstat

\titleb{6.}{Numerically Effective Line Bundles}
Many problems of algebraic geometry (e.g.\ problems of classification of
algebraic surfaces or higher dimensional varieties) lead in a natural
way to the study of line bundles satisfying semipositivity conditions.
It turns out that semipositivity in the sense of curvature (at least,
as far as smooth metrics are considered) is not a very satisfactory
notion. A more flexible notion perfectly suitable for algebraic purposes
is the notion of {\it numerical effectivity}. The goal of this section
is to give a few fundamental algebraic definitions and to discuss their
differential geometric counterparts. We first suppose that $X$ is a
projective algebraic manifold, $\dim X=n$.

\begstat{(6.1) Definition} A holomorphic line bundle $L$ over a projective
manifold $X$ is said to be numerically effective, nef for short, if
$L\cdot C=\int_C c_1(L)\ge 0$ for every curve $C\subset X$.
\endstat

If $L$ is nef, it can be shown that $L^p\cdot Y=\int_Y c_1(L)^p\ge 0$
for any $p$-dimensional subvariety $Y\subset X$ (see e.g.\ [Har70]).
In relation with this, let us recall the Nakai-Moishezon ampleness
criterion: a line bundle $L$ is ample if and only if $L^p\cdot Y>0$ for
every $p$-dimensional subvariety~$Y$. From this, we easily infer

\begstat{(6.2) Proposition} Let $L$ be a line bundle on a projective
algebraic manifold~$X$, on which an ample line bundle $A$ and a hermitian
metric $\omega$ are given. The following properties are equivalent:
\smallskip
\item{\rm a)} $L$ is nef$\,;$
\smallskip
\item{\rm b)} for any integer $k\ge 1$, the line bundle $kL+A$ is ample$\,;$
\smallskip
\item{\rm c)} for every $\varepsilon>0$, there is a smooth metric
$h_\varepsilon$ on $L$ such that $\ii\Theta_{h_\varepsilon}(L)\ge
-\varepsilon\omega$.
\vskip0pt
\endstat

\begproof{} a) $\Rightarrow$ b). If $L$ is nef and $A$ is ample then
clearly $kL+A$ satisfies the Nakai-Moishezon criterion, hence
$kL+A$ is ample.
\smallskip
\noindent b) $\Rightarrow$ c). Condition c) is independent of the choice
of the hermitian metric, so we may select a metric $h_A$ on $A$ with
positive curvature and set $\omega=\ii\Theta(A)$. If $kL+A$ is ample,
this bundle has a metric $h_{kL+A}$ of positive curvature. Then the
metric $h_L=(h_{kL+A}\otimes h_A^{-1})^{1/k}$ has curvature
$$\ii\Theta(L)={1\over k}\big(\ii\Theta(kL+A)-\ii\Theta(A)\big)\ge
-{1\over k}\ii\Theta(A)\,;$$
in this way the negative part can be made smaller than $\varepsilon\,\omega$
by taking $k$ large enough.
\smallskip
\noindent c) $\Rightarrow$ a). Under hypothesis c), we get
$L\cdot C=\int_C{\ii\over2\pi}\Theta_{h_\varepsilon}(L)\ge-
{\varepsilon\over2\pi}\int_C\omega$ for every curve $C$ and every
$\varepsilon>0$, hence $L\cdot C\ge 0$ and $L$ is nef.\qed
\endproof

Let now $X$ be an arbitrary compact complex manifold. Since there need
not exist any curve in $X$, Property 6.2~c) is simply taken as a definition
of nefness ([DPS94]):

\begstat{(6.3) Definition} A line bundle $L$ on a compact complex
manifold $X$ is said to be nef if for every $\varepsilon>0$, there is a
smooth hermitian metric $h_\varepsilon$ on $L$ such that
$\ii\Theta_{h_\varepsilon}(L)\ge -\varepsilon\omega$.
\endstat

In general, it is not possible to extract a smooth limit $h_0$ such that
$\ii\Theta_{h_0}(L)\ge 0$. The following simple example is given in
[DPS94] (Example~1.7). Let $E$ be a non trivial extension
$0\to\cO\to E\to\cO\to 0$ over an elliptic curve $C$ and let $X=P(E)$
be the corresponding ruled surface over $C$. Then $L=\cO_{P(E)}(1)$ is
nef but does not admit any smooth metric of nonnegative curvature.
This example answers negatively a question raised by [Fuj83].

Let us now introduce the important concept of {\it Kodaira-Iitaka dimension}
of a line bundle.

\begstat{(6.4) Definition} If $L$ is a line bundle, the Kodaira-Iitaka
dimension $\kappa(L)$ is the supremum of the rank of the canonical maps
$$\Phi_m:X\ssm B_m\lra P(V_m^\star),~~~~x\longmapsto H_x=\{\sigma\in
V_m\,;\,\sigma(x)=0\},~~~~m\ge 1$$
with $V_m=H^0(X,mL)$ and $B_m=\bigcap_{\sigma\in V_m}\sigma^{-1}(0)={}$base
locus of $V_m$. In case $V_m=\{0\}$ for all $m\ge 1$, we set $\kappa(L)=
-\infty$.\newline
A line bundle is said to be big if $\kappa(L)=\dim X$.
\endstat

The following lemma is well-known (the proof is a rather elementary
consequence of the Schwarz lemma).

\begstat{(6.5) Serre-Siegel lemma {\rm([Ser54], [Sie55])}} Let $L$ be any
line bundle on a compact complex manifold. Then we have
$$h^0(X,mL)\le O(m^{\kappa(L)})\qquad\hbox{for $m\ge 1$,}$$
and $\kappa(L)$ is the smallest constant for which this estimate holds.\qed
\endstat

We now discuss the various concepts of positive cones in the space of
numerical classes of line bundles, and establish a simple dictionary
relating these concepts to corresponding concepts in the context of
differential geometry.

Let us recall that an integral cohomology class in
$H^2(X,\bbbz)$ is the first Chern class of a holomorphic (or algebraic)
line bundle if and only if it lies in the
{\it Neron-Severi} group
$$\NS(X)=\Ker\big(H^2(X,\bbbz)\to H^2(X,\cO_X)\big)$$
(this fact is just an elementary consequence of the exponential exact
sequence $0\to\bbbz\to\cO\to\cO^\star\to 0$). If $X$ is compact
K\"ahler, as we will suppose from now on in this section,
this is the same as saying that the class is of type $(1,1)$
with respect to Hodge decomposition.

Let $\NS_\bbbr(X)$ be the real vector space $\NS(X)\otimes\bbbr\subset
H^2(X,\bbbr)$. We define four convex cones
$$\eqalign{
&N_\amp(X)\subset N_\eff(X)\subset \NS_\bbbr(X),\cr
&N_\nef(X)\subset N_\psef(X)\subset\NS_\bbbr(X)\cr}$$
which are, respectively, the {\it convex cones} generated by Chern classes
$c_1(L)$ of ample and effective line bundles, resp.\ the {\it closure of
the convex cones} generated by numerically effective and pseudo-effective
line bundles; we say that $L$ is effective if $mL$ has a section for some
$m>0$, i.e.\ if $\cO(mL)\simeq\cO(D)$ for some effective divisor $D\,$; and
we say that $L$ pseudo-effective if $c_1(L)$ is the cohomology class of some
closed positive current~$T$, i.e.\ if $L$ can be equipped with a singular
hermitian metric $h$ with $T={\ii\over 2\pi}\Theta_h(L)\ge 0$ as a current.
For each of the ample, effective, nef and pseudo-effective
cones, the first Chern class $c_1(L)$ of a line bundle $L$ lies in the
cone if and only if $L$ has the corresponding property (for $N_\psef$
use the fact that the space of positive currents of mass $1$ is
weakly compact; the case of all other cones is obvious).

\begstat{(6.6) Proposition} Let $(X,\omega)$ be a compact K\"ahler manifold.
The numerical cones satisfy the following properties.
\smallskip
\item{\rm a)} $N_\amp=N_\amp^\circ\subset N_\nef^\circ$,~~
$N_\nef\subset N_\psef$.
\smallskip
\item{\rm b)} If moreover $X$ is projective algebraic, we have
$N_\amp=N_\nef^\circ$ $($therefore $\ovl N_\amp=N_\nef)$,
and $\ovl N_\eff=N_\psef$.
\smallskip\noindent
If $L$ is a line bundle on $X$ and $h$ denotes a hermitian metric on $L$,
the following properties are equivalent:
\smallskip
\item{\rm c)} $c_1(L)\in N_\amp$
$\Leftrightarrow$ $\exists\varepsilon>0$, $\exists h$ smooth such that
$\ii\Theta_h(L)\ge\varepsilon\omega$.
\smallskip
\item{\rm d)} $c_1(L)\in N_\nef$
$\Leftrightarrow$ $\forall\varepsilon>0$, $\exists h$ smooth such that
$\ii\Theta_h(L)\ge-\varepsilon\omega$.
\smallskip
\item{\rm e)} $c_1(L)\in N_\psef$
$\Leftrightarrow$ $\exists h$ possibly singular such that
$\ii\Theta_h(L)\ge 0$.
\smallskip
\item{\rm f)} If moreover $X$ is projective algebraic, then\newline
$c_1(L)\in N_\eff^\circ$~$\Leftrightarrow\kappa(L)=\dim X$\newline
\phantom{$c_1(L)\in N_\eff^\circ$}~$\Leftrightarrow\exists\varepsilon>0$,
$\exists h$ possibly singular such that
$\ii\Theta_h(L)\ge\varepsilon\omega$.
\vskip0pt
\endstat

\begproof{} c) and d) are already known and e) is a definition.
\smallskip\noindent
a) The ample cone $N_\amp$ is always open by definition and contained
in $N_\nef$, so the first inclusion is obvious ($N_\amp$ is
of course empty if $X$ is not projective algebraic). Let us now prove
that $N_\nef\subset N_\psef$. Let $L$ be a line bundle with
$c_1(L)\in N_\nef$. Then for every $\varepsilon>0$, there is
a current $T_\varepsilon={\ii\over 2\pi}\Theta_{h_\varepsilon}(L)\ge
-\varepsilon\omega$. Then $T_\varepsilon+\varepsilon\omega$ is
a closed positive current and the family is uniformly bounded
in mass for $\varepsilon\in{}]0,1]$, since
$$\int_X(T_\varepsilon+\varepsilon\omega)\wedge\omega^{n-1}=
\int_X c_1(L)\wedge\omega^{n-1}+\varepsilon\int_X\omega^n.$$
By weak compactness, some subsequence converges to a weak limit
$T\ge 0$ and $T\in c_1(L)$ (the cohomology class $\{T\}$ of a current
is easily shown to depend continuously on $T$ with respect to the weak
topology; use e.g.\ Poincar\'e duality to check this).
\smallskip\noindent
b) If $X$ is projective, the equality $N_\amp=N_\nef^\circ$
is a simple consequence of 6.2~b) and of the fact that ampleness (or
positivity) is an open property. It remains to show that
$N_\psef\subset\ovl N_\eff$. Let $L$ be a line bundle with
$c_1(L)\in N_\psef$ and let $h_L$ be a singular hermitian on $L$
such that $T={\ii\over 2\pi}\Theta(L)\ge 0$. Fix a point $x_0\in X$
such that the Lelong number of $T$ at $x_0$ is zero, and take a sufficiently
positive line bundle $A$ (replacing $A$ by a multiple if necessary), such
that $A-K_X$ has a singular metric $h_{A-K_X}$ of curvature
$\ge\varepsilon\omega$ and such that $h_{A-K_X}$ is smooth on
$X\ssm\{x_0\}$ and has an isolated logarithmic pole of Lelong number $\ge n$
at~$x_0$. Then apply Corollary 5.13 to $F=mL+A-K_X$ equipped with
the metric $h_L^{\otimes m}\otimes h_{A-K_X}$. Since the weight
$\varphi$ of this metric has a Lelong number $\ge n$ at $x_0$ and
a Lelong number equal to the Lelong number of $T={\ii\over2\pi}\Theta(L)$
at nearby points, $\limsup_{x\to x_0}\nu(T,x)=\nu(T,x_0)=0$,
Corollary 5.13 implies that $H^0(X,K_X+F)=H^0(X,mL+A)$ has a section
which does not vanish at~$x_0$. Hence there is an effective divisor $D_m$
such that $\cO(mL+A)=\cO(D_m)$ and $c_1(L)={1\over m}\{D_m\}-{1\over m}
c_1(A)=\lim{1\over m}\{D_m\}$ is in $\ovl N_\eff$.\qed
\smallskip

\noindent
f) Fix a nonsingular ample divisor $A$. If $c_1(L)\in N_\eff^\circ$, there
is an integer $m>0$ such that $c_1(L)-{1\over m}c_1(A)$ is still effective,
hence for $m,p$ large we have $mpL-pA=D+F$ with an effective divisor $D$
and a numerically trivial line bundle $F$. This implies
$\cO(kmpL)=\cO(kpA+kD+kF)\supset\cO(kpA+kF)$, hence $h^0(X,kmpL)\ge h^0(X,kpA+
kF)\sim (kp)^nA^n/n!$ by the Riemann-Roch formula. Therefore $\kappa(L)=n$.

If $\kappa(L)=n$, then $h^0(X,kL)\ge ck^n$ for $k\ge k_0$ and $c>0$.
The exact cohomology sequence
$$0\lra H^0(X,kL-A)\lra H^0(X,kL)\lra H^0(A,kL_{\restriction A})$$
where $h^0(A,kL_{\restriction A})=O(k^{n-1})$ shows that
$kL-A$ has non zero sections for $k$ large. If $D$ is the divisor of
such a section, then $kL\simeq\cO(A+D)$. Select a smooth metric
on $A$ such that ${\ii\over2\pi}\Theta(A)\ge\varepsilon_0\omega$ for
some $\varepsilon_0>0$, and take the singular metric
on $\cO(D)$ with weight function $\varphi_D=\sum\alpha_j\log|g_j|$
described in Example~3.13. Then the metric with weight $\varphi_L=
{1\over k}(\varphi_A+\varphi_D)$ on $L$ yields
$${\ii\over 2\pi}\Theta(L)={1\over k}
\Big({\ii\over2\pi}\Theta(A)+[D]\Big)\ge(\varepsilon_0/k)\,\omega,$$
as desired.

Finally, the curvature condition $\ii\Theta_h(L)\ge\varepsilon\omega$ in the
sense of currents yields by definition $c_1(L)\in N_\psef^\circ$. Moreover,
b) implies $N_\psef^\circ=N_\eff^\circ$.\qed
\endproof

Before going further, we need a lemma.

\begstat{(6.7) Lemma} Let $X$ be a compact K\"ahler $n$-dimensional
manifold, let $L$ be a nef line bundle on~$X$, and let $E$ be an
arbitrary holomorphic vector bundle. Then
$h^q(X,cO(E)\otimes\cO(kL))=o(k^n)$
as $k\to+\infty$, for every $q\ge1$. If $X$ is projective algebraic, the
following more precise bound holds:
$$h^q(X,\cO(E)\otimes\cO(kL))=O(k^{n-q}),\qquad\forall q\ge0.$$
\endstat

\begproof{} The K\"ahler case will be proved in Section~12, as a
consequence of the holomorphic Morse inequalities. In the projective
algebraic case, we proceed by induction on $n=\dim X$. If $n=1$ the result
is clear, as well as if $q=0$. Now let $A$ be a nonsingular ample divisor
such that $E\otimes\cO(A-K_X)$ is Nakano positive. Then the Nakano
vanishing theorem applied to the vector bundle $F=E\otimes\cO(kL+A-K_X)$
shows that $H^q(X,\cO(E)\otimes\cO(kL+A))=0$ for all $q\ge 1$.
The exact sequence
$$0\to\cO(kL)\to\cO(kL+A)\to\cO(kL+A)_{\restriction A}\to0$$
twisted by $E$ implies
$$H^q(X,\cO(E)\otimes\cO(kL))\simeq H^{q-1}(A,\cO(E_{\restriction A}
\otimes\cO(kL+A)_{\restriction A}),$$
and we easily conclude by induction since $\dim A=n-1$. Observe that the
argument does not work any more if $X$ is not algebraic. It seems to be
unknown whether the $O(k^{n-q})$ bound still holds in that case.\qed
\endproof

\begstat{(6.8) Corollary} If $L$ is nef, then $L$ is big $($i.e.\
$\kappa(L)=n)$ if and only if $L^n>0$. Moreover, if $L$ is nef and big,
then for every $\delta>0$, $L$ has a singular metric $h=e^{-2\varphi}$
such that $\,\max_{x\in X}\nu(\varphi,x)\le\delta$ and
$\ii\Theta_h(L)\ge\varepsilon\,\omega$ for some $\varepsilon>0$.
The metric $h$ can be chosen to be smooth on the complement of a fixed
divisor $D$, with logarithmic poles along $D$.
\endstat

\begproof{} By Lemma~6.7 and the Riemann-Roch formula, we have
$h^0(X,kL)=\chi(X,kL)+o(k^n)=k^nL^n/n!+o(k^n)$, whence the first
statement. If $L$ is big, the proof made in 6.5~f) shows that
there is a singular metric $h_1$ on $L$ such that
$${\ii\over 2\pi}\Theta_{h_1}(L)={1\over k}
\Big({\ii\over2\pi}\Theta(A)+[D]\Big)$$
with a positive line bundle $A$ and an effective divisor $D$.
Now, for every $\varepsilon>0$, there is a smooth metric $h_\varepsilon$
on $L$ such that ${\ii\over2\pi}\Theta_{h_\varepsilon}(L)\ge-\varepsilon
\omega$, where $\omega={\ii\over2\pi}\Theta(A)$. The convex combination
of metrics $h'_\varepsilon=h_1^{k\varepsilon}h_\varepsilon^{1-k\varepsilon}$
is a singular metric with poles along $D$ which satisfies
$${\ii\over2\pi}\Theta_{h'_\varepsilon}(L)\ge \varepsilon(\omega+[D])-
(1-k\varepsilon)\varepsilon\omega\ge k\varepsilon^2\omega.$$
Its Lelong numbers are $\varepsilon\nu(D,x)$ and they can be made
smaller than $\delta$ by choosing $\varepsilon>0$ small.\qed
\endproof

We still need a few elementary facts about the numerical dimension of nef
line bundles.

\begstat{(6.9) Definition} Let $L$ be a nef line bundle on a compact
K\"ahler manifold~$X$. One defines the numerical dimension of $L$ to be
$$\nu(L)=\max\big\{k=0\ld n\,;\,c_1(L)^k\ne0~\hbox{in}~
H^{2k}(X,\bbbr)\big\}.$$
\endstat

By Corollary 6.8, we have $\kappa(L)=n$ if and only if $\nu(L)=n$.
In general, we merely have an inequality.

\begstat{(6.10) Proposition} If $L$ is a nef line bundle on a compact
K\"ahler manifold, then $\kappa(L)\le\nu(L)$.
\endstat

\begproof{} By induction on $n=\dim X$. If $\nu(L)=n$ or $\kappa(L)=n$
the result is true, so we may assume $r:=\kappa(L)\le n-1$ and $k:=\nu(L)\le
n-1$. Fix $m>0$ so that $\Phi=\Phi_{|mL|}$ has generic rank $r$.
Select a nonsingular ample divisor $A$ in $X$ such that the restriction of
$\Phi_{|mL|}$ to $A$ still has rank $r$ (for this, just take $A$
passing through a point $x\notin B_{|mL|}$ at which $\rank(d\Phi_x)=r<n$,
in such a way that the tangent linear map $d\Phi_{x\restriction T_{A,x}}$
still has rank~$r$). Then $\kappa(L_{\restriction A})\ge r=\kappa(L)$
(we just have an equality because there might exist sections in
$H^0(A,mL_{\restriction A})$ xhich do not extend to~$X$). On the other hand,
we claim that $\nu(L_{\restriction A})=k=\nu(L)$. The inequality
$\nu(L_{\restriction A})\ge\nu(L)$ is clear. Conversely, if we set
$\omega={\ii\over2\pi}\Theta(A)>0$, the cohomology class $c_1(L)^k$ can
be represented by a closed positive current of bidegree $(k,k)$
$$T=\lim_{\varepsilon\to0}\Big({\ii\over 2\pi}\Theta_{h_\varepsilon}(L)+
\varepsilon\omega\Big)^k$$
after passing to some subsequence (there is a uniform bound for the mass
thanks to the K\"ahler assumption, taking wedge products with
$\omega^{n-k}$). The current $T$ must be non zero since $c_1(L)^k\ne0$ by
definition of $k=\nu(L)$. Then $\{[A]\}=\{\omega\}$ as cohomology classes,
and
$$\int_A c_1(L_{\restriction A})^k\wedge\omega^{n-1-k}=
\int_X c_1(L)^k\wedge[A]\wedge\omega^{n-1-k}=
\int_X T\wedge\omega^{n-k}>0.$$
This implies $\nu(L_{\restriction A})\ge k$, as desired. The induction
hypothesis with $X$ replaced by $A$ yields
$$\kappa(L)\le\kappa(L_{\restriction A})\le
\nu(L_{\restriction A})\le\nu(L).\eqno\square$$
\endproof

\begstat{(6.11) Remark} \rm It may happen that $\kappa(L)<\nu(L)$: take e.g.
$$L\to X=X_1\times X_2$$
equal to the total tensor product of an ample line bundle $L_1$ on a
projective manifold $X_1$ and of a unitary flat line bundle $L_2$ on an
elliptic curve $X_2$ given by a representation $\pi_1(X_2)\to U(1)$ such
that no multiple $kL_2$ with $k\ne0$ is trivial. Then
$H^0(X,kL)=H^0(X_1,kL_1)\otimes H^0(X_2,kL_2)=0$ for $k>0$, and thus
$\kappa(L)=-\infty$. However $c_1(L)=\pr_1^\star c_1(L_1)$ has
numerical dimension equal to $\dim X_1$. The same example shows
that the Kodaira dimension may increase by restriction to a
subvariety (if $Y=X_1\times\{{\rm point}\}$, then
$\kappa(L_{\restriction Y})=\dim Y$).\qed
\endstat

We now derive an algebraic version of the Nadel vanishing theorem in the
context of nef line bundles. This algebraic vanishing theorem has been
obtained independently by Kawamata [Kaw82] and Viehweg [Vie82], who both
reduced it to the Kodaira-Nakano vanishing theorem by cyclic covering
constructions. Since then, a number of other proofs have been given,
one based on connections with logarithmic singularities [EV86], another
on Hodge theory for twisted coefficient systems [Kol85], a third one
on the Bochner technique [Dem89] (see also [EV92] for a general survey,
and [Eno93] for an extension to the compact K\"ahler case). Since the
result is best expressed in terms of multiplier ideal sheaves (avoiding
then any unnecessary desingularization in the statement), we feel that
the direct approach via Nadel's vanishing theorem is probably the most
natural one.

If $D=\sum\alpha_jD_j\ge 0$ is an effective $\bbbq$-divisor, we define
the {\it multiplier ideal sheaf} $\cI(D)$ to be equal to $\cI(\varphi)$
where $\varphi=\sum\alpha_j|g_j|$ is the corresponding psh function
defined by generators $g_j$ of $\cO(-D_j)\,$; as we saw in Remark 5.9,
the computation of $\cI(D)$ can be made algebraically by using
desingularizations $\mu:\wt X\to X$ such that $\mu^\star D$ becomes a
divisor with normal crossings on~$\wt X$.

\begstat{(6.12) Kawamata-Viehweg vanishing theorem} Let $X$ be a projective
algebraic manifold and let $F$ be a line bundle over $X$ such that some
positive multiple $mF$ can be written $mF=L+D$ where $L$ is a
nef line bundle and $D$ an effective divisor. Then
$$H^q\big(X,\cO(K_X+F)\otimes\cI(m^{-1}D)\big)=0~~~~
{\it for}~~q>n-\nu(L).$$
\endstat

\begstat{(6.13) Special case} If $F$ is a nef line bundle, then
$$H^q\big(X,\cO(K_X+F)\big)=0~~~~{\it for}~~q>n-\nu(F).$$
\endstat

\begproof{of Theorem 6.12.} First suppose that $\nu(L)=n$, i.e.\ that
$L$ is big. By the proof of 6.5~f), there is a singular hermitian
metric on $L$ such that the corresponding weight $\varphi_{L,0}$
has algebraic singularities and
$$\ii\Theta_0(L)=2\ii d'd''\varphi_L\ge\varepsilon_0\omega$$
for some $\varepsilon_0>0$. On the other hand, since $L$ is nef, there
are metrics given by weights $\varphi_{L,\varepsilon}$ such that
${\ii\over 2\pi}\Theta_\varepsilon(L)\ge\varepsilon\omega$ for
every $\varepsilon>0$, $\omega$ being a K\"ahler metric. Let
$\varphi_D=\sum\alpha_j\log|g_j|$ be the weight of the singular metric on
$\cO(D)$ described in Example~3.13. We define a singular metric on $F$ by
$$\varphi_F={1\over m}
\big((1-\delta)\varphi_{L,\varepsilon}+\delta\varphi_{L,0}+\varphi_D\big)$$
with $\varepsilon\ll\delta\ll 1$, $\delta$ rational. Then $\varphi_F$
has algebraic singularities, and by taking $\delta$ small enough we
find $\cI(\varphi_F)=\cI({1\over m}\varphi_D)=\cI({1\over m}D)$.
In fact, $\cI(\varphi_F)$ can be computed by taking integer parts
of $\bbbq$-divisors (as explained in Remark~5.9), and adding
$\delta\varphi_{L,0}$ does not change the integer part of the rational
numbers involved when $\delta$ is small. Now
$$\eqalign{dd^c\varphi_F
&={1\over m}\big((1-\delta)dd^c\varphi_{L,\varepsilon}+
\delta dd^c\varphi_{L,0}+dd^c\varphi_D\big)\cr
&\ge{1\over m}\big(-(1-\delta)\varepsilon\omega+
\delta\varepsilon_0\omega+[D]\ge{\delta\varepsilon\over m}\omega,\cr}$$
if we choose $\varepsilon\le\delta\varepsilon_0$. Nadel's theorem thus
implies the desired vanishing result for all $q\ge 1$.

Now, if $\nu(L)<n$, we use hyperplane sections and argue by induction
on $n=\dim X$. Since the sheaf $\cO(K_X)\otimes\cI(m^{-1}D)$ behaves
fonctorially with respect to modifications (and since the $L^2$
cohomology complex is ``the same'' upstairs and downstairs),
we may assume after blowing-up that $D$ is a divisor with normal
crossings. By Remark~5.9, the multiplier ideal sheaf
$\cI(m^{-1}D)=\cO(-\lfloor m^{-1}D\rfloor)$ is locally free.
By Serre duality, the expected vanishing is equivalent to
$$H^q(X,\cO(-F)\otimes\cO(\lfloor m^{-1}D\rfloor))=0\qquad
\hbox{\rm for}~q<\nu(L).$$
Then select a nonsingular ample divisor $A$ such that $A$ meets all
components $D_j$ transversally. Select $A$ positive enough so that
$\cO(A+F-\lfloor m^{-1}D\rfloor)$ is ample. Then
$H^q(X,\cO(-A-F)\otimes\cO(\lfloor m^{-1}D\rfloor))=0$ for $q<n$
by Kodaira vanishing, and the exact sequence
$0\to\cO_X(-A)\to\cO_X\to(i_A)_\star\cO_A\to 0$ twisted by
$\cO(-F)\otimes\cO(\lfloor m^{-1}D\rfloor)$ yields an isomorphism
$$H^q(X,\cO(-F)\otimes\cO(\lfloor m^{-1}D\rfloor))\simeq
H^q(A,\cO(-F_{\restriction A})\otimes\cO(\lfloor m^{-1}
D_{\restriction A}\rfloor).$$
The proof of 5.8 showed that $\nu(L_{\restriction A})=\nu(L)$, hence
the induction hypothesis implies that the cohomology group on $A$ in
the right hand side is zero for $q<\nu(L)$.\qed
\endproof

\begstat{(6.14) Remark} \rm Enoki [Eno92] proved that the special
case 6.13 of the Kawamata-Viehweg vansihing theorem still holds when
$X$ is a compact K\"ahler manifold. The idea is to replace the
induction on dimension (hyperplane section argument) by a clever use
of the Aubin-Calabi-Yau theorem. We conjecture that the general case
6.12 is also valid for compact K\"ahler manifolds.\qed
\endstat

\titleb{7.}{Seshadri Constants and the Fujita Conjecture}
In questions related to the search of effective bounds for the existence
of sections in multiples of ample line bundles, one is lead to study the
``local positivity'' of such bundles. In fact, there is a simple
way of giving a precise measurement of positivity. We more or less
follow the original exposition given in [Dem90], and discuss some
interesting results obtained in [EL92] and [EKL94].

\begstat{(7.1) Definition} Let $L$ be a nef line bundle over a projective
algebraic manifold~$X$. The Seshadri constant of $L$ at a point $x\in X$
is the nonnegative real number
$$\varepsilon(L,x)=\sup\big\{\varepsilon\ge 0\,;\,\mu^\star
L-\varepsilon E~\hbox{is nef}\big\},\leqno(7.1')$$
where $\pi:\wt X\to X$ is the blow-up of $X$ at $x$ and $E$ is the
exceptional divisor. An equivalent definition is
$$\varepsilon(L,x)=\inf_{C\ni x}~{L\cdot C\over \nu(C,x)},\leqno(7.1'')$$
where the infimum is taken over all irreducible curves $C$ passing
through $x$ and $\nu(C,x)$ is the multiplicity of $C$ at $x$.
\endstat

To get the equivalence between the two definitions, we just observe that
for any irreducible curve $\wt C\subset\wt X$ not contained in the
exceptional divisor and $C=\pi(\wt C)$, then (exercise to the reader!)
$$(\pi^\star L-\varepsilon E)\cdot \wt C=L\cdot C-\varepsilon\nu(C,x).$$
The infimum
$$\varepsilon(L) =\inf_{x\in X}~\varepsilon (L,x) = \inf_C~{L\cdot C
\over \nu(C)}~~~~\hbox{\rm where}~~\nu(C)=\max_{x\in C}~\nu(C,x)\leqno(7.2)$$
will be called the {\it global Seshadri constant} of $L$. It is well known
that $L$ is ample if and only if $\varepsilon(L)>0$ (Seshadri's criterion
for ampleness, see [Har70] Chapter~1). It is useful to think of the Seshadri
constant $\varepsilon(L,x)$ as measuring how positive $L$ is along
curves passing through~$x$. The following exercise presents some
illustrations of this intuitive idea.

\begstat{(7.3) Exercise} \rm
\item{\rm a)} If $L$ is an ample line bundle such that $mL$ is
{\it very ample}, then $\varepsilon(L)\ge{1\over m}$. (This is the
elementary half of Seshadri's criterion for ampleness).
\smallskip
\item{\rm b)} For any two nef line bundles $L_1,L_2$, show that
$$\varepsilon(L_1+L_2,x)\ge \varepsilon(L_1,x)+\varepsilon(L_2,x)\qquad
\hbox{for all $x\in X$.}\eqno\square$$
\endstat

If $L$ is a nef line bundle, we are also interested in singular
metrics with isolated logarithmic poles (e.g.\ in view of applying
Corollary 5.13): we say that a logarithmic pole $x$ of the weight $\varphi$
is {\it isolated} if $\varphi$ is finite and continuous on $V\ssm \{x\}$
for some neighborhood $V$ of~$x$ and we define
$$\gamma(L,x)=\sup\left.\cases{
\hbox{\rm $\gamma\in\bbbr_+$ such that $L$ has a singular metric}\cr
\hbox{\rm with $\ii\Theta(L)\ge0$ and with an isolated}\cr
\hbox{\rm logarithmic pole of coefficient $\gamma$ at $x$}\cr}\right\}.
\leqno(7.4)$$
If there are no such metrics, we set $\gamma(L,x)=0$.

The numbers $\varepsilon(L,x)$ and $\gamma(L,x)$ will be seen to
carry a lot of useful information about the global sections of
$L$ and its multiples $kL$. To make this precise,
we first introduce some further definitions. Let
$s(L,x)$ be the largest integer $s\in\bbbn$ such that the global
sections in $H^0(X,L)$ generate all $s$-jets $J^s_xL = \cO_x(L)/
\gm_x^{s+1}\cO_x(L)$. If $L_x$ is not generated, i.e.\ if all
sections of $L$ vanish at $x$, we set $s(L,x)=-\infty$. We also
introduce the limit value
$$\sigma (L,x) = \limsup_{k\to +\infty} {1\over k} s(kL,x) =
\sup_{k\in\bbbn^\star} {1\over k} s(kL,x)\leqno(7.5)$$
if $s(kL,x)\ne-\infty$ for some $k$, and $\sigma(L,x)=0$ otherwise.
The limsup is actually equal to the sup thanks to the superadditivity
property
$$s(L_1+L_2,x) \ge s(L_1,x) + s(L_2,x).$$
The limsup is in fact a limit as soon as $kL$ spans at $x$
for $k\ge k_0$, e.g.\ when $L$ is ample.

\begstat{(7.6) Theorem} Let $L$ be a line bundle over $X$.
\smallskip
\item{\rm a)} If $L$ is nef then
$\varepsilon(L,x)\ge\gamma(L,x)\ge\sigma(L,x)$ for every $x\in X$.
\smallskip
\item{\rm b)} If $L$ is ample then
$\varepsilon(L,x)=\gamma(L,x)=\sigma(L,x)$ for every $x\in X$.
\smallskip
\item{\rm c)} If $L$ is big and nef, the equality holds for all $x$
outside some divisor $D$ in~$X$.
\vskip0pt
\endstat

\begproof{} Fix a point $x\in X$ and a coordinate system $(z_1\ld z_n)$
centered at $x$. If $s=s(kL,x)$, then $H^0(X,kL)$
generates all
$s$-jets at $x$ and we can find holomorphic sections $f_1\ld f_N$ whose
$s$-jets are all monomials $z^\alpha$, $|\alpha|=s$. We define a global
singular metric on $L$ by
$$|\xi | = \bigg( \sum_{1\le j\le N} |f_j(z)\cdot \xi
^{-k}|^2\bigg)^{-1/2k},~\xi \in L_z\leqno(7.7)$$
associated to the weight function $\varphi(z)={1\over 2k}\log\sum\big|
\theta(f_j(z))\big|^2$ in any trivialization $L_{\restriction\Omega}\simeq
\Omega\times\bbbc$. Then $\varphi $ has an isolated logarithmic pole of
coefficient $s/k$ at $x$, thus
$$\gamma (L,x) \ge {1\over k} s(kL,x)$$
and in the limit we get $\gamma(L,x)\ge\sigma(L,x)$.

Now, suppose that $L$ has a singular metric with an isolated log pole of
coefficient $\ge\gamma$ at $x$. Set ${\ii\over 2\pi}\Theta(L)=
dd^c\varphi$  on a neighborhood $\Omega$ of $x$ and let $C$
be an irreducible curve passing through $x$. Then all weight functions
associated to the metric of $L$ must be locally integrable along $C$
(since $\varphi $ has an isolated pole at $x$). We infer
$$L\cdot C=\int_C{\ii\over2\pi}\Theta(L)\ge\int_{C\cap\Omega}dd^c\varphi
\ge\gamma\,\nu(C,x)$$
because the last integral is larger than the Lelong number of the
current $[C]$ with respect to the weight $\varphi$ and we
may apply our comparison theorem 2.15 with the ordinary Lelong number
associated to the weight $\log|z-x|$. Therefore
$$\varepsilon (L,x) = \inf {L\cdot C\over \nu(C,x)} \ge \sup \gamma
=\gamma (L,x).$$

Finally, we show that $\sigma(L,x)\ge\varepsilon(L,x)$ when $L$ is
ample. This is done essentially by same arguments as in the proof
of Seshadri's criterion, as explained in [Har70]. Consider
the blow-up $\pi:\smash{\wt X}\to X$
at point $x$, the exceptional divisor $E=\pi^{-1}(x)$ and the
line bundles $F_{p,q} =\cO(p\,\pi^\star L-q\,E)$ over $\smash{\wt X}$,
where $p,q>0$. Recall that $\cO(-E)_{\restriction E}$ is the canonical
line bundle $\cO_E(1)$ over $E\simeq\bbbp^{n-1}$, in particular we have
$E^n=\cO_E(-1)^{n-1}=(-1)^{n-1}$. For any irreducible curve
$\smash{\wt C}\subset\smash{\wt X}$, either $\smash{\wt C}\subset E$ and
$$F_{p,q}\cdot \wt C = \cO(-q\,E)\cdot \wt C = q\,\cO_E(1)\cdot\wt C=
q\deg\wt C$$
or $\pi (\wt C) = C$ is a curve and
$$F_{p,q}\cdot C=p\,L\cdot C - q\,\nu(C,x) \ge \big(p-q/\varepsilon
(L,x)\big)L\cdot C.$$
Thus $F_{p,q}$ is nef provided that $p\ge q/\varepsilon (L,x)$.
Since $F_{p,q}$ is ample when $p/q$ is large, a simple interpolation
argument shows that $F_{p,q}$ is ample for $p>q/\varepsilon(L,x)$.
In that case, the Kodaira-Serre vanishing theorem gives
$$H^1(\wt X,k\,F_{p,q}) = H^1\big(\wt X,\cO(kp\,\pi^\star L-kq\,E)\big)=0$$
for $k$ large. Hence we get a surjective map
$$H^0(\wt X,kp\,\pi^\star L) \lraww  H^0\Big(\wt X,
\cO(kp\,\pi^\star L)\otimes(\cO/\cO(-kq\,E)\big)\Big)\simeq
J_x^{kq-1}(kp\,L),$$
that is, $H^0(X,kp\,L)$ generates all $(kq-1)$ jets at $x$.
Therefore $p>q/\varepsilon (L,x)$ implies $s(kp\,L,x) \ge kq-1$
for $k$ large, so $\sigma (L,x)\ge q/p$. At the limit we get
$\sigma(L,x)\ge\varepsilon(L,x)$.

Assume now that $L$ is nef and big and that $\varepsilon(L,x)>0$. By the
proof of 6.6~f), there exist an integer $k_0\ge 1$ and effective divisors
$A,D$ such that $k_0L\simeq A+D$ where $A$ is ample. Then $a\,\pi^\star A-E$
is ample for $a$ large. Hence there are
integers $a,b>0$ such that $a\,\pi^\star A-b\,E-K_{\wt X}$
is ample. When $F_{p,q}$ is nef, the sum with any positive multiple
$k\,F_{p,q}$ is still ample and the Akizuki-Kodaira-Nakano vanishing
theorem gives
$$H^1(\wt X,k\,F_{p,q}+a\,\pi^\star A-b\,E)
=H^1\big(\wt X,(kp+k_0a)\,\pi^\star L-a\,\pi^\star D-(kq+b)E\big)=0$$
when we substitute $A=k_0L-D$. As above, this implies that we have a
surjective map
$$H^0\big(X,(kp+k_0a)\,L-a\,D\big)\lraww
J^{kq+b-1}_x\big((kp+k_0a)\,L-a\,D\big)$$
when $p\ge q/\varepsilon(L,x)$. Since $\cO(-aD)\subset\cO$, we infer
$s\big((kp+k_0a)L,x\big)\ge kq+b-1$ at every point $x\in X\ssm D$
and at the limit $\sigma(L,x)\ge\varepsilon(L,x)$.\qed
\endproof

\begstat{(7.8) Remark} \rm Suppose that the line bundle $L$ is ample. The
same arguments show that if $\pi :\wt X \to X$ is the blow-up at two points
$x,y$ and if $E_x+E_y$ is the exceptional divisor,
then $F_{p,q} =p\,\pi^\star L-q\,E_x-E_y$ is
ample for $p>q/\varepsilon (L,x) + 1/\varepsilon (L,y)$. In that  case,
$H^0(X,kp\,L)$ generates $J^{kq-1}_x(kp\,L)\oplus
J^{k-1}_y(kp\,L)$ for $k$ large. Take $p>q/\varepsilon
(L,x)+1/\varepsilon (L)$ and let $y$ run over $X\ssm  \{x\}$.
For $k$ large, we obtain sections $f_j\in H^0(X,kp\,L)$ whose
jets at $x$ are all monomials $z^\alpha $, $|\alpha|=kq-1$, and with
no other common zeros. Moreover, Formula (7.7) produces a metric on $L$
which is smooth and has positive definite curvature on $X\ssm
\{x\}$, and which has a log pole of coefficient $(kq-1)/kp$ at $x$.
Therefore the supremum $\gamma (L,x) =\sup\{\gamma \}$ is always
achieved by metrics which are smooth and have positive definite curvature
on $X\ssm \{x\}$.\qed
\endstat

\begstat{(7.9) Remark} \rm If $Y$ is a $p$-dimensional algebraic subset
of $X$ passing through~$x$, then
$$L^p\cdot Y\ge \varepsilon (L,x)^p \nu(Y,x)$$
where $L^p\cdot Y=\int_Yc_1(L)^p$ and $\nu(Y,x)$ is the multiplicity of
$Y$ at $x$. In fact, if $L$ is ample, we can take a metric on $L$
which is smooth on $X\ssm \{x\}$ and defined on a neighborhood $\Omega $
of $x$ by a weight function $\varphi $ with a log pole of coefficient
$\gamma $ at $x$. By the comparison theorem for Lelong numbers, we get
$$L^p\cdot Y\ge\int_{Y\cap\Omega}(dd^c\varphi)^p\ge \gamma^p\nu(Y,x)$$
and $\gamma $ can be chosen arbitrarily close to $\varepsilon
(L,x)$. If $L$ is nef, we apply the inequality to
$k\,L+M$ with $M$ ample and take the limit as $k\to
+\infty$.\qed
\endstat

The Seshadri constants $\varepsilon (L,x)$ and $\varepsilon(L)=\inf
\varepsilon (L,x)$ are especially interesting because they
provide effective results concerning the existence of sections of
the so called {\it adjoint line bundle} $K_X+L$. The following proposition
illustrates this observation.

\begstat{(7.10) Proposition} Let $L$ be a big nef line bundle over~$X$.
\smallskip
\item{\rm a)} If $\varepsilon(L,x)>n+s$, then $H^0(X,K_X+L)$
generates all $s$-jets at~$x$.
\smallskip
\item{\rm b)} If $\varepsilon(L)>2n$, then $K_X+L$ is very ample.
\endstat

\begproof{} a) By the proof of Theorem~7.6, the line bundle
$\pi^\star L-q\,E$ is nef for $q\le\varepsilon (L,x)$. Moreover,
its $n$-th self intersection is equal to $L^n+(-q)^nE^n=L^n-q^n$ and
as $L^n\ge\varepsilon(L,x)^n$ by remark 3.5, we see that
$\pi^\star L-q\,E$ is big for $q<\varepsilon (L,x)$.
The Kawamata-Viehweg vanishing theorem~5.2 then gives
$$H^1(\wt X,K_{\wt X}+\pi^\star L-q\,E)=
H^1\big(\wt X,\pi^\star K_X+\pi^\star L-(q-n+1)E\big)=0,$$
since $K_{\wt X}=\pi^\star K_X+(n-1)E$. Thus we get a
surjective map
$$\cmalign{
H^0\big(\wt X,\pi^\star&K_X+\pi^\star L)&\lraww
H^0\big(\wt X,\pi^\star\cO(K_X&+ L)\otimes
\cO/\cO(-(q-n+1)E)\big)\cr
&\|&&\|\cr
\hfill H^0(X,&K_X+L)&\lraww\hfill J^{q-n}_x(&K_X+L)\cr}$$
provided that $\varepsilon (L,x)>q$. The first statement is proved. To
show that $K_X+L$ is very ample, we blow up at two points $x,y$. The
line bundle $\pi^\star L-n\,E_x-n\,E_y$ is ample for $1/\varepsilon
(L,x) + 1/\varepsilon (L,y) < 1/n$, a sufficient condition for this is
$\varepsilon (L)>2n$. Then we see that
$$H^0(X,K_X+L)\lra(K_X+L)_x\oplus(K_X+L)_y$$
is also surjective.\qed
\endproof

\begstat{(7.11) Exercise} \rm Derive Proposition~7.10 directly from
Corollary 5.13, assuming that $L$ is ample and using the equality
$\varepsilon(L,x)=\gamma(L,x)$.\qed
\endstat

In relation with these questions, Fujita [Fuj87,~88] made the following
interesting conjecture about adjoint line bundles.

\begstat{(7.12) Conjecture (Fujita)} If $L$ is an ample line bundle, then
$K_X+mL$ is generated by global sections for $m\ge n+1$ and very ample for
$m\ge n+2$.
\endstat

Using Mori theory, Fujita proved that $K_X+mL$ is nef for $m\ge n+1$
and ample for $m\ge n+2$, but these results are of course much weaker
than the conjecture; they can be derived rather easily by a direct
application of the Kawamata-Viehweg vanishing theorem (see Section~8).
Observe that if the conjecture holds, it would be actually optimal, as shown
by the case of projective space: if $X=\bbbp^n$, then $K_X=\cO(-n-1)$,
hence the bounds $m=n+1$ for global generation and $m=n+2$ for very
ampleness are sharp when $L=\cO(1)$. The case of curves $(n=1)$ is
easily settled in the affirmative:

\begstat{(7.13) Exercise} \rm If $X$ is a curve and $L$ is a line bundle
of positive degree on $X$, show by Riemann-Roch that
$\varepsilon(L,x)=\sigma(L,x)=\deg L$ at every point. Show in this case
that Fujita's conjecture follows from Proposition~7.10.\qed
\endstat

In the case of surfaces $(n=2)$, Fujita's conjecture has been proved by
I.~Reider in [Rei88], as a special case of a stronger numerical criterion for
very ampleness (a deep generalization of Bombieri's work [Bom73] on
pluricanonical embeddings of surfaces of general type). Reider's method
is based on the Serre construction for rank 2 bundles and on the Bogomolov
instability criterion. Since then, various other proofs have been obtained
(Sakai [Sak88], Ein-Lazarsfeld [EL93]). We will give in the next section an
algebraic proof of Reider's result based on vanishing theorems, following
closely [EL93]. In higher dimensions, the global generation part of the
statement has been proved by Ein-Lazarsfeld [EL93] for $n=3$. However,
up to now, there is no strong indication that the conjecture should be true
in higher dimensions.

In a naive attempt to prove Fujita's conjecture using part a) of
Proposition~7.10, it is natural to ask whether one has
$\varepsilon(L,x)\ge 1$ when $L$ is an ample line bundle.
Unfortunately, simple examples due to R.~Miranda show that
$\varepsilon(L,x)$ may be arbitrarily small as soon as $\dim X\ge2$.

\begstat{(7.14) Proposition (R.~Miranda)} Given $\varepsilon>0$, there
exists a rational surface $X$, a point $x\in X$ and an ample line bundle
$L$ on $X$ such that $\varepsilon(L,x)\le\varepsilon$.
\endstat

\begproof{} Let $C\subset\bbbp^2$ be an irreducible curve of large degree
$d$ with a point $x$ of multiplicity $m$. Let $C'$ be another irreducible
curve of the same degree meeting $C$ transversally. Blow-up the points
of $C\cap C'$ to obtain a rational surface $X$, admitting a map
$\varphi:X\to\bbbp^1$ (the map $\varphi$ is simply given by
$\varphi(z)=P'(z)/P(z)$ where $P=0$, $P'=0$ are equations of $C$ and
$C'$). The fibers of $\varphi$ are the curves in the pencil spanned by
$C$ and $C'$. If $C'$ is chosen general enough, then all these fibers are
irreducible. As the fibers are disjoint, we have $C^2=C\cdot C'=0$ in $X$.
Now, let $E$ be one of the components of the exceptional divisor
of~$X$. Fix an integer $a\ge 2$. It follows from the Nakai-Moishezon
criterion that the divisor $L=aC+E$ is ample: in fact $L^2=2a C\cdot E+E^2
=2a-1$, $L\cdot E=a-1$, $L\cdot\Gamma=L\cdot C=1$ if $\Gamma$ is a
fiber of $\varphi\,$; all other irreducible curves $\Gamma\subset X$
must satisfy $\varphi(\Gamma)=\bbbp^1$, hence they have non empty
intersection with $C$ and $L\cdot\Gamma\ge a$. However $\nu(C,x)=m$ by
construction, hence $\varepsilon(L,x)\le 1/m$.\qed
\endproof

Although Seshadri constants fail to be bounded below in a uniform
way, the following statement (which would imply the Fujita conjecture
at generic points) is expected to be true.

\begstat{(7.15) Conjecture} If $L$ is a big nef line bundle, then
$\varepsilon(L,x)\ge 1$ for $x$ generic.
\endstat

By a generic point, we mean a point in the complement of a proper algebraic
subvariety, or possibly, in a countable union of proper algebraic
subvarieties. Quite recently, major progress has been made on this
question by Ein-Lazarsfeld [EL92] and Ein-K\"uchle-Lazarsfeld [EKL94].

\begstat{(7.16) Theorem {\rm([EL92])}} Let $L$ be an ample line bundle
on a smooth projective surface $X$. Then $\varepsilon(L,x)\ge 1$ for all
except perhaps countably many points $x\in X$, and moreover if $L^2>1$
then the set of exceptional points is finite. If $L^2\ge 5$ and
$L\cdot C\ge 2$ for all curves $C\subset X$, then $\varepsilon(L,x)\ge 2$
for all but finitely many $x\in X$.\qed\eject
\endstat

\begstat{(7.17) Theorem {\rm([EKL94])}} Let $L$ be a big nef line bundle
on a projective $n$-dimensional manifold~$X$. Then $\varepsilon(L,x)\ge 1/n$
for all $x$ outside a countable union of proper algebraic subsets.\qed
\endstat

\titleb{8.}{Algebraic Approach to the Fujita Conjecture}
This section is devoted to a proof of various results related to the
Fujita conjecture. The main ideas occuring here are inspired by a
recent work of Y.T.~Siu [Siu94a]. His method, which is algebraic in
nature and quite elementary, consists in a combination of the
Riemann-Roch formula together with Nadel's vanishing theorem (in fact,
only the algebraic case is needed, thus the original Kawamata-Viehweg
vanishing theorem would be sufficient). Slightly later, Angehrn and
Siu [AS94], [Siu94b] introduced other closely related methods, producing
better bounds for the global generation question; since their method
is rather delicate, we can only refer the ready to the above references.
In the sequel, $X$ denotes a projective algebraic $n$-dimensional manifold.
The first observation is the following well-known consequence of the
Riemann-Roch formula.

\begstat{(8.1) Special case of Riemann-Roch} Let $\cJ\subset\cO_X$ be a
coherent ideal sheaf on $X$ such that the subscheme $Y=V(\cJ)$ has
dimension $d$ $($with possibly some lower dimensional components$)$.
Let $[Y]=\sum\lambda_j[Y_j]$ be the effective algebraic cycle of dimension
$d$ associated to the $d$ dimensional components of~$Y$ $($taking
into account multiplicities $\lambda_j$ given by the ideal $\cJ)$. Then
for any line bundle $F$, the Euler characteristic
$$\chi(Y,\cO(F+mL)_{\restriction Y})=\chi(X,\cO(F+mL)\otimes\cO_X/\cJ)$$
is a polynomial $P(m)$ of degree $d$ and leading coefficient
$L^d\cdot[Y]/d!$
\endstat

The second fact is an elementary lemma about numerical polynomials
(polynomials with rational coefficients, mapping $\bbbz$ into $\bbbz$).

\begstat{(8.2) Lemma} Let $P(m)$ be a numerical polynomial of degree $d>0$
and leading coefficient $a_d/d!$, $a_d\in\bbbz$, $a_d>0$. Suppose that
$P(m)\ge 0$ for $m\ge m_0$. Then
\smallskip
\item{\rm a)} For every integer $N\ge0$, there exists $m\in[m_0,m_0+Nd]$
such that~\hbox{$P(m)\ge N$}.
\smallskip
\item{\rm b)} For every $k\in\bbbn$, there exists $m\in[m_0,m_0+kd]$ such
that \hbox{$P(m)\ge a_dk^d/2^{d-1}$}.
\smallskip
\item{\rm c)} For every integer $N\ge 2d^2$, there exists
$m\in[m_0,m_0+N]$ such that \hbox{$P(m)\ge N$}.
\vskip0pt
\endstat

\begproof{} a) Each of the $N$ equations $P(m)=0$, $P(m)=1$,
$\ldots$, $P(m)=N-1$ has at most $d$ roots, so there must be an integer
$m\in[m_0,m_0+dN]$ which is not a root of these.
\smallskip\noindent
b) By Newton's formula for iterated differences $\Delta P(m)=P(m+1)-P(m)$,
we get
$$\Delta^dP(m)=\sum_{1\le j\le d}(-1)^j{d\choose j}P(m+d-j)=a_d,
\qquad\forall m\in\bbbz.$$
Hence if $j\in\big\{0,2,4\ld 2\lfloor d/2\rfloor\big\}\subset[0,d]$ is the
even integer achieving the maximum of $P(m_0+d-j)$ over this finite set,
we find
$$2^{d-1}P(m_0+d-j)=\left({d\choose 0}+{d\choose 2}+\ldots\right)P(m_0+d-j)
\ge a_d,$$
whence the existence of an integer $m\in[m_0,m_0+d]$ with $P(m)\ge
a_d/2^{d-1}$. The case $k=1$ is thus proved. In general, we apply the
above case to the polynomial $Q(m)=P(km-(k-1)m_0)$, which has leading
coefficient $a_dk^d/d!$
\smallskip\noindent
c) If $d=1$, part a) already yields the result. If $d=2$, a look at the
parabola shows that
$$\max_{m\in[m_0,m_0+N]}P(m)\ge\cases{
a_2N^2/8&if $N$ is even,\cr
a_2(N^2-1)/8&if $N$ is odd;\cr}$$
thus $\max_{m\in[m_0,m_0+N]}P(m)\ge N$ whenever $N\ge 8$. If $d\ge 3$,
we apply b) with $k$ equal to the smallest integer such
that $k^d/2^{d-1}\ge N$, i.e.\ $k=\lceil 2(N/2)^{1/d}\rceil$, where
$\lceil x\rceil\in\bbbz$ denotes the round-up of $x\in\bbbr$.
\hbox{Then $kd\le (2(N/2)^{1/d}+1)d\le N$} whenever $N\ge 2d^2$, as a short
computation shows.\qed
\endproof

We now apply Nadel's vanishing theorem pretty much in the same way as
Siu [Siu94a], but with substantial simplifications in the technique and
improvements in the bounds. Our method yields simultaneously a simple
proof of the following basic result.

\begstat{(8.3) Theorem} If $L$ is an ample line bundle over a projective
$n$-fold~$X$, then the adjoint line bundle $K_X+(n+1)L$ is nef.
\endstat

By using Mori theory and the base point free theorem ([Mor82],
[Kaw84]), one can even show that $K_X+(n+1)L$ is semiample, i.e., there
exists a positive integer $m$ such that $m(K_X+(n+1)L)$ is generated by
sections (see [Kaw85] and [Fuj87]). The proof rests on the observation
that $n+1$ is the maximal length of extremal rays of smooth projective
$n$-folds. Our proof of (8.3) is different and will be given simultaneously
with the proof of Th.~(8.4) below.

\begstat{(8.4) Theorem} Let $L$ be an ample line bundle and let $G$ be
a nef line bundle on a projective $n$-fold~$X$. Then the following
properties hold.
\smallskip
\item{\rm a)} $2K_X+mL+G$ generates simultaneous jets of order
$s_1\ld s_p\in\bbbn$ at arbitrary points $x_1\ld x_p\in X$, i.e., there
is a surjective map
$$H^0(X,2K_X+mL+G)\lraww\bigoplus_{1\le j\le p}\cO(2K_X+mL+G)\otimes
\cO_{X,x_j}/\gm_{X,x_j}^{s_j+1},$$
provided that $\displaystyle m\ge 2+\sum_{1\le j\le p}{3n+2s_j-1\choose
n}$.\hfill\break In particular $2K_X+mL+G$ is very ample for
$\displaystyle m\ge 2+{3n+1\choose n}$.
\smallskip
\item{\rm b)} $2K_X+(n+1)L+G$ generates simultaneous jets of order
$s_1\ld s_p$ at arbitrary points $x_1\ld x_p\in X$ provided that
the intersection numbers $L^d\cdot Y$ of $L$ over all
$d$-dimensional algebraic subsets $Y$ of $X$ satisfy
$$L^d\cdot Y>{2^{d-1}\over\lfloor n/d\rfloor^d}\sum_{1\le j\le p}
{3n+2s_j-1\choose n}.$$
\endstat

\begproof{} The proofs of (8.3) and (8.4$\,$a,$\,$b) go along the
same lines, so we deal with them simultaneously (in the case of
(8.3), we simply agree that $\{x_1\ld x_p\}=\emptyset$). The idea is
to find an integer (or rational number) $m_0$ and a singular hermitian
metric $h_0$ on $K_X+m_0L$ with strictly positive curvature current
$\ii\Theta_{h_0}\ge\varepsilon\omega$, such that $V(\cI(h_0))$ is
$0$-dimensional and the weight $\varphi_0$ of $h_0$ satis\-fies
$\nu(\varphi_0,x_j)\ge n+s_j$ for all~$j$. As $L$ and $G$ are nef,
$(m-m_0)L+G$ has for all $m\ge m_0$ a metric $h'$ whose
curvature $\ii\Theta_{h'}$ has arbitrary small negative part (see
[Dem90]), e.g., $\ii\Theta_{h'}\ge -{\varepsilon\over 2}\omega$. Then
$\ii\Theta_{h_0}+\ii\Theta_{h'}\ge{\varepsilon \over 2}\omega$ is again
positive definite. An application of Cor~(1.5) to
$F=K_X+mL+G=(K_X+m_0 L)+((m-m_0)L+G)$ equipped with the metric
$h_0\otimes h'$ implies the existence of the desired sections in
$K_X+F=2K_X+mL+G$ for $m\ge m_0$.

Let us fix an embedding $\Phi_{|\mu L|}:X\to\bbbp^N$, $\mu\gg 0$, given
by sections $\lambda_0\ld\lambda_N\in H^0(X,\mu L)$, and let $h_L$ be
the associated metric on $L$ of positive definite curvature form
$\omega={\ii\over 2\pi}\Theta(L)$. In order to obtain the desired
metric $h_0$ on $K_X+m_0 L$, we fix $a\in\bbbn^\star$ and use a double
induction process to construct singular metrics $(h_{k,\nu})_{\nu\ge
1}$ on $aK_X+b_kL$ for a non increasing sequence of positive integers
$b_1\ge b_2\ge\ldots\ge b_k\ge\ldots\,$. Such a sequence much be
stationary and $m_0$ will just be the stationary limit $m_0=\lim
b_k/a$. The metrics $h_{k,\nu}$ are taken to satisfy the following
properties:
\smallskip
\item{$\alpha\!)$} $h_{k,\nu}$ is an algebraic metric of the form
$$\|\xi\|_{h_{k,\nu}}^2={|\tau_k(\xi)|^2\over
\big(\sum_{1\le i\le\nu,\,0\le j\le N}\big|\tau_k^{(a+1)\mu}(\sigma_i^{a\mu}
\cdot\lambda_j^{(a+1)b_k-am_i})\big|^2\big)^{1/(a+1)\mu}},$$
defined by sections $\sigma_i\in H^0(X,(a+1)K_X+m_iL)$, $m_i<{a+1\over a}b_k$,
$1\le i\le\nu$, where $\xi\mapsto\tau_k(\xi)$ is an arbitrary local
trivialization of $aK_X+b_kL\,$; note that $\sigma_i^{a\mu}\cdot
\lambda_j^{(a+1)b_k-am_i}$ is a section of
$$a\mu((a+1)K_X+m_iL)+((a+1)b_k-am_i)\mu L=(a+1)\mu(aK_X+b_kL).$$
\item{$\beta\kern-1pt)$\kern-1pt} $\ord_{x_j}(\sigma_i)\ge(a+1)(n+s_j)$ for
all $i,j\,$;
\smallskip
\item{$\gamma)$} $\cI(h_{k,\nu+1})\supset\cI(h_{k,\nu})$ and
$\cI(h_{k,\nu+1})\ne\cI(h_{k,\nu})$ whenever the zero variety
$V(\cI(h_{k,\nu}))$ has positive dimension.
\medskip\noindent
The weight $\varphi_{k,\nu}={1\over 2(a+1)\mu}\log\sum\big|\tau_k^{(a+1)\mu}
(\sigma_i^{a\mu}\cdot\lambda_j^{(a+1)b_k-am_i})\big|^2$ of $h_{k,\nu}$ is
plurisubharmonic and the condition $m_i<{a+1\over a}b_k$ implies
$(a+1)b_k-am_i\ge 1$, thus the difference
$\varphi_{k,\nu}-{1\over2(a+1)\mu}\log\sum|\tau(\lambda_j)|^2$ is also
plurisubharmonic. Hence ${\ii\over 2\pi}\Theta_{h_{k,\nu}}(aK_X+b_kL)=
{\ii\over\pi}d'd''\varphi_{k,\nu}\ge{1\over(a+1)}\omega$. Moreover, condition
$\beta)$ clearly implies $\nu(\varphi_{k,\nu},x_j)\ge a(n+s_j)$. Finally,
condition $\gamma)$ combined with the strong Noetherian property of coherent
sheaves ensures that the sequence $(h_{k,\nu})_{\nu\ge 1}$ will finally
produce a zero dimensional subscheme $V(\cI(h_{k,\nu}))$. We agree that
the sequence $(h_{k,\nu})_{\nu\ge 1}$ stops at this point, and we denote
by $h_k=h_{k,\nu}$ the final metric, such that $\dim V(\cI(h_k))=0$.

For $k=1$, it is clear that the desired metrics $(h_{1,\nu})_{\nu\ge 1}$
exist if $b_1$ is taken large enough (so large, say, that
$(a+1)K_X+(b_1-1)L$ generates jets of order $(a+1)(n+\max s_j)$ at every
point; then the sections $\sigma_1\ld\sigma_\nu$ can be chosen with
$m_1=\ldots=m_\nu=b_1-1$). Suppose that the metrics $(h_{k,\nu})_{\nu\ge 1}$
and $h_k$ have been constructed and let us proceed with the construction
of $(h_{k+1,\nu})_{\nu\ge 1}$. We do this again by induction on~$\nu$,
assuming that $h_{k+1,\nu}$ is already constructed and that $\dim
V(\cI(h_{k+1,\nu}))>0$. We start in fact the induction with $\nu=0$, and
agree in this case that $\cI(h_{k+1,0})=0$ (this would correspond to an
infinite metric of weight identically equal to $-\infty$). By Nadel's
vanishing theorem applied to
$$F_m=aK_X+mL=(aK_X+b_kL)+(m-b_k)L$$
with the metric $h_k\otimes(h_L)^{\otimes m-b_k}$, we get
$$H^q(X,\cO((a+1)K_X+mL)\otimes\cI(h_k))=0\qquad
\hbox{\rm for $q\ge 1$, $m\ge b_k$.}$$
As $V(\cI(h_k))$ is $0$-dimensional, the sheaf $\cO_X/\cI(h_k)$ is a
skyscraper sheaf, and the exact sequence $0\to\cI(h_k)\to\cO_X\to
\cO_X/\cI(h_k)\to 0$ twisted with the invertible sheaf $\cO((a+1)K_X+mL)$
shows that
$$H^q(X,\cO((a+1)K_X+mL))=0\qquad \hbox{\rm for $q\ge 1$, $m\ge b_k$.}$$
Similarly, we find
$$H^q(X,\cO((a+1)K_X+mL)\otimes\cI(h_{k+1,\nu}))=0\qquad
\hbox{\rm for $q\ge 1$, $m\ge b_{k+1}$}$$
(also true for $\nu=0$, since $\cI(h_{k+1,0})=0$), and when
\hbox{$m\ge\max(b_k,b_{k+1})=b_k$}, the exact sequence $0\to\cI(h_{k+1,\nu})
\to\cO_X\to\cO_X/\cI(h_{k+1},\nu)\to 0$ implies
$$H^q(X,\cO((a+1)K_X+mL)\otimes\cO_X/\cI(h_{k+1,\nu}))=0\qquad
\hbox{\rm for $q\ge 1$, $m\ge b_k$.}$$
In particular, since the $H^1$ group vanishes, every section $u'$ of
$(a+1)K_X+mL$ on the subscheme $V(\cI(h_{k+1,\nu}))$ has an extension
$u$ to $X$. Fix a basis $u_1'\ld u_N'$ of the sections on
$V(\cI(h_{k+1,\nu}))$ and take arbitrary extensions $u_1\ld u_N$ to
$X$. Look at the linear map assigning the collection of jets of order
$(a+1)(n+s_j)-1$ at all points~$x_j$
$$u=\sum_{1\le j\le N}a_ju_j\longmapsto
\bigoplus J^{(a+1)(n+s_j)-1}_{x_j}(u).$$
Since the rank of the bundle of $s$-jets is ${n+s\choose n}$, the target
space has dimension
$$\delta=\sum_{1\le j\le p}{n+(a+1)(n+s_j)-1\choose n}.$$
In order to get a section $\sigma_{\nu+1}=u$ satisfying condition $\beta)$
with non trivial restriction $\sigma'_{\nu+1}$ to $V(\cI(h_{k+1,\nu}))$, we
need at least $N=\delta+1$ independent sections $u'_1\ld u'_N$. This
condition is achieved by applying Lemma~(8.2) to the numerical polynomial
$$\eqalign{P(m)
&=\chi(X,\cO((a+1)K_X+mL)\otimes\cO_X/\cI(h_{k+1,\nu}))\cr
&=h^0(X,\cO((a+1)K_X+mL)\otimes\cO_X/\cI(h_{k+1,\nu}))\ge 0,\qquad
m\ge b_k.\cr}$$
The polynomial $P$ has degree $d=\dim V(\cI(h_{k+1,\nu}))>0$.
We get the existence of an integer $m\in[b_k,b_k+\eta]$ such that
$N=P(m)\ge\delta+1$ with some explicit integer $\eta\in\bbbn$
(for instance $\eta=n(\delta+1)$ always works by (8.2$\,$a), but we will
also use the other possibilities to find an optimal choice in each case).
Then we find a section $\sigma_{\nu+1}\in H^0(X,(a+1)K_X+mL)$
with non trivial restriction $\sigma'_{\nu+1}$ to $V(\cI(h_{k+1,\nu}))$,
vanishing at order $\ge(a+1)(n+s_j)$ at each point $x_j$. We just set
$m_{\nu+1}=m$, and the condition $m_{\nu+1}<{a+1\over a}b_{k+1}$ is
satisfied if $b_k+\eta<{a+1\over a}b_{k+1}$. This shows that we can take
inductively
$$b_{k+1}=\left\lfloor{a\over a+1}(b_k+\eta)\right\rfloor+1.$$
By definition, $h_{k+1,\nu+1}\le h_{k+1,\nu}$, hence $\cI(h_{k+1,\nu+1})
\supset\cI(h_{k+1,\nu})$. We necessarily have $\cI(h_{k+1,\nu+1})\ne
\cI(h_{k+1,\nu})$, for $\cI(h_{k+1,\nu+1})$ contains the ideal sheaf
associated with the zero divisor of $\sigma_{\nu+1}$, whilst
$\sigma_{\nu+1}$ does not vanish identically on~$V(\cI(h_{k+1,\nu}))$.
Now, an easy computation shows that the iteration of
$b_{k+1}=\lfloor{a\over a+1}(b_k+\eta)\rfloor+1$ stops at $b_k=a(\eta+1)+1$ for
any large initial value $b_1$. In this way, we obtain a metric $h_\infty$
of positive definite curvature on $aK_X+(a(\eta+1)+1)L$, with
$\dim V(\cI(h_\infty))=0$ and $\nu(\varphi_\infty,x_j)\ge a(n+s_j)$ at
each point~$x_j$.
\endproof

\begproof{of $(8.3)$.} In this case, the set $\{x_j\}$ is taken to be
empty, thus $\delta=0$. By (8.2$\,$a), the condition $P(m)\ge 1$ is
achieved for some $m\in[b_k,b_k+n]$ and we can take $\eta=n$. As $\mu L$ is
very ample, there exists on $\mu L$ a metric with an isolated logarithmic
pole of Lelong number $1$ at any given point~$x_0$ (e.g., the algebraic
metric defined with all sections of $\mu L$ vanishing at~$x_0$). Hence
$$F_a'=aK_X+(a(n+1)+1)L+n\mu L$$
has a metric $h_a'$ such that $V(\cI(h_a'))$ is zero dimensional and
contains $\{x_0\}$. By Cor~(1.5), we conclude that
$$K_X+F_a'=(a+1)K_X+(a(n+1)+1+n\mu)L$$
is generated by sections, in particular $K_X+{a(n+1)+1+n\mu\over a+1}L$
is nef. As $a$ tends to $+\infty$, we infer that $K_X+(n+1)L$ is nef.\qed
\endproof

\begproof{of $(8.4\,${\rm a).}} Here, the choice $a=1$ is sufficient for
our purposes. Then
$$\delta=\sum_{1\le j\le p}{3n+2s_j-1\choose n}.$$
If $\{x_j\}\ne\emptyset$, we have $\delta+1\ge{3n-1\choose n}+1\ge 2n^2$
for $n\ge 2$. Lemma (8.2$\,$c) shows that $P(m)\ge\delta+1$ for some
$m\in[b_k,b_k+\eta]$ with $\eta=\delta+1$. We can start in fact the
induction procedure $k\mapsto k+1$ with $b_1=\eta+1=\delta+2$, because the
only property needed for the induction step is the vanishing
property
$$H^0(X,2K_X+mL)=0\qquad\hbox{for $q\ge 1$, $m\ge b_1$},$$
which is true by the Kodaira vanishing theorem and the ampleness
of $K_X+b_1L$ (here we use Fujita's result (8.3), observing that $b_1>n+1$).
Then the recursion formula $b_{k+1}=\lfloor{1\over 2}(b_k+\eta)\rfloor+1$
yields $b_k=\eta+1=\delta+2$ for all $k$, and (8.4$\,$a) follows.\qed
\endproof

\begproof{of $(8.4\,${\rm b).}} Quite similar to (8.4$\,$a), except that we
take $\eta=n$, $a=1$ and $b_k=n+1$ for all~$k$. By Lemma (8.2$\,$b), we have
$P(m)\ge a_dk^d/2^{d-1}$ for some integer $m\in[m_0,m_0+kd]$,
where $a_d>0$ is the coefficient of highest degree in~$P$. By Lemma (8.1)
we have $a_d\ge\inf_{\dim Y=d}L^d\cdot Y$. We take $k=\lfloor n/d\rfloor$. The
condition $P(m)\ge\delta+1$ can thus be realized for some $m\in[m_0,m_0+kd]
\subset[m_0,m_0+n]$ as soon as
$$\inf_{\dim Y=d}L^d\cdot Y~\lfloor n/d\rfloor^d/2^{d-1}>\delta,$$
which is equivalent to the condition given in (8.4$\,$b).\qed
\endproof

\begstat{(8.5) Corollary} Let $X$ be a smooth projective $n$-fold, let
$L$ be an ample line bundle and $G$ a nef line bundle over~$X$. Then
$m(K_X+(n+2)L)+G$ is very ample for $m\ge{3n+1\choose n}-2n$.
\endstat

\begproof{} Apply Th.~(8.4$\,$a) with
$G'=a(K_X+(n+1)L)+G$, so that
$$2K_X+mL+G'=(a+2)(K_X+(n+2)L)+(m-2n-4-a)L+G,$$
and take $m=a+2n+4\ge 2+{3n+1\choose n}$.\qed
\endproof

The main drawback of the above technique is that multiples of $L$ at
least equal to $(n+1)L$ are required to avoid zeroes of the Hilbert
polynomial. In particular, it is not possible to obtain directly
a very ampleness criterion for $2K_X+L$ in the statement of (8.4$\,$b).
Nevertheless, using different ideas from Angehrn-Siu [AS94], [Siu94b] has
obtained such a criterion. We derive here a slightly weaker version,
thanks to the following elementary Lemma.

\begstat{(8.6) Lemma} Assume that for some integer $\mu\in\bbbn^\star$
the line bundle $\mu F$ gene\-rates simultaneously all jets of order
$\mu(n+s_j)+1$ at any point $x_j$ in a subset $\{x_1\ld x_p\}$ of~$X$.
Then $K_X+F$ generates simultaneously all jets of order $s_j$ at~$x_j$.
\endstat

\begproof{} Take the algebraic metric on $F$ defined by a basis of
sections $\sigma_1\ld\sigma_N$ of $\mu F$ which vanish at order
$\mu(n+s_j)+1$ at all points~$x_j$. Since we are still free to choose
the homogeneous term of degree $\mu(n+s_j)+1$ in the Taylor expansion
at~$x_j$, we find that $x_1\ld x_p$ are isolated zeroes of
$\bigcap\sigma_j^{-1}(0)$. If $\varphi$ is the weight of the metric
of $F$ near~$x_j$, we thus have $\varphi(z)\sim (n+s_j+{1\over\mu})
\log|z-x_j|$ in suitable coordinates. We replace $\varphi$ in a
neighborhood of $x_j$ by
$$\varphi'(z)=\max\big(\varphi(z)\,,\,|z|^2-C+(n+s_j)\log|z-x_j|\big)$$
and leave $\varphi$ elsewhere unchanged (this is possible by taking
$C>0$ very large). Then $\varphi'(z)=|z|^2-C+(n+s_j)\log|z-x_j|$ near
$x_j$, in particular $\varphi'$ is strictly plurisubharmonic
near~$x_j$. In this way, we get a metric $h'$ on $F$ with semipositive
curvature everywhere on $X$, and with positive definite curvature on a
neighborhood of $\{x_1\ld x_p\}$. The conclusion then follows directly
from H\"ormander's $L^2$ estimates (5.1) and (5.2).\qed
\endproof

\begstat{(8.7) Theorem} Let $X$ be a smooth projective $n$-fold, and let
$L$ be an ample line bundle over~$X$. Then $2K_X+L$ generates simultaneous
jets of order $s_1\ld s_p$ at arbitrary points $x_1\ld x_p\in X$ provided
that the intersection numbers $L^d\cdot Y$ of $L$ over all
\hbox{$d$-dimensional} algebraic subsets $Y$ of $X$ satisfy
$$L^d\cdot Y>{2^{d-1}\over\lfloor n/d\rfloor^d}
\sum_{1\le j\le p}{(n+1)(4n+2s_j+1)-2\choose n},\qquad 1\le d\le n.$$
\endstat

\begproof{} By Lemma~(8.6) applied with $F=K_X+L$ and $\mu=n+1$, the
desired jet generation of $2K_X+L$ occurs if $(n+1)(K_X+L)$ generates jets
of order \hbox{$(n+1)(n+s_j)+1$} at~$x_j$. By Lemma~(8.5) again with
\hbox{$F=aK_X+(n+1)L$} and $\mu=1$, we see by backward induction on $a$
that we need the simultaneous generation of jets of order
\hbox{$(n+1)(n+s_j)+1+(n+1-a)(n+1)$} at~$x_j$. In particular, for
$2K_X+(n+1)L$ we need the generation of jets of order
\hbox{$(n+1)(2n+s_j-1)+1$}. Theorem (8.4$\,$b) yields the desired
condition.\qed
\endproof

We now list a few immediate consequences of Theorem~8.4, in
connection with some classical questions of algebraic geometry.

\begstat{(8.8) Corollary} Let $X$ be a projective $n$-fold of
general type with $K_X$ ample. Then $mK_X$ is very ample for $m\ge
m_0={3n+1\choose n}+4$.
\endstat

\begstat{(8.9) Corollary} Let $X$ be a Fano $n$-fold, that is,
a $n$-fold such that $-K_X$ is ample. Then $-mK_X$ is very ample for
$m\ge m_0={3n+1\choose n}$.
\endstat

\begproof{} Corollaries 8.8, 8.9 follow easily from Theorem 8.4$\,$a)
applied to \hbox{$L=\pm K_X$}. Hence we get pluricanonical embeddings
$\Phi:X\to\bbbp^N$ such that
$\Phi^\star\cO(1)=\pm m_0K_X$. The image $Y=\Phi(X)$ has
degree
$$\deg(Y)=\int_Y c_1\big(\cO(1)\big)^n=
\int_X c_1\big(\pm m_0K_X\big)^n=m_0^n|K_X^n|.$$
It can be easily reproved from this that there are only finitely
many deformation types of Fano $n$-folds, as well as of $n$-folds
of general type with $K_X$ ample, corresponding to a given
discriminant $|K_X^n|$ (from a theoretical viewpoint, this result is
a consequence of Matsusaka's big theorem [Mat72] and [KoM72], but
the bounds which can be obtained from it are probably extremely huge).
In the Fano case, a fundamental result obtained indepently by
Koll\'ar-Miyaoka-Mori [KoMM92] and Campana [Cam92] shows that the
discriminant $K_X^n$ is in fact bounded by a constant $C_n$ depending
only on~$n$. Therefore, one can find an explicit bound $C'_n$ for the
degree of the embedding~$\Phi$, and it follows that there are only
finitely many families of Fano manifolds in each dimension.\qed
\endproof

In the case of surfaces, much more is known. We will content ourselves
with a brief account of recent results. If $X$ is a surface, the failure of
an adjoint bundle $K_X+L$ to be globally generated or very ample is
described in a very precise way by the following result of I.~Reider [Rei88].

\begstat{(8.10) Reider's Theorem} Let $X$ be a smooth projective surface
and let $L$ be a nef line bundle on~$X$.
\smallskip
\item{\rm a)} Assume that $L^2\ge 5$ and let $x\in X$ be a given point.
Then $K_X+L$ has a section which does not vanish at~$x$, unless there is
an effective divisor $D\subset X$ passing through $x$ such that either
$$\eqalign{
&L\cdot D=0\quad\hbox{and}\quad D^2=-1\,;\qquad\hbox{or}\cr
&L\cdot D=1\quad\hbox{and}\quad D^2=0.\cr}$$
\smallskip
\item{\rm b)} Assume that $L^2\ge 10$. Then any two points $x,y\in X$
$($possibly infinitely near$)$ are separated by sections of $K_X+L$,
unless there is an effective divisor $D\subset X$ passing through $x$
and $y$ such that either
$$\eqalignno{
&L\cdot D=0\quad\hbox{and}\quad D^2=-1~\hbox{or}~{}-2\,;
\qquad\hbox{or}&\cr
&L\cdot D=1\quad\hbox{and}\quad D^2=0\phantom{-{}}~\hbox{or}~{}-1\,;
\qquad\hbox{or}&\cr
&L\cdot D=2\quad\hbox{and}\quad D^2=0.&\square\cr}$$
\vskip0pt
\endstat

\begstat{(8.11) Corollary} Let $L$ be an ample line bundle on a smooth
projective surface~$X$. Then $K_X+3L$ is globally generated and $K_X+4L$
is very ample. If $L^2\ge 2$ then $K_X+2L$ is globally generated and
$K_X+3L$ is very ample.\qed
\endstat

The case of higher order jets can be treated similarly. The most general
result in this direction has been obtained by Beltrametti and Sommese
[BeS93].

\begstat{(8.12) Theorem {\rm([BeS93])}} Let $X$ be a smooth projective surface
and let $L$ be a nef line bundle on~$X$. Let $p$ be a positive integer
such that $L^2>4p$. Then for every $0$-dimensional subscheme $Z\subset X$
of length $h^0(Z,\cO_Z)\le p$ the restriction
$$\rho_Z:H^0(X,\cO_X(K_X+L))\lra H^0(Z,\cO_Z(K_X+L))$$
is surjective, unless there is an effective divisor $D\subset X$
intersecting the support $|Z|$ such that
$$L\cdot D-p\le D^2<{1\over 2}L\cdot D.\eqno\square$$
\endstat

\begproof{(Sketch).} The proof the above theorems rests in an essential way
on the construction of rank 2 vector bundles sitting in an exact sequence
$$0\to \cO_X\to E\to L\otimes\cI_Z\to 0.$$
Arguing by induction on the length of $Z$, we may assume that $Z$ is a
$0$-dimensional subscheme such that $\rho_Z$ is not surjective, but such
that $\rho_{Z'}$ is surjective for every proper subscheme $Z'\subset Z$.
The existence of $E$ is obtained by a classical construction of Serre
(unfortunately, this construction only works in dimension $2$).
The numerical condition on $L^2$ in the hypotheses ensures that
$c_1(E)^2-4\,c_2(E)>0$, hence $E$ is unstable in the sense of
Bogomolov. The existence of the effective divisor $D$ asserted in
8.10 or 8.12 follows. We refer to [Rei88], [BeS93] and [Laz93] for details.
The reader will find in [FdB93] a proof of the Bogomolov inequality
depending only on the Kawamata-Viehweg vanishing theorem.\qed
\endproof

\begstat{(8.13) Exercise} \rm The goal of the exercise is to prove the
following weaker form of Theorems 8.10 and 8.12, by a simple direct method
based on Nadel's vanishing theorem:
\smallskip
\item{}{\it Let $L$ be a nef line bundle on a smooth projective surface~$X$.
Fix points $x_1\ld x_N$ and corresponding multiplicities $s_1\ld s_N$, and
set $p=\sum (2+s_j)^2$. Then $H^0(X,K_X+L)$ generates simultaneously
jets of order $s_j$ at all points $x_j$ provided that
$L^2>p$ and $L\cdot C>p$ for all curves $C$ passing through one of the
points~$x_j$.}
\smallskip
\item{\rm a)} Using the Riemann-Roch formula, show that the condition $L^2>p$
implies the existence of a section of a large multiple $mL$ vanishing
at order ${}>m(2+s_j)$ at each of the points.
\item{\rm b)} Construct a sequence of singular hermitian metrics on $L$
with positive definite curvature, such that the weights $\varphi_\nu$ have
algebraic singularities, $\nu(\varphi_\nu,x_j)\ge 2+s_j$ at each point,
and such that for some integer $m_1>0$ the multiplier ideal sheaves
satisfy $\cI(m_1\varphi_{\nu+1})\supsetneq\cI(m_1\varphi_\nu)$ if
$V(\cI(\varphi_\nu))$ is not $0$-dimensional near some $x_j$.
\smallskip\noindent
{\it Hint\/}: a) starts the procedure. Fix $m_0>0$ such that $m_0L-K_X$ is
ample. Use Nadel's vanishing theorem to show that
$$H^q(X,\cO((m+m_0)L)\otimes\cI(\lambda m\varphi_\nu))=0
\qquad\hbox{\rm for all $q\ge 1$, $m\ge 0$, $\lambda\in[0,1]$.}$$
Let $D_\nu$ be the effective $\bbbq$-divisor describing the 1-dimensional
singularities of $\varphi_\nu$. Then $\cI(\lambda m\varphi_\nu)
\subset\cO(-\lfloor\lambda mD_\nu\rfloor)$ and the quotient has
$0$-dimensional support, hence
$$H^q(X,\cO((m+m_0)L)\otimes\cO(-\lfloor\lambda m D_\nu\rfloor))=0
\qquad\hbox{\rm for all $q\ge 1$, $m\ge 0$, $\lambda\in[0,1]$.}$$
By Riemann-Roch again prove that
$$h^0(X,\cO((m+m_0)L)\otimes\cO/\cO(-\lfloor\lambda
mD_\nu\rfloor))={m^2\over 2}(2\lambda\, L\cdot D_\nu
-\lambda^2 D_\nu^2)+O(m).\leqno(\star)$$
As the left hand side of $(\star)$ is increasing with $\lambda$, one
must have $D_\nu^2\le L\cdot D_\nu$. If~$V(\cI(\varphi_\nu))$ is not
$0$-dimensional at $x_j$, then the coeffi\-cient of some component of
$D_\nu$ passing through $x_j$ is at least~$1$, hence
$$2\, L\cdot D_\nu -D_\nu^2\ge L\cdot D_\nu\ge p+1.$$
Show the existence of an integer $m_1>0$ independent of $\nu$ such that
$$h^0(X,\cO((m+m_0)L)\otimes\cO/\cO(-\lfloor mD_\nu\rfloor))>
\sum_{1\le j\le N}{(m+m_0)(2+s_j)+2\choose 2}$$
for $m\ge m_1$, and infer the existence of a suitable section of $(m_1+
m_0)L$ which is not in $H^0(X,\cO((m_1+m_0)L-\lfloor m_1D_\nu\rfloor))$.
Use this section to construct $\varphi_{\nu+1}$ such that
$\cI(m_1\varphi_{\nu+1})\supsetneq\cI(m_1\varphi_\nu)$.
\endstat

\titleb{9.}{Regularization of Currents\newline
and Self-intersection Inequalities}
Let $X$ be a compact complex $n$-dimensional manifold. It will be
convenient to work with currents that are not necessarily positive,
but such that their negative part is locally bounded.
We say that a bidimension $(p,p)$ current $T$ is {\it almost
positive} if there exists a smooth form $v$ of bidegree $(n-p,n-p)$
such that $T+v\ge 0$. Similarly, a function $\varphi$ on $X$ is said to
be {\it almost psh} if $\varphi$ is locally equal to the sum of a psh
function and of a smooth function; then the $(1,1)$-current
$dd^c\varphi$ is almost positive; conversely, if a locally integrable
function $\varphi$ is such that $dd^c\varphi$ is almost positive, then
$\varphi$ is equal a.e.\ to an almost psh function. If $T$ is closed
and almost positive, the Lelong numbers  $\nu(T,x)$ are well defined,
since the negative part always contributes to zero.

\begstat{(9.1) Theorem} Let $T$ be a closed almost positive $(1,1)$-current
and let $\alpha$ be a smooth real $(1,1)$-form in the the same
$dd^c$-cohomology class as $T$, i.e. $T=\alpha+dd^c\psi$ where $\psi$
is an almost psh function. Let $\gamma$ be a continuous real $(1,1)$-form
such that $T\ge\gamma$. Suppose that $\cO_{T_X}(1)$ is equipped with a
smooth hermitian metric such that the curvature form satisfies
$${\ii\over 2\pi}\Theta\big(\cO_{T_X}(1)\big)+\pi^\star u\ge 0$$
with $\pi:P(T^\star X)\to X$ and with some nonnegative smooth $(1,1)$-form
$u$ on~$X$. Fix a hermitian metric $\omega$ on~$X$. Then for every~$c>0$,
there is a sequence of closed almost positive $(1,1)$-currents
$T_{c,k}=\alpha+dd^c\psi_{c,k}$ such that $\psi_{c,k}$
is smooth on $X\ssm E_c(T)$ and decreases to $\psi$ as $k$ tends to~$+\infty$
$($in particular, $T_{c,k}$ is smooth on $X\ssm E_c(T)$ and converges weakly
to $T$ on~$X)$, and
$$T_{c,k}\ge\gamma-\lambda_{c,k}u-\varepsilon_k\omega$$
where
\smallskip
\item{\rm a)} $\lambda_{c,k}(x)$ is a decreasing sequence of continuous
functions on $X$ such that\newline $\lim_{k\to+\infty}\lambda_{c,k}(x)=
\min\big(\nu(T,x),c\big)$ at every point,
\smallskip
\item{\rm b)} $\lim_{k\to+\infty}\varepsilon_k=0$,
\smallskip
\item{\rm c)} $\nu(T_{c,k},x)=\big(\nu(T,x)-c\big)_+$
at every point $x\in X$.
\endstat

Here $\cO_{T_X}(1)$ is the canonical line bundle associated with $T_X$ over
the hyperplane bundle~$P(T^\star X)$. Observe that the theorem gives in
particular approximants $T_{c,k}$ which are smooth everywhere on $X$ if
$c$ is taken such that $c>\max_{x\in X}\nu(T,x)$. The equality
in~c) means that the procedure kills all Lelong numbers that are
$\le c$ and shifts all others downwards by $c$. Hence Theorem~9.1 is an
analogue over manifolds of Kiselman's procedure [Kis78~,79] for killing
Lelong numbers of a psh function on an open subset of~$\bbbc^n$.

\begproof{(Sketch).} We refer to [Dem92] for a detailed proof. We only sketch
a special case for which the main idea is simple to explain. The special
case we wish to consider is the following: $X$ is projective algebraic,
$u={\ii\over2\pi}\Theta(G)$ is the curvature form of a nef
$\bbbq$-divisor $G$, and $T$ has the form
$$T=\alpha+dd^c\psi,\qquad\psi=\log(|f_1^2|+\ldots+|f_N|^2\big)^{1/2},$$
where $\alpha$ is smooth and the $f_j$'s are sections of some $C^\infty$
hermitian line bundle $L$ on~$X$. The lower bound for $T$ is
$$T\ge\gamma:=\alpha-{\ii\over2\pi}\Theta(L),$$
hence $T$ is almost positive. Somehow, the general situation can be
reduced locally to this one by an approximation theorem for currents
based on H\"ormander's $L^2$ estimates, in the form given by
Ohsawa-Takegoshi ([OT87], [Ohs88]); the main point is to show that any
closed positive $(1,1)$-current is locally a weak limit of effective
$\bbbq$-divisors which have roughly the same Lelong numbers as the
given current, up to small errors converging to~$0$; the proof is then
completed by means of rather tricky gluing techniques for psh functions
(see [Dem92]). Now, let $A$ be an ample divisor and let
$\omega={\ii\over 2\pi}\Theta(A)$ be a positive curvature form for~$A$.
After adding $\varepsilon A$ to $G$ ($\varepsilon\in\bbbq^+$), which is
the same as adding $\varepsilon\omega$ to $u$, we may assume that $u$
is positive definite and that $\cO_{T_X}(1)+\pi^\star G$ is ample
($\cO_{T_X}(1)$ is relatively ample, so adding something ample from $X$
is enough to make it ample on $P(T^\star_X)$). The Lelong numbers of
$T$ are given by the simple formula
$$\nu(T,x)=\min_{1\le i\le N}\ord_x(f_i).$$
The basic idea is to decrease the Lelong numbers of $T$ by adding
some ``derivatives'' of the $f_i$'s in the sum of squares. The derivatives
have of course to be computed as global sections on~$X$. For this, we
introduce two positive integers $m,p$ which will be selected later more
carefully, and we consider the $m$-jet sections $J^mf_i^{\otimes p}$, viewed
as sections of the $m$-jet vector bundle $J^mL^{\otimes p}$. The bundle
$J^mL^{\otimes p}$ has a filtration whose graded terms are $S^\nu T^\star_X
\otimes L^{\otimes p}$, $0\le\nu\le m$, hence $(J^mL^{\otimes p})^\star
\otimes L^{\otimes p}$ has a filtration with graded terms $S^\nu T_X$,
$0\le\nu\le m$. Our hypothesis $T_X\otimes\cO(G)$ ample implies that
$(J^mL^{\otimes p})^\star\otimes\cO(pL+mG)$ is ample.
Hence, raising to some symmetric power $S^\nu$ where $\nu$ is a multiple
of a denominator for the $\bbbq$-divisor $G$, we obtain that
$S^\nu(J^mL^{\otimes p})^\star\otimes\cO(\nu pL+\nu mG)$ is generated
by global sections $P_1\ld P_{N'}$ for $\nu$ large (we view the $P_j$'s
as some kind of ``differential operators'' acting on sections of~$L$).
Now, the $\nu$-th symmetric product $S^\nu(J^mf_i^{\otimes p})$ takes
values in $S^\nu(J^mL^{\otimes p})$, hence the pairing with the dual
bundle yields a section
$$F^{p,\nu}_{m,i,j}:=P_j\cdot S^\nu(J^mf_i^{\otimes p})\in
H^0(X,\cO(\nu pL+\nu mG)).$$
Since the $P_j$ generate sections and $J^mf_i^{\otimes p}$ vanishes at
order $(p\,\ord_x(f_i)-m)_+$ at any point~$x$, we get
$$\min_j\ord_x(F^{p,\nu}_{m,i,j})=\nu(p\,\ord_x(f_i)-m)_+.$$
We set
$$\eqalign{
T_{c,k}&=\alpha+dd^c\psi_{c,k},\cr
\psi_{c,k}&={1\over 2\nu p}\log\Big(\big(\sum_i|f_i|^2\big)^{\nu p}+
\sum_{1\le\mu\le m}k^{-\mu}
\sum_{i,j}\big|F^{p,\nu}_{\mu,i,j}\big|^2\Big)\cr}$$
Clearly, $\psi_{c,k}$ is decreasing and converges to $\psi$, and its
Lelong numbers are
$$\nu(\psi_{c,k},x)=\min_i\Big(\ord_x(f_i)-{m\over p}\Big)_+
=\Big(\nu(\psi,x)-{m\over p}\Big)_+.$$
Select $c={m\over p}$ (or take a very close rational approximation
if an arbitrary real number $c\in\bbbr^+$ is to be considered). As
$F^{p,\nu}_{m,i,j}$ takes values in $\cO(\nu pL+\nu mG)$, an easy
computation yields
$$T_{c,k}\ge\alpha-{\ii\over 2\pi}\Theta\Big(L+{m\over p}G\Big)=\gamma-
{m\over p}u=\gamma-cu.$$
However, at points where $\nu(T,x)<c$, there is already a term with
${\mu\over p}\simeq\nu(T,x)$ such that the contribution of the terms
indexed by $\mu$ give a non zero contribution. The terms corresponding to
higher indices will then not contribute much to the lower bound of
$dd^c\psi_{c,k}$ since in the summation $k^{-\mu+1}$ becomes negligible
against $k^{-\mu}$ as $k\to+\infty$. This gives a
lower bound $T_{c,k}\ge\gamma-\lambda_{c,k}-\varepsilon_k\omega$
of the expected form.\qed
\endproof

\begstat{(9.2) Corollary} Let $\Theta$ be a closed almost positive current
of bidimension $(p,p)$ and let $\alpha_1\ld\alpha_q$ be closed almost
positive $(1,1)$-currents such that $\alpha_1\wedge\ldots\wedge\alpha_q
\wedge \Theta$ is well defined by application of criteria~$2.3$ or $2.5$,
when $\alpha_j$ is written locally as~$\alpha_j=dd^cu_j$. Then
$$\{\alpha_1\wedge\ldots\wedge\alpha_q\wedge \Theta\}=
\{\alpha_1\}\cdots\{\alpha_q\}\cdot\{\Theta\}.$$
\endstat

\begproof{} Theorem~9.1 and the monotone continuity
theorem for Monge-Amp\`ere operators show that
$$\alpha_1\wedge\ldots\wedge\alpha_q\wedge \Theta=\lim_{k\to+\infty}
\alpha_1^k\wedge\ldots\wedge\alpha_q^k\wedge \Theta$$
where $\alpha_j^k\in\{\alpha_j\}$ is smooth. Since the result is
by definition true for smooth forms, we conclude by the weak
continuity of cohomology class assignment.\qed
\endproof

Now, let $X$ be a compact K\"ahler manifold equipped with a K\"ahler
metric~$\omega$. The {\it degree} of a closed positive current $\Theta$
with respect to $\omega$ is by definition
$$\deg_\omega\Theta=\int_X\Theta\wedge\omega^p,~~~~
\hbox{\rm bidim}\,\Theta=(p,p).\leqno(9.3)$$
In particular, the degree of a $p$-dimensional analytic set $A\subset X$
is its volume $\int_A\omega^p$ with respect to~$\omega$. We are
interested in the following problem.

\begstat{(9.4) Problem} Let $T$ be a closed positive $(1,1)$-current
on~$X$. Is it possible to derive a bound for the codimension~$p$
components in the Lelong upperlevel sets $E_c(T)$ in terms of the
cohomology class $\{T\}\in H^2_{DR}(X,\bbbr)\,{\rm ?}$
\endstat

Let $\Xi\subset X$ be an arbitrary subset. We introduce the sequence
$$0=b_1\le\ldots\le b_n\le b_{n+1}$$
of ``{\it jumping values}" of $E_c(T)$ over $\Xi$, defined by
the property that the dimension of $E_c(T)$ in a neighborhood of $\Xi$
drops by one unit when $c$ gets larger than~$b_p$, namely
$$b_p=\inf\big\{c>0\,;\,\codim\big(E_c(T),x\big)\ge p,~\forall x\in\Xi\big\}.$$
Then, when $c\in{}]b_p,b_{p+1}]$, we have $\codim E_c(T)=p$
in a neighborhood of $\Xi$, with at least one component of codimension $p$
meeting $\Xi$. Let $(Z_{p,k})_{k\ge 1}$ be the family of all these
$p$-codimensional components $\big($occurring in any of the sets $E_c(T)$
for $c\in{}]b_p,b_{p+1}]\big)$, and let
$$\nu_{p,k}=\min_{x\in Z_{p,k}}\nu(T,x)\in{}]b_p,b_{p+1}]$$
be the generic Lelong number of $T$ along~$Z_{p,k}$. Then we have
the following self-intersection inequality.

\begstat{(9.5) Theorem} Suppose that $X$ is K\"ahler and that
$\cO_{T_X}(1)$ has a hermitian metric such that ${\ii\over 2\pi}
\Theta\big(\cO_{T_X}(1)\big)+\pi^\star u \ge 0$, where $u$ is a
smooth closed semipositive $(1,1)$-form. Let $\Xi \subset X$ be an
arbitrary subset, let $T$ be a closed positive current of bidegree
$(1,1)$, and let $(b_p)$, $(Z_{p,k}$ be the corresponding jumping
values and $p$-codimensional components of $E_c(T)$ meeting~$\Xi$.
Assume either that $\Xi=X$ or that the cohomology class $\{T\}\in
H^{1,1}(X)$ is nef $($i.e.\ in the closure of the K\"ahler cone$)$.
Then, for each $p=1\ld n$, the De Rham cohomology class
$(\{T\}+b_1\{u\})\cdots(\{T\}+b_p\{u\})$ can be represented by a closed
positive current $\Theta_p$ of bidegree $(p,p)$ such that
$$\Theta_p\ge\sum_{k\ge 1}(\nu_{p,k}-b_1)\ldots(\nu_{p,k}-b_p)\,[Z_{p,k}]
+(T_\ac+b_1u)\wedge\ldots\wedge(T_\ac+b_pu)$$
where $T_\ac\ge 0$ is the absolutely continuous part in the
Lebesgue decomposition of $T$ $($decomposition of the coefficients of $T$
into absolutely continuous and singular mea\-sures$)$, $T=T_\ac+T_\sing$.
\endstat

By neglecting the second term in the right hand side and taking the
wedge product with $\omega^{n-p}$, we get the following interesting
consequence:

\begstat{(9.6) Corollary} If $\omega$ is a K\"ahler metric on $X$ and if
$\{u\}$ is a nef cohomology class such that $c_1\big(\cO_{T_X}(1)
\big)+\pi^\star\{u\}$ is nef, the degrees of the components
$Z_{p,k}$ with respect to $\omega$ satisfy the estimate
$$\eqalign{
\sum_{k=1}^{+\infty}(\nu_{p,k}-b_1)\ldots(\nu_{p,k}-b_p)
&\int_X[Z_{p,k}]\wedge\omega^{n-p}\cr
&\le\big(\{T\}+b_1\{u\}\big)\cdots\big(\{T\}+b_p\{u\}\big)
\cdot\{\omega\}^{n-p}.\cr}$$
\endstat

As a special case, if $D$ is an effective divisor and $T=[D]$, we get a
bound for the degrees of the $p$-codimensional singular strata of $D$ in
terms of a polynomial of degree $p$ in the cohomology class~$\{D\}$; the
multiplicities \hbox{$(\nu_{p,k}-b_1)\ldots(\nu_{p,k}-b_p)$} are then positive
integers. The case when $X$ is $\bbbp^n$ or a homogeneous manifold
is especially simple: then $T_X$ is generated by sections and
we can take $u=0\,$; the bound is thus simply $\{D\}^p\cdot
\{\omega\}^{n-p}$; the same is true more generally as soon as $\cO_{T_X}(1)$
is nef. The main idea of the proof is to kill the Lelong numbers of $T$ up
to the level~$b_j\,$; then the singularities of the resulting current $T_j$
occur only in codimension $j$ and it becomes possible to define the wedge
product $T_1\wedge\ldots\wedge T_p$ by means of Proposition 2.3. Here are
the details:

\begproof{of Theorem 9.5.} First suppose that $\Xi=X$. We argue by
induction on~$p$. For $p=1$, Siu's decomposition formula shows that
$$T=\sum\nu_{1,k}[Z_{1,k}]+R,$$
and we have $R\ge T_\ac$ since the other part has singular measures
as coefficients. The result is thus true with $\Theta_1=T$. Now, suppose
that $\Theta_{p-1}$ has been constructed. For $c>b_p$, the current
$T_{c,k}=\alpha+dd^c\psi_{c,k}$
produced by Theorem~9.1 is such that the codimension of the set of
poles $\psi_{c,k}^{-1}(-\infty)=E_c(T)$ is at least~$p$ at every point
$x\in X$ (recall that $\Xi=X$). Then Proposition 2.5 shows that
$$\Theta_{p,c,k}=\Theta_{p-1}\wedge(T_{c,k}+c\,u+\varepsilon_k\omega)$$
is well defined. If $\varepsilon_k$ tends to zero slowly enough,
$T_{c,k}+c\,u+\varepsilon_k\omega$ is positive by (9.1$\,$a), so
$\Theta_{p,c,k}\ge 0$. Moreover, by Corollary~9.2, the
cohomology class of $\Theta_{p,c,k}$ is $\{\Theta_{p-1}\}\cdot
(\{T\}+c\{u\}+\varepsilon_k\{\omega\})$, converging to
$\{\Theta_{p-1}\}\cdot(\{T\}+c\{u\})$. Since the mass
$\int_X\Theta_{p,c,k}\wedge\omega^{n-p}$ remains uniformly bounded,
the family $(\Theta_{p,c,k})_{c\in{}]b_p,b_p+1],k\ge 1}$ is relatively
compact in the weak topology. We define
$$\Theta_p=\lim_{c\to b_p+0}\lim_{k\to+\infty}\Theta_{p,c,k},$$
possibly after extracting some weakly convergent subsequence. Then
$\{\Theta_p\}=\{\Theta_{p-1}\}\cdot(\{T\}+b_p\{u\})$, and so
$\{\Theta_p\}=(\{T\}+b_1\{u\})\cdots(\{T\}+b_p\{u\})$. Moreover,
it is well-known (and easy to check) that Lelong numbers are upper
semi-continuous with respect to weak limits of currents. Therefore
$$\eqalign{
\nu(\Theta_p,x)&\ge\limsup_{c\to b_p+0}\limsup_{k\to+\infty}~
\nu\big(\Theta_{p-1}\wedge(T_{c,k}+c\,u+\varepsilon_k\omega),x\big)\cr
&\ge\nu(\Theta_{p-1},x)\times
\limsup_{c\to b_p+0}\limsup_{k\to+\infty}~\nu(T_{c,k},x)\cr
&\ge\nu(\Theta_{p-1},x)\big(\nu(T,x)-b_p\big)_+\cr}$$
by application of Proposition 2.16 and (9.1$\,$a). Hence by induction
we get
$$\nu(\Theta_p,x)\ge
\big(\nu(T,x)-b_1\big)_+\ldots\big(\nu(T,x)-b_p\big)_+,$$
in particular, the generic Lelong number of $\Theta_p$ along $Z_{p,k}$
is at least equal to the product $(\nu_{p,k}-b_1)\ldots(\nu_{p,k}-b_p)$.
This already implies
$$\Theta_p\ge\sum_{k\ge 1}(\nu_{p,k}-b_1)\ldots(\nu_{p,k}-b_p)\,[Z_{p,k}].$$
Since the right hand side is Lebesgue singular, the desired inequality
will be proved if we show in addition that
$$\Theta_{p,\ac}\ge(T_\ac+b_1u)\wedge\ldots\wedge(T_\ac+b_pu),$$
or inductively, that $\Theta_{p,\ac}\ge\Theta_{p-1,\ac}
\wedge (T_\ac+b_pu)$. In order to do this, we simply have to
make sure that $\lim_{k\to+\infty}T_{c,k,\ac}=T_\ac$ almost
everywhere and use induction again. But our arguments are not affected
if we replace $T_{c,k}$ by $T'_{c,k}=\alpha+dd^c\psi'_{c,k}$ where
$\psi'_{c,k}=\max\{\psi,\psi_{c,k}-A_k\}$ and $(A_k)$ is a sequence
converging quickly to~$+\infty$. Lemma 9.7 below shows that a
suitable choice of $A_k$ gives $\lim (dd^c\psi'_{c,k})_\ac
=(dd^c\psi)_\ac$ almost everywhere. This concludes the proof in the case
$\Xi=X$.

When $\Xi\ne X$, a slight difficulty appears: in fact, there may remain in
$T_{c,k}$ some poles of codimension $\le p-1$, corresponding to components
of $E_c(T)$ which do not meet $\Xi$ (since we completely
forgot these components in the definition of the jumping values). It follows
that the wedge product $\Theta_{p-1}\wedge(T'_{c,k}+c\,u+\varepsilon_k\omega)$
is no longer well defined. In this case, we proceed as follows.
The assumption that $T$ is nef implies that are smooth functions
$\wt\psi_k$ such that the cohomology class $\{T\}$ has a representative
$\alpha+dd^c\wt\psi_k\ge-\varepsilon_k\omega$. We replace $T'_{c,k}$ in
the above arguments by $T'_{c,k,\nu}=\alpha+dd^c\psi'_{c,k,\nu}$ where
$$\psi'_{c,k,\nu}=\max\{\psi,\psi_{c,k}-A_k,\wt\psi_k-\nu\}.$$
Then certainly $T'_{c,k,\nu}+c\,u+\varepsilon_k\omega\ge0$ and we can
define a closed positive current $\Theta_{p-1}\wedge(T'_{c,k,\nu}+
c\,u+\varepsilon_k\omega)$ without any difficulty since $\psi'_{c,k,\nu}$
is locally bounded. We first extract a weak limit $\Theta_{p,c,k}$ as
$\nu\to+\infty$. By monotone continuity of Monge-Amp\`ere operators, we
find
$$\Theta_{p,c,k}=\Theta_{p-1}\wedge(T'_{c,k}+c\,u+\varepsilon_k\omega)$$
in the neighborhood of $\Xi$ where this product is well defined. All other
arguments are the same as before.
\endproof

\begstat{(9.7) Lemma} Let $\Omega\subset\bbbc^n$ be an open subset and
let $\varphi$ be an arbitrary psh function on $\Omega$. Set
$\varphi_\nu=\max(\varphi,\psi_\nu)$ where $\psi_\nu$ is a decreasing
sequence of psh functions converging to $-\infty$, each
$\psi_\nu$ being locally bounded in $\Omega$ $($or perhaps only
in the complement of an analytic subset of codimension~$\ge p)$.
Let $\Theta$ be a closed positive current of bidegree $(p-1,p-1)$.
If $\Theta\wedge dd^c\varphi_\nu$ converges to a weak
limit~$\Theta'$, then
$$\Theta'_\ac\ge\Theta_\ac\wedge (dd^c\varphi)_\ac.$$
\endstat

\begproof{} Let $(\rho_\varepsilon)$ (resp.\ $(\wt\rho_\varepsilon)$)
be a family of regularizing kernels on $\bbbc^n$ (resp.\ on $\bbbr^2$), and
let $\max_\varepsilon(x,y)=({\max}\star\wt\rho_\varepsilon)(x,y)$ be
a regularized max function. For $\varepsilon>0$ small enough, the function
$$\varphi_{\nu,\varepsilon}={\textstyle\max_\varepsilon}
(\varphi\star\rho_\varepsilon,\psi_\nu\star\rho_\varepsilon)$$
is psh and well defined on any preassigned open set
$\Omega'\compact\Omega$. As $\varphi_{\nu,\varepsilon}$ decreases to
$\varphi_\nu$ when $\varepsilon$ decreases to $0$, proposition~10.2
shows that
$$\lim_{\varepsilon\to 0}\Theta\wedge dd^c\varphi_{\nu,\varepsilon}
=\Theta\wedge dd^c\varphi_\nu$$
in the weak topology. Let $(\beta_j)$ be a sequence of test forms
which is dense in the space of test forms of bidegree $(n-p,n-p)$ and
contains strongly positive forms with arbitrary large compact support
in $\Omega$. Select $\varepsilon_\nu>0$ so small that
$$\langle\Theta\wedge dd^c\varphi_{\nu,\varepsilon_\nu}
-\Theta\wedge dd^c\varphi_\nu\,,\,\beta_j\rangle\le{1\over\nu}
{}~~~~\hbox{\rm for}~~j\le\nu.$$
Then the sequence $\Theta\wedge dd^c\varphi_{\nu,\varepsilon_\nu}$
is locally uniformly bounded in mass and converges weakly to the same limit
$\Theta'$ as $\Theta\wedge dd^c\varphi_\nu$. Moreover, at every
point $x\in\Omega$ such that $\varphi(x)>-\infty$, we have
$\varphi_{\nu,\varepsilon_\nu}(x)\ge\varphi(x)>\psi_\nu\star
\rho_{\varepsilon_\nu}(x)+1$ for $\nu$ large, because
$\lim_{\nu\to-\infty}\psi_\nu=-\infty$ locally uniformly. Hence
$\varphi_{\nu,\varepsilon_\nu}=\varphi\star\rho_{\varepsilon_\nu}$
on a neighborhood of $x$ (which may depend on $\nu$) and
$dd^c\varphi_{\nu,\varepsilon_\nu}(x)=(dd^c\varphi)
\star\rho_{\varepsilon_\nu}(x)$ for $\nu\ge\nu(x)$. By the Lebesgue
density theorem, if $\mu$ is a measure of absolutely continuous part
$\mu_\ac$, the sequence $\mu\star\rho_{\varepsilon_\nu}(x)$ converges
to $\mu_\ac(x)$ at almost every point. Therefore
$\lim dd^c\varphi_{\nu,\varepsilon_\nu}(x)=
(dd^c\varphi)_\ac(x)$ almost everywhere For any strongly positive test
form $\alpha=i\alpha_1\wedge\ovl\alpha_1\wedge\ldots\wedge i\alpha_{n-p}
\wedge\ovl\alpha_{n-p}$ of bidegree $(n-p,n-p)$ on $\Omega$, we get
$$\eqalign{
\int_\Omega\Theta'\wedge\alpha&=\lim_{\nu\to+\infty}\int_\Omega\Theta
\wedge dd^c\varphi_{\nu,\varepsilon_\nu}\wedge\alpha\cr
&\ge\liminf_{\nu\to+\infty}\int_\Omega\Theta_\ac
\wedge dd^c\varphi_{\nu,\varepsilon_\nu}\wedge\alpha
\ge\int_\Omega\Theta_\ac\wedge dd^c\varphi_\ac
\wedge\alpha.\cr}$$
Indeed, the first inequality holds because $dd^c
\varphi_{\nu,\varepsilon_\nu}$ is smooth, and the last one results
from Fatou's lemma. This implies
$\Theta'_\ac\ge\Theta_\ac\wedge(dd^c\varphi)_\ac$
and Lemma~9.7 follows.\qed
\endproof

\titleb{10.}{Use of Monge-Amp\`ere Equations}
The goal of the next two sections is to find numerical criteria for an
adjoint line bundle $K_X+L$ to be generated by sections (resp.\ very
ample, $s$-jet ample). The conditions ensuring these conclusions should
be ideally expressed in terms of explicit lower bounds for the
intersection numbers $L^p\cdot Y$, where $Y$ runs over $p$-dimensional
subvarieties of $X$ (as the form of Reider's theorem suggests in the case
of surfaces). Unfortunately, the simple algebraic approach described in
the proof of Theorems 8.4 and 8.5 does not seem to be applicable to get
criteria for the very ampleness of $K_X+mL$ in the range $m\le n+1$.
We will now explain how this can be achieved by an alternative analytic
method. The essential idea is to construct directly the psh weight
function $\varphi$ needed in Nadel's vanishing theorem by solving a
Monge-Amp\`ere equation (Aubin-Calabi-Yau theorem). Quite recently,
Ein and Lazarsfeld [EL94] have developed purely algebraic methods
which yield similar results; however, up to now, the bounds obtained
with the algebraic method are not as good as with the analytic approach.

Let us first recall a special case of the well-known theorem of Aubin-Yau
related to the Calabi conjecture.  The special case we need is the following
fundamental existence result about solutions of Monge-Amp\`ere equations.

\begstat{(10.1) Theorem {\rm([Yau78], see also [Aub78])}} Let $X$ be a compact
complex $n$-dimen\-sional manifold with a smooth K\"ahler metric $\omega$.
Then for any smooth volume form $f>0$ with $\int_X f=\int_X \omega^n$, there
exists a unique K\"ahler metric $\wt\omega=\omega+dd^c\psi$ in the
K\"ahler class $\{\omega\}$ such that $\wt\omega^n=(\omega+dd^c\psi)^n=f$.
\endstat

There are several equivalent ways of formulating this result. Usually, one
starts with a $(1,1)$-form $\gamma$ representing the first Chern class
$c_1(X)=c_1(\Lambda^nT^\star_X)$. Then, under the normalization
$\int_X f=\int_X\omega^n$, there is a unique volume form $f$ on $X$
which, viewed as a hermitian form on $\Lambda^nT_X$, yields
${\ii\over2\pi}\Theta_f(\Lambda^nT_X)=\gamma$. Then Theorem 10.1 is actually
equivalent to finding a K\"ahler metric $\wt\omega$ in the K\"ahler class
$\{\omega\}$, such that Ricci$(\wt\omega)=\gamma$. We will not use this
viewpoint here, and will be essentially concerned instead with the
Monge-Amp\`ere equation $(\omega+dd^c\psi)^n=f$.

There are two different ways in which the Monge-Amp\`ere equation will be
used. The most essential idea is that the Monge-Amp\`ere equation can be
used to produce weights with logarithmic singularities (as needed for
the application of Corollary~5.13), when the right hand side $f$ is taken
to be a linear combination of Dirac measures (in fact, $f$ has to be
smooth so we rather find solutions $\psi_\varepsilon$ corresponding
to smooth approximations $f_\varepsilon$ of the Dirac measures). This
will be explained later.

Another useful consequence of the Monge-Amp\`ere equation is a general
version of convexity inequality due to Hovanski [Hov79] and Teissier
[Tei79,~82], which is a natural generalization of the usual Hodge
index theorem for surfaces. This inequality is reproved along similar
lines in [BBS89], where it is applied to the study of projective
$n$-folds of log-general type. For the sake of completeness, we include
here a different and slightly simpler proof, based on
the Aubin-Yau theorem~10.1 instead of the Hodge index theorem.
Our proof also has the (relatively minor) advantage of working over
arbitrary  K\"ahler manifolds.

\begstat{(10.2) Proposition} The following inequalities hold in any
dimension~$n$.
\smallskip
\item{\rm a)} If $\alpha_1\ld\alpha_n$ are semipositive $(1,1)$-forms
on $\bbbc^n$, then
$$\alpha_1\wedge\alpha_2\wedge\ldots\wedge\alpha_n\ge
(\alpha_1^n)^{1/n}(\alpha_2^n)^{1/n}\ldots(\alpha_n^n)^{1/n}.$$
\item{\rm b)} If $u_1\ld u_n$ are nef cohomology classes of type
$(1,1)$ on a K\"ahler manifold $X$ of dimension $n$, then
$$u_1\cdot u_2\cdots u_n\ge(u_1^n)^{1/n}(u_2^n)^{1/n}\ldots(u_n^n)^{1/n}.$$
\vskip0pt
\endstat

By a nef cohomology class of type $(1,1)$, we mean a class in
the closed convex cone of $H^{1,1}(X,\bbbr)$ generated by K\"ahler classes,
that is, a class $\{u\}$ admitting representatives $u_\varepsilon$ with
$u_\varepsilon\ge-\varepsilon\omega$ for every $\varepsilon>0$.
For instance, inequality b) can be applied to $u_j=c_1(L_j)$ when
$L_1\ld L_n$ are nef line bundles over a projective manifold.

\begproof{} Observe that a) is a pointwise inequality between
$(n,n)$-forms whereas b) is an inequality of a global nature for the cup
product intersection form. We first show that a) holds when only two of
the forms $\alpha_j$ are distinct, namely that
$$\alpha^p\wedge\beta^{n-p}\ge(\alpha^n)^{p/n}(\beta^n)^{(n-p)/n}$$
for all $\alpha,\beta\ge 0$. By a density argument, we may suppose
$\alpha,\beta>0$. Then there is a simultaneous orthogonal basis in which
$$\alpha=i\sum_{1\le j\le n}\lambda_j\,dz_j\wedge d\ovl z_j,~~~~
\beta=i\sum_{1\le j\le n}dz_j\wedge d\ovl z_j$$
with $\lambda_j>0$, and a) is equivalent to
$$p!(n-p)!\sum_{j_1<\ldots<j_p}\lambda_{j_1}\ldots\lambda_{j_p}\ge
n!\,(\lambda_1\ldots\lambda_n)^{p/n}.$$
As both sides are homogeneous of degree $p$ in $(\lambda_j)$, we may
assume $\lambda_1\ldots\lambda_n=1$.
Then our inequality follows from the inequality between the arithmetic and
geometric means of the numbers $\lambda_{j_1}\ldots\lambda_{j_p}$.
Next, we show that statements a) and b) are equivalent in any
dimension $n$.

\noindent a) $\Longrightarrow$ b).  By density, we may suppose that
$u_1\ld u_n$ are K\"ahler classes.  Fix a positive $(n,n)$ form
$f$ such that $\int_X f=1$.  Then Theorem 10.1 implies that there is a
K\"ahler metric $\alpha_j$ representing $u_j$ such that $\alpha_j^n=u_j^n f$.
Inequality a) combined with an integration over $X$ yields
$$u_1\cdots u_n=\int_X\alpha_1\wedge\ldots\wedge\alpha_n\ge
(u_1^n)^{1/n}\ldots(u_n^n)^{1/n}\int_X f.$$
b) $\Longrightarrow$ a). The forms $\alpha_1\ld\alpha_n$
can be considered as constant $(1,1)$-forms on any complex torus
$X=\bbbc^n/\Gamma$. Inequality b) applied to the associated cohomology
classes $u_j\in H^{1,1}(X,\bbbr)$ is then equivalent to a).

Finally we prove a) by induction on $n$, assuming the result already
proved in dimension $n-1$. We may suppose that $\alpha_n$ is positive
definite, say \hbox{$\alpha_n=i\sum dz_j\wedge d\ovl z_j$} in a suitable
basis. Denote by $u_1\ld u_n$ the associated cohomology classes on the
abelian variety $X=\bbbc^n/\bbbz[i]^n$. Then $u_n$ has integral periods,
so some multiple of $u_n$ is the first Chern class of a very ample line
bundle $\cO(D)$ where $D$ is a smooth irreducible divisor in $X$.
Without loss of generality, we may suppose $u_n=c_1(\cO(D))$. Thus
$$u_1\cdots u_{n-1}\cdot u_n=u_{1\restriction D}\cdots u_{n-1\restriction D}$$
and by the induction hypothesis we get
$$u_1\cdots u_n\ge(u_{1\restriction D}^{n-1})^{1/(n-1)}\ldots
(u_{n-1\restriction D}^{n-1})^{1/(n-1)}.$$
However $u_{j\restriction D}^{n-1}=u_j^{n-1}\cdot u_n\ge
(u_j^n)^{(n-1)/n}(u_n^n)^{1/n}$, since a) and b) are equivalent and
a) is already proved in the case of two forms. b) follows in
dimension~$n$, and therefore a) holds in $\bbbc^n$.\qed
\endproof

\begstat{(10.3) Remark} \rm In case $\alpha_j$ (resp.\ $u_j$) are
positive definite, the equality holds in 10.2 (a,b) if and only if
$\alpha_1\ld\alpha_n$ (resp.\ $u_1\ld u_n$) are proportional.  In our
inductive proof, the restriction morphism $H^{1,1}(X,\bbbr)\lra
H^{1,1}(D,\bbbr)$ is injective for $n\ge 3$ by the hard Lefschetz theorem,
hence it is enough to consider the case of $\alpha^p\wedge\beta^{n-p}$.
The equality between arithmetic and geometric means occurs only when all
numbers $\lambda_{j_1}\ldots\lambda_{j_p}$ are equal, so all $\lambda_j$
must be equal and $\alpha=\lambda_1\beta$, as desired.
More generally, one can show (exercise to the reader!) that there is an
inequality
$$\eqalign{
\alpha_1\wedge\ldots\wedge\alpha_p&\wedge\beta_1\wedge\ldots\wedge
\beta_{n-p}\ge\cr
&\ge(\alpha_1^p\wedge\beta_1\wedge\ldots\wedge\beta_{n-p})^{1/p}
\ldots(\alpha_p^p\wedge\beta_1\wedge\ldots\wedge\beta_{n-p})^{1/p}\cr}$$
for all $(1,1)$-forms $\alpha_j$, $\beta_k\ge 0$, and a similar inequality
for products of nef cohomology classes $u_j$, $v_k$.\qed
\endstat

We now show how the Aubin-Calabi-Yau theorem can be applied to construct
singular metrics on ample (or more generally
big and nef) line bundles. We first suppose that $L$ is an ample line
bundle over a projective $n$-fold $X$ and that $L$ is equipped with a
smooth metric of positive curvature. We consider the K\"ahler metric
$\omega={i\over 2\pi}\Theta(L)$. Any form $\wt\omega$ in the K\"ahler class
of $\omega$ can be written as $\wt\omega=\omega+dd^c\psi$,
i.e.\ is the curvature form of $L$ after multiplication of the original
metric by a smooth weight function $e^{-\psi}$. By lemma~5.1, the
Monge-Amp\`ere equation
$$(\omega+dd^c\psi)^n=f\leqno(10.4)$$
can be solved for $\psi$, whenever $f$ is a smooth $(n,n)$-form with
$f>0$ and $\int_X f=L^n$.  In order to produce logarithmic poles at given
points $x_1\ld x_N\in X$, the main idea is to let $f$ converge to a
Dirac measure at $x_j\,$; then $\wt\omega$ will be shown to converge
to a closed positive $(1,1)$-current with non zero Lelong number at $x_j$.

Let $(z_1\ld z_n)$ be local coordinates defined on some neighborhood
$V_j$ of $x_j$, and let
$$\alpha_{j,\varepsilon}=dd^c
\big(\chi(\log|z_j-x_j|/\varepsilon)\big)\leqno(10.5)$$
where $\chi:\bbbr\to\bbbr$ is a smooth convex increasing function such
that $\chi(t)=t$ for $t\ge 0$ and $\chi(t)=-1/2$ for $t\le -1$.
Then $\alpha_{j,\varepsilon}$ is a smooth positive $(1,1)$-form, and
$\alpha_{j,\varepsilon}=dd^c\log|z_j-x_j|$ in the complement of the ball
$|z_j-x_j|\le\varepsilon$. It follows that $\alpha_{j,\varepsilon}^n$
has support in the ball $|z_j-x_j|\le\varepsilon$, and
Stokes' formula gives
$$\int_{B(x_j,\varepsilon)}\alpha_{j,\varepsilon}^n=
\int_{B(x_j,\varepsilon)}\big(dd^c\log|z_j-x_j|\big)^n=1.\leqno(10.6)$$
Hence $\alpha_{j,\varepsilon}^n$ converges weakly to the Dirac measure
$\delta_{x_j}$ as $\varepsilon$ tends to $0$. For all positive
numbers $\tau_j>0$ such that $\sigma:=\sum\tau_j^n<L^n=\int_X\omega^n$,
Theorem 10.1 gives a solution of the Monge-Amp\`ere equation
$$\omega_\varepsilon^n=\sum_{1\le j\le N}\tau_j^n\,\alpha_{j,\varepsilon}^n+
\Big(1-{\sigma\over L^n}\Big)\omega^n~~~~\hbox{\rm with}~~
\omega_\varepsilon=\omega+dd^c\psi_\varepsilon,
\leqno(10.7)$$
since the right-hand side of the first equation is $>0$ and has the
correct integral value $L^n$ over $X$. The solution $\psi_\varepsilon$ is
merely determined up to a constant. If $\gamma$ is an arbitrary
K\"ahler metric on~$X$, we can normalize $\psi_\varepsilon$ in such a
way that  $\int_X\psi_\varepsilon\,\gamma^n=0$.

\begstat{(10.8) Lemma} There is a sequence $\varepsilon_\nu$ converging
to zero such that $\psi_{\varepsilon_\nu}$ has a limit $\psi$ in $L^1(X)$
and such that the sequence of $(1,1)$-forms
$\smash{\omega_{\varepsilon_\nu}}$ converges weakly towards a closed
positive current $T$ of type~$(1,1)$. Moreover, the cohomology class
of $T$ is equal to $c_1(L)$ and $T=\omega+dd^c\psi$.
\endstat

\begproof{} The integral $\int_X \omega_\varepsilon\wedge\gamma^{n-1}=
L\cdot\{\gamma\}^{n-1}$ remains bounded, so we can find a sequence
$\varepsilon_\nu$ converging to zero such that the subsequence
$\omega_{\varepsilon_\nu}$ converges weakly towards a closed positive
current $T$ of bidegree $(1,1)$. The cohomology class of a current is
continuous with respect to the weak topology (this can be seen by
Poincar\'e duality). The cohomology class of $T$ is thus equal to
$c_1(L)$. The function $\psi_\varepsilon$ satisfies the equation
${1\over\pi}\Delta\psi_\varepsilon={\rm tr}_\gamma(\omega_\varepsilon-\omega)$
where $\Delta$ is the Laplace operator associated to~$\gamma$. Our
normalization of $\psi_\varepsilon$ implies
$$\psi_\varepsilon=\pi\,G\,{\rm tr}_\gamma(\omega_\varepsilon-\omega),$$
where $G$ is the Green operator of $\Delta$. As $G$ is
a compact operator from the Banach space of bounded Borel measures into
$L^1(X)$, we infer that some subsequence $(\psi_{\varepsilon_\nu})$ of
our initial subsequence converges to a limit $\psi$ in $L^1(X)$. By the
weak continuity of $dd^c$, we get
$T=\lim(\omega+dd^c\psi_{\varepsilon_\nu})=
\omega+dd^c\psi$.\qed
\endproof

Let $\Omega\subset X$ be an open coordinate patch such that $L$ is
trivial on a neighborhood of $\ovl\Omega$, and let $e^{-h}$ be the weight
representing the initial hermitian metric on $L_{\restriction\ovl\Omega}$.
Then $dd^c h=\omega$ and $dd^c
(h+\psi_\varepsilon)=\omega_\varepsilon$, so the function
$\varphi_\varepsilon=h+\psi_\varepsilon$ defines a psh weight
on $L_{\restriction\Omega}$, as well as its limit $\varphi=h+\psi$.
By the continuity of $G$, we also infer from the proof of Lemma~10.8
that the family $(\psi_\varepsilon)$ is bounded in $L^1(X)$. The usual
properties of subharmonic functions then show that there is a uniform
constant $C$ such that $\varphi_\varepsilon\le C$ on $\ovl\Omega$. We use
this and equation (10.7) to prove that the limit $\varphi$ has logarithmic
poles at all points $x_j\in\Omega$, thanks to Bedford and Taylor's
maximum principle for solutions of Monge-Amp\`ere equations~[BT76]:

\begstat{(10.9) Lemma} Let $u,v$ be smooth $($or continuous$)$ psh
functions on $\ovl\Omega$, where $\Omega$ is a bounded open set in
$\bbbc^n$. If
$$u_{\restriction\partial\Omega}\ge v_{\restriction\partial\Omega}~~~
\hbox{\rm and}~~~(dd^c u)^n\le(dd^c v)^n~~~
\hbox{\rm on}~~\Omega,$$
then $u\ge v$ on $\Omega$.\qed
\endstat

In the application of Lemma 10.9, we suppose that $\Omega$ is a neighborhood
of $x_j$ and take
$$u=\tau_j\big(\chi(\log|z_j-x_j|/\varepsilon)+
\log\varepsilon\big)+C_1,~~~~v=\varphi_\varepsilon,$$
where $C_1$ is a large constant. Then for $\varepsilon>0$ small enough
$$\eqalign{
u_{\restriction\partial\Omega}&=\tau_j\log|z_j-x_j|+C_1,~~~~
v_{\restriction\partial\Omega}\le C,\cr
(dd^cv)^n&=\omega_\varepsilon^n\ge\tau_j^n\,
\alpha_{j,\varepsilon}^n=(dd^cu)^n
{}~~~~\hbox{\rm on}~~\Omega.\cr}$$
For $C_1$ sufficiently large, we infer $u\ge v$ on $\Omega$, hence
$$\varphi_\varepsilon\le\tau_j\log(|z_j-x_j|+\varepsilon)+C_2~~~~
\hbox{\rm on}~~\Omega.$$

\begstat{(10.10) Corollary} The psh weight $\varphi=h+\psi$ on
$\smash{L_{\restriction\Omega}}$ associated to the limit function
$\psi=\lim\smash{\psi_{\varepsilon_\nu}}$ satisfies
$dd^c\varphi=T$. Moreover, $\varphi$
has logarithmic poles at all points $x_j\in\Omega$ and
$$\varphi(z)\le\tau_j\log|z_j-x_j|+O(1)~~~~\hbox{\rm at}~~x_j.\eqno\square$$
\endstat

\begstat{(10.11) Remark} \rm The choice of the coefficients $\tau_j$ is
made according to the order $s_j$ of jets at $x_j$ which sections in
$H^0(X,K_X+L)$ should generate. Corollary 5.13 requires
$\nu(\varphi,x_j)\ge n+s_j$, hence we need only take $\tau_j=n+s_j$.
Accordingly, $L$ must satisfy the numerical condition
$$L^n>\sigma=\sum(n+s_j)^n.$$
In particular, for the question of global generation, we need only take one
point with $s_1=0$, thus $\sigma=n^n$, and for the question of very ampleness
we need either two points with $s_1=s_2=0$ or one point with $s_1=1$, thus
$\sigma=\max(2n^n,(n+1)^n)=(n+1)^n$. In fact $\sigma=2n^n$ is enough to
get very ampleness: in the the case of two infinitely near points $x_1=x_2=0$
in the direction $\partial\over\partial z_n$ (say), we can replace (10.5)
and (10.7) respectively by
$$\leqalignno{
\alpha'_\varepsilon&=dd^c\chi\Big({1\over 2}\log{
|z_1|^2+\ldots+|z_{n-1}|^2+|z_n|^4\over\varepsilon^2}\Big),&(10.5')\cr
\omega_\varepsilon^n&=n^n(\alpha'_\varepsilon)^n+
\Big(1-{2n^n\over L^n}\Big)\omega^n~~~~\hbox{\rm with}~~
\omega_\varepsilon=\omega+dd^c\psi_\varepsilon.&(10.7')\cr}$$
In this case, we have $\int_X(\alpha'_\varepsilon)^n=2$. Arguments similar
to those used in the proof of Corollary 10.10 show that
$$\varphi(z)\le n\log\big(|z_1|+\ldots+|z_{n-1}|+|z_n|^2\big)+O(1),$$
whence $\cI(\varphi)_0\subset(z_1\ld z_{n_1},z_n^2)$ by the result
of Exercise 5.10. Corollary 5.12 can then be used to obtain
separation of infinitesimally near points.\qed
\endstat

\noindent{\bf Case of a big nef line bundle}. All our arguments were
developed under the assumption that $L$ is ample, but if $L$ is only
nef and big, we can proceed in the following way. Let $A$ be a fixed
ample line bundle with smooth curvature form $\gamma={\ii\over 2\pi}
\Theta(A)>0$. As $mL+A$ is ample for any $m\ge 1$, by Theorem~5.1 there
exists a smooth hermitian metric on $L$  depending on $m$, such that
$\omega_m={\ii\over2\pi}\Theta(L)_m+{1\over m}{\ii\over2\pi}\Theta(A)>0$ and
$$\omega_m^n={(L+{1\over m}A)^n\over A^n}\,\gamma^n.\leqno(10.12)$$
However, a priori we cannot control the asymptotic behaviour of $\omega_m$
when $m$ tends to infinity, so we introduce the sequence of non necessarily
positive $(1,1)$-forms $\omega'_m={\ii\over2\pi}\Theta(L)_1+{1\over m}
{\ii\over2\pi}\Theta(A)\in\{\omega_m\}$, which is uniformly bounded in
$C^\infty(X)$ and converges to $\Theta(L)_1$. Then we solve the
Monge-Amp\`ere equation
$$\omega_{m,\varepsilon}^n=\sum_{1\le j\le N}\tau_j^n\,\alpha_{j,\varepsilon}^n
+\Big(1-{\sigma\over(L+{1\over m}A)^n}\Big)\omega_m^n\leqno(10.13)$$
with $\omega_{m,\varepsilon}=\omega'_m+dd^c
\psi_{m,\varepsilon}$ and some smooth function $\psi_{m,\varepsilon}$ such
that $\int_X\psi_{m,\varepsilon}\gamma^n=0\,$; this is again possible
by Yau's theorem~5.1. The numerical condition needed on $\sigma$ to solve
(10.13) is obviously satisfied for all $m$ if we suppose
$$\sigma=\sum\tau_j^n<L^n<\Big(L+{1\over m}A\Big)^n.$$
The same arguments as before show that there exist a
convergent subsequence $\lim_{\nu\to+\infty}\psi_{m_1,\varepsilon_\nu}=
\psi$ in $L^1(X)$ and a closed positive $(1,1)$-current \hbox{$T=\lim
\omega_{m_1,\varepsilon_\nu}$}
$={\ii\over2\pi}\Theta(L)_1+dd^c\psi\in c_1(L)$
such that Corollary~10.10 is still valid; in this case, $h$ is taken to
be the weight function corresponding to $\Theta(L)_1$. Everything thus
works as in the ample case.

\nobreak\titleb{11.}{Numerical Criteria for Very Ample Line Bundles}
In the last section, we explained a construction of psh weights admitting
singularities at all points in a given finite set. Our next goal is to
develope a technique for bounding the Lelong numbers at other points. An
explicit numerical criterion can then be derived from Corollary 5.13.

First suppose that $L$ is an ample line bundle over~$X$. The idea is to
apply the self-intersection inequality~9.5 to the $(1,1)$-current
$T=\lim\omega_{\varepsilon_\nu}$
produced by equation~(10.7), and to integrate the inequality with respect
to the K\"ahler form $\omega={\ii\over2\pi}\Theta(L)$. Before doing this,
we need to estimate the excess of intersection in terms of $T_{\ac}^n$.

\begstat{(11.1) Proposition} The absolutely continuous part $T_{\ac}$
of $T$ satisfies
$$T_{\ac}^n\ge\Big(1-{\sigma\over L^n}\Big)\omega^n~~~~
\hbox{\it a.e.\ on}~~X.$$
\endstat

\begproof{} The result is local, so we can work in an open set
$\Omega$ which is relatively compact in a coordinate patch of $X$.
Let $\rho_\delta$ be a family of
smoothing kernels. By a standard lemma based on the comparison between
arithmetic and geometric means (see e.g.\ [BT76], Proposition~5.1), the
function $A\mapsto(\det A)^{1/n}$ is concave on the cone of nonnegative
hermitian $n\times n$ matrices. Thanks to this concavity property we get
$$\big[\big(\omega_\varepsilon\star\rho_\delta(x)\big)^n\big]^{1/n}\ge
(\omega_\varepsilon^n)^{1/n}\star\rho_\delta(x)\ge
\Big(1-{\sigma\over L^n}\Big)^{1/n}(\omega^n)^{1/n}\star\rho_\delta(x),$$
thanks to equation (6.5). As $\varepsilon_\nu$ tends to $0$,
$\omega_{\varepsilon_\nu}\star\rho_\delta$ converges to
$T\star\rho_\delta$ in the strong topology of $C^\infty(\Omega)$, thus
$$\big((T\star\rho_\delta)^n\big)^{1/n}\ge
\Big(1-{\sigma\over L^n}\Big)^{1/n}(\omega^n)^{1/n}\star\rho_\delta~~~~
\hbox{\rm on}~~\Omega.$$
Now, take the limit as $\delta$ goes to $0$. By the Lebesgue density
theorem $T\star\rho_\delta(x)$ converges almost everywhere
to $T_{\ac}(x)$ on $\Omega$, so we are done.\qed
\endproof

According to the notation used in \S~9, we consider an arbitrary
subset $\Xi\subset X$ and introduce the jumping values
$$b_p=\inf\big\{c>0\,;\,\codim\big(E_c(T),x\big)\ge p,~\forall x\in\Xi\big\}.$$
By Proposition~11.1 and Inequality 10.2~a), we have
$$T_{\ac}^j\wedge\omega^{n-j}\ge
\Big(1-{\sigma\over L^n}\Big)^{j/n}\omega^n.\leqno(11.2)$$
Now, suppose that the ``formal vector bundle'' $T_X\otimes\cO(aL)$ is
nef, i.e.\ that the $\bbbr$-divisor $\cO_{T_X}(1)+a\,\pi^\star L$ is nef for
some constant $a\ge 0$. We can then apply Theorem~9.5 with $u=a\,\omega$ and
$$\{\Theta_p\}=(1+b_1a)\cdots(1+b_pa)\{\omega^p\}\,;$$
by taking the wedge product of $\Theta_p$ with $\omega^{n-p}$, we get
$$\eqalign{(1+b_1a)\ldots(1+b_pa)\int_X\omega^n
&{}\ge\sum_{k\ge 1}(\nu_{p,k}-b_1)\ldots(\nu_{p,k}-b_p)
\int_X[Z_{p,k}]\wedge\omega^{n-p}\cr
&~~{}+\int_X(T_{\ac}+b_1a\,\omega)\wedge\ldots\wedge
(T_{\ac}+b_pa\,\omega)\wedge\omega^{n-p}.\cr}$$
Combining this inequality with (11.2) for $T_{\ac}^{p-j}$ yields
$$\eqalign{
(1+b_1a)\ldots(1+b_pa)\,L^n&\ge
\sum_{k\ge 1}(\nu_{p,k}-b_1)\ldots(\nu_{p,k}-b_p)\,L^{n-p}\cdot Z_{p,k}\cr
&{}+\sum_{0\le j\le p}S_j^p(b)\,a^j\Big(1-{\sigma\over L^n}\Big)^{(p-j)/n}
L^n,\cr}$$
where $S_j^p(b)$, $1\le j\le p$, denotes the elementary symmetric function
of degree $j$ in $b_1\ld b_p$ and $S_0^p(b)=1$. As $\prod(1+b_ja)=\sum
S_j^p(b)a^j$, we get
$$\leqalignno{
{}~~~~\smash{\sum_{k\ge 1}}(\nu_{p,k}-b_1)&{}\ldots(\nu_{p,k}-b_p)\,
L^{n-p}\cdot Z_{p,k}\cr
&\le\sum_{0\le j\le p}S_j^p(b)\,a^j\Big(1-\Big(1-{\sigma\over L^n}
\Big)^{(p-j)/n}\Big)L^n.&(11.3)\cr}$$
If $L$ is only supposed to be big and nef, we follow essentially the same
arguments and replace $\omega$ in all our inequalities by
$\omega_m={\ii\over 2\pi}(\Theta(L)_m+{1\over m}\Theta(A))$ with $A$ ample
(see section~6). Note that all $(n,n)$-forms $\omega_m^n$ were defined to
be proportional to  $\gamma^n=({\ii\over 2\pi}\Theta(A))^n$, so
inequality~11.1 becomes in the limit
$$T_{\ac}^n\ge\Big(1-{\sigma\over L^n}\Big)
{L^n\over A^n}\gamma^n=\Big(1-{\sigma\over L^n}\Big)
{L^n\over(L+{1\over m}A)^n}\omega_m^n.$$
The intersection inequality (11.3) is the expected generalization of
Proposition~8.2 in arbitrary codimension. In this inequality, $\nu_{p,k}$
is the generic Lelong number of $T$ along~$Z_{p,k}$, and $Z_{p,k}$ runs over
all \hbox{$p$-codimensional} components $Y$ of $\bigcup_{c>b_p}E_c(T)$
intersecting $\Xi\,$; by definition of $b_j$ we have $\max_k\nu_{p,k}=b_{p+1}$.
Hence we obtain:

\begstat{(11.4) Proposition} Let $L$ be a big nef line bundle such that
$\cO_{T_X}(1)+a\,\pi^\star L$ is nef, and let $T\in c_1(L)$
be the positive curvature current obtained by concentrating the Monge-Amp\`ere
mass $L^n$ into a finite sum of Dirac measures with total mass~$\sigma$, plus
some smooth positive density spread over $X$ $($equation $(10.7))$.
Then the jumping values $b_p$ of the Lelong number of $T$ over an arbitrary
subset $\Xi\subset X$ satisfy the recursive inequalities
$$(b_{p+1}-b_1)\ldots(b_{p+1}-b_p)\le{1\over\min_Y L^{n-p}\cdot Y}
\sum_{0\le j\le p-1}S_j^p(b)\,a^j\sigma_{p-j}\,,\leqno(11.5)$$
where $\sigma_j=\big(1-(1-\sigma/L^n)^{j/n}\big)L^n$, and where $Y$
runs over all \hbox{$p$-codimensional} subvarieties of $X$
intersecting~$\Xi$.
\endstat

Observe that $\sigma_j$ is increasing in $j\,$; in particular
$\sigma_j<\sigma_n=\sigma$ for $j\le n-1$. Moreover, the convexity of the
exponential function shows that
$$t\mapsto{1\over t}\big(1-(1-\sigma/L^n)^t\big)L^n$$
is decreasing, thus $\sigma_j>\sigma_p j/p$ for $j<p\,$; in particular
$\sigma_j>\sigma j/n$ for $j\le n-1$. We are now in a position to prove
the following general result.

\begstat{(11.6) Main Theorem} Let $X$ be a projective $n$-fold and let
$L$ be a big nef line bundle over~$X$. Fix points $x_j\in X$ and
multiplicities $s_j\in\bbbn$. Set
$$\eqalign{
\sigma_0&=\sum_{1\le j\le N}(n+s_j)^n
\qquad\hbox{resp.\ $\sigma_0=2n^n$ if $N=1$ and $s_1=1$,}
\qquad\hbox{and}\cr
\sigma_p&=\big(1-(1-\sigma_0/L^n)^{p/n}\Big)L^n,\qquad 1\le p\le n-1.\cr}$$
Suppose that the formal vector bundle $T_X\otimes\cO(aL)$ is nef
for some $a\ge0$, that $L^n>\sigma_0$, and that there is a sequence
$0=\beta_1<\ldots<\beta_n\le 1$ with
$$L^{n-p}\cdot Y>(\beta_{p+1}-\beta_1)^{-1}\ldots(\beta_{p+1}-
\beta_p)^{-1}\sum_{0\le j\le p-1}S_j^p(\beta)\,a^j\sigma_{p-j}\leqno(11.7)$$
for every subvariety $Y\subset X$ of codimension $p=1,2\ld n-1$
passing through one of the points $x_j$. Then there is a surjective map
$$H^0(X,K_X+L)\lraww\bigoplus_{1\le j\le N}\cO(K_X+L)_{x_j}\otimes
\big(\cO_{X,x_j}/\gm_{X,x_j}^{s_j+1}\big),$$
i.e.\ $H^0(X,K_X+L)$ generates simultaneously all jets of order $s_j$
at $x_j$.
\endstat

As the notation is rather complicated, it is certainly worth examining
the particular case of surfaces and $3$-folds, for the problem of
getting global generation (taking $\sigma_0=n^n$), resp.\ very ampleness
(taking $\sigma_0=2n^n$). If $X$ is a surface, we
find $\sigma_0=4$ (resp.\ $\sigma_0=8$), and we take $\beta_1=0$, $\beta_2=1$.
This gives only two conditions, namely
$$L^2>\sigma_0,~~~~L\cdot C>\sigma_1\leqno(11.7^{n=2})$$
for every curve $C$ intersecting $\Xi$. These bounds are not
very far from those obtained with Reider's theorem, although they
are not exactly as sharp. If $X$ is a $3$-fold,
we have $\sigma_0=27$ (resp.\ $\sigma_0=64$), and we take
$\beta_1=0<\beta_2=\beta<\beta_3=1$. Therefore our condition
is that there exists $\beta\in{}]0,1[$ such that
$$L^3>\sigma_0,\qquad L^2\cdot S>\beta^{-1}\sigma_1,\qquad
L\cdot C>(1-\beta)^{-1}(\sigma_2+\beta\,a\sigma_1)\leqno(11.7^{n=3})$$
for every curve $C$ or surface $S$ intersecting $\Xi$. If we take
$\beta$ to be of the order of magnitude of $a^{-2/3}$, these bounds show
that the influence of $a$ on the numerical conditions for $L$
is at worst $a^{1/3}$ in terms of the homogeneous quantities $(L^p\cdot
Y)^{1/p}$. In higher dimensions, a careful choice of the $\beta_p$'s
shows that the influence of $a$ on these quantities is always less
than $O(a^{1-\delta_n})$, $\delta_n>0$.

\begproof{of Theorem 11.6.} Select $\tau_j>n+s_j$ so that $L^{n-p}\cdot Y$
still satisfies the above lower bound with the corresponding value
$\sigma=\sum\tau_j^n>\sigma_0$. Then apply Theorem~11.4 with
$\Xi=\{x_1\ld x_N\}$. Inequality (11.5) shows inductively that $b_p<\beta_p$
for $p\ge 2$, so $b_n<1$ and therefore $x_j$ is an isolated point
in $E_1(T)$. On the other hand, the Monge-Amp\`ere equation gives us a
weight function $\varphi$ admitting a logarithmic pole at each point $x_j$,
in such a way that $\nu(\varphi,x_j)\ge\tau_j>n+s_j$ (Corollary 10.10).
However $T=dd^c\varphi$ need not be positive definite. In order to apply
our vanishing theorems, we still have to make the curvature positive
definite everywhere. Since $L$ is nef and big, it has a singular metric
for which the curvature current $T_0=dd^c\varphi_0$ satisfies
$\nu(T_0,x)<1$ everywhere and $T_0\ge\varepsilon\gamma$ for some
K\"ahler metric (Corollary~6.8). Then $T'=(1-\delta)T+\delta T_0$ still
has Lelong numbers $\nu(T',x_j)>n+s_j$ for $\delta>0$ small, and the
inequality $\nu(T',x)\le(1-\delta)\nu(T,x)+\delta$ implies that we have
$\nu(T',x)<1$ for $x$ near~$x_j$, $x\ne x_j$. Corollaries 5.12 and 5.13
imply the Theorem. When $N=1$ and $s_1=1$, we can use Remark 10.11
to show that the value $\sigma_0=2n^n$ is admissible in place of
$(n+1)^n$.\qed
\endproof

\begstat{(11.8) Corollary} Let $X$ be a smooth algebraic surface, and let $L$
be a big nef line bundle over $X$. Then on a given subset $\Xi\subset X$
$$\kern-15pt\table{~~~
&\hrulefill&\hrulefill&\hrulefill&\hrulefill&\hrulefill&\cr
&\verl&~K_X+L&\verl~~\hbox{\it is spanned}~\,
&\verl~\,\hbox{\it separates points}~
&\verl~~\hbox{\it generates $s$-jets}~~&\verl\cr
&\hrulefill&\hrulefill&\hrulefill&\hrulefill&\hrulefill&\cr
&\verl&\verl\qquad\quad L^2>~~~&\verl\quad\quad~\,4&\verl~~~8~~~\verl~~~9~~~
\verl~~\,12&\verl\quad\quad~(2+s)^2&\verl\cr
&\ptstep~~\smash{\raise-3pt\hbox{\it when}}~~
&\hrulefill&\hrulefill&\hrulefill&\hrulefill&\cr
&\verl&\verl~~~\,\forall C,\,~~~L\cdot C>~~
&\verl\quad\quad~\,2&\verl~~~6~~~\verl~~~5~~~\verl~~~4
&\verl\quad~~2+3s+s^2&\verl\cr
&\hrulefill&\hrulefill&\hrulefill&\hrulefill&\hrulefill&\cr}$$
for all curves $C\subset X$ intersecting $\Xi$.
In particular, if $L$ is ample, $K_X+mL$ is always globally spanned
for $m\ge 3$ and very ample for $m\ge 5$.
\endstat

\begproof{} For $s$-jets, we have $\sigma_0=(2+s)^2$, so we find the
condition
$$L^2>(2+s)^2,~~~~L\cdot C>\big(1-(1-(2+s)^2/L^2)^{1/2}\big)L^2.$$
The last constant decreases with $L^2$ and is thus at most equal to
the value obtained when $L^2=(2+s)^2+1\,$; the integral part is then
precisely $2+3s+s^2$.\qed
\endproof

\begstat{(11.9) Remark} \rm Ein and Lazarsfeld [EL94] have recently
obtained an algebraic proof of the Main Theorem, in which they can even
weaken the condition that $T_X\otimes\cO(aL)$ is nef into the condition
that $-K_X+aL$ is nef (by taking determinants, $T_X\otimes\cO(aL)$ nef
$\Rightarrow$ $-K_X+naL$ nef). However, it seems that the lower bounds they
get for $L^{n-p}\cdot Y$ are substantially larger than the ones given
by the analytic method. Also, the assumption on $T_X$ is sufficient to get
universal bounds for $2K_X+L$, as we will see later.\qed
\endstat

For the applications, we introduce a convenient definition of higher jet
generation, following an idea of [BSo93].

\begstat{(11.10) Definition} We say that $L$ generates $s$-jets on a given
subset $\Xi\subset X$ if $H^0(X,L)\lraww \bigoplus J^{s_j}_{x_j}L$ is onto
for any choice of points $x_1\ld x_N\in\Xi$ and integers $s_1\ld s_N$ with
$\sum(s_j+1)=s+1$. We say that $L$ is $s$-jet ample if the above property
holds for $\Xi=X$.
\endstat

With this terminology, $L$ is $0$-jet ample if and only if $L$ is
generated by global sections and $1$-jet ample if and only if $L$
is very ample. In order that $K_X+L$ generates $s$-jets on $\Xi$, the
constant to be used in Theorem~11.6 is $\sigma_0=\max\sum(n+s_j)^n$ over all
decompositions $s+1=\sum(s_j+1)$. In fact, if we set $t_j=s_j+1$, the
following lemma gives $\sigma_0=(n+s)^n$, that is, the maximum is reached
when only one point occurs.

\begstat{(11.11) Lemma} Let $t_1\ld t_N\in{}[1,+\infty[$. Then
$$\sum_{1\le j\le N}(n-1+t_j)^n\le\Big(n-1+\sum_{1\le j\le N}t_j\Big)^n.$$
\endstat

\begproof{} The right hand side is a polynomial with nonnegative
coefficients and the coefficient of a monomial $t_j^k$ involving exactly
one variable is the same as in the left hand side (however, the constant
term is smaller). Thus the difference
is increasing in all variables and we need only consider the case
$t_1=\ldots=t_N=1$. This case follows from the obvious inequality
$$n^nN=n^n+{n\choose 1}n^{n-1}(N-1)\le(n+N-1)^n.\eqno\square$$
\endproof

\begstat{(11.12) Corollary} If all intersection numbers $L^p\cdot Y$
satisfy the inequalities in Theorem 11.6 with $\sigma_0=(n+s)^n$, then
$K_X+L$ is $s$-jet ample.
\endstat

In order to find universal conditions for $K_X+L$ to be very ample, our
main theorem would require a universal value $a$ depending only on
$n=\dim_\bbbc X$ such that $T_X\otimes\cO(aL)$ is always nef. However,
this is clearly impossible as the example of curves already shows: if~$X$
is a curve of genus $g$ and $L$ has degree $1$, then $T_X\otimes\cO(aL)$
is nef if and only if $a\ge 2g-2$. As we will see in \S~13, an explicit
value $a$ depending only on the intersection numbers $L^n$ and
$L^{n-1}\cdot K_X$ exists, but this value is very large. Here, these
difficulties can be avoided by means of the following simple lemma.

\begstat{(11.13) Lemma} Let $F$ be a very ample line bundle over $X$.
Then the vector bundle $T_X\otimes\cO(K_X+nF)$ is nef and generated by
global sections.
\endstat

\begproof{} By the assumption that $F$ is very ample, the $1$-jet bundle
$J^1F$ is generated by its sections. Consider the exact sequence
$$0\lra T^\star_X\otimes F\lra J^1F\lra F\lra 0$$
where $\rank(J^1F)=n+1$ and $\det(J^1F)=K_X+(n+1)F$. The $n$-th exterior
power $\Lambda^n(J^1F)=(J^1F)^\star\otimes\det(J^1F)$ is also generated
by sections. Hence, by dualizing the sequence, we get a surjective morphism
$$\Lambda^n(J^1F)\lra (T^\star_X\otimes F)^\star\otimes
\det(J^1F)=T_X\otimes\cO(K_X+nF).$$
Therefore $T_X\otimes\cO(K_X+nF)$ is generated by sections and, in
particular, it is nef.\qed
\endproof

Before going further, it is convenient to introduce the following
quantitative measurement of ampleness:

\begstat{(11.14) Definition} Let $F$ be a nef line bundle and let
$\Xi\subset X$ be an arbitrary subset. We define
$$\mu(F,\Xi)=\min_{1\le p\le n}\min_{\dim Y=p,\,Y\cap\Xi\ne\emptyset}
(F^p\cdot Y)^{1/p}$$
where $Y$ runs over all $p$-dimensional subvarieties intersecting~$\Xi$.
The main properties of this invariant are\/$:$
\smallskip\noindent
$\bu$ Linearity with respect to $F\,:~~\forall k\ge 0,~~
\mu(kF,\Xi)=k\,\mu(F,\Xi)\,;$
\smallskip\noindent
$\bu$ Nakai-Moishezon criterion\/$:~~F$ is ample if and only if $\mu(F,X)>0$.
\endstat

The next idea consists in the following iteration trick: Lemma~11.13
suggests that a universal lower bound for the nefness of $T_X\otimes\cO(aL')$
can be achieved with $L'=K_X+L$ if $L$ is sufficiently ample. Then it follows
from the Main Theorem 11.6 that $K_X+L'=2K_X+L$ is very ample under suitable
numerical conditions. Lemma~11.13 applied with $F=2K_X+L$ shows that the
bundle \hbox{$T_X\otimes\cO\big((2n+1)K_X+nL\big)$} is nef, and thus
\hbox{$T_X\otimes\cO\big((2n+1)L''\big)$} is nef with
\hbox{$L''=K_X+{1\over 2}L\le L'$}. Hence we see
that the Main Theorem can be iterated. The special value $a=2n+1$ will
play an important role.

\begstat{(11.15) Lemma} Let $L'$ be an ample line bundle over $X$. Suppose
that the vector bundle $T_X\otimes\cO\big((2n+1)L'\big)$ is nef. Then, for
any $s\ge 1$, $K_X+L'$ is $s$-jet ample as soon as $\mu(L',X)\ge 3(n+s)^n$.
$($By Remark 10.11, this is still true with $3(n+s)^n$ replaced by $6n^n$
if $s=1.)$
\endstat

\begproof{} We apply the Main Theorem 11.6 with $\sigma_0=(n+s)^n$ and
$a=2n+1$. Thanks to the inequality $\sigma_j\le\sigma_0$ and to the
identity
$$\sum_{0\le j\le p-1}S_j^p(\beta)=(1+\beta_1a)\ldots(1+\beta_pa),$$
(recall that $\beta_1=0$), we find lower bounds $L^{n-p}\cdot Y>M_p$ with
$$M_p\le{(1+\beta_1a)\ldots(1+\beta_pa)\,\sigma_0\over
(\beta_{p+1}-\beta_1)\ldots(\beta_{p+1}-\beta_p)}.$$
If we choose an increasing sequence $0=\beta_1<\ldots<\beta_n=1$ such that
$\beta_p/\beta_{p+1}$ is increasing, we get $\beta_j/\beta_{p+1}<
\beta_{j+n-p-1}/\beta_n=\beta_{j+n-p-1}$, hence
$$\eqalign{
&\prod_{1\le j\le p}(\beta_{p+1}-\beta_j)=\beta_{p+1}^p\prod_{1\le j\le p}
\Big(1-{\beta_j\over\beta_{p+1}}\Big)\ge\beta_{p+1}^p\prod_{1\le j\le p}
(1-\beta_{j+n-p-1}),\cr
&M_p\le\prod_{1\le j\le n-1}{(1+\beta_ja)\over(1-\beta_j)}\beta_{p+1}^{-p}
\sigma_0.\cr}$$
A suitable choice is $\beta_p=n^{-n(n-p)/(p-1)}$. We then find
$$M_p\le C_n n^{n(n-p-1)}(n+s)^n\le C_n(n+s)^{n(n-p)}$$
where
$$C_n=\prod_{1\le p\le n-1}{\big(1+(2n+1)n^{-n(n-p)/(p-1)}\big)
\over\big(1-n^{-n(n-p)/(p-1)}\big)}.$$
Numerical calculations left to the reader show that that $C_n<3$ for all
$n\ge 2$ and that $C_n=3-4\log n/n+O(1/n)$ as $n\to+\infty$. Lemma 11.15
follows.\qed
\endproof

\begstat{(11.16) Lemma} Let $F$ be a line bundle which generates $s$-jets
at every point. Then $F^p\cdot Y\ge s^p$ for every $p$-dimensional
subvariety $Y\subset X$.
\endstat

\begproof{} Fix an arbitrary point $x\in Y$. Then consider the
singular metric on $F$ given by
$$||\xi||^2={|\xi|^2\over\sum|u_j(z)|^2},$$
where $(u_1\ld u_N)$ is a basis of $H^0(X,F\otimes\gm_x^s)$. By our
assumption, these sections have an isolated common zero of order $s$
at~$x$. Hence $F$ possesses a singular metric such that the weight
$\varphi={1\over 2}\log\sum|u_j|^2$ is psh and has an
isolated logarithmic pole of Lelong number $s$ at~$x$. By the comparison
inequality (3.6) with $\psi(z)=\log|z-x|$, we get
$$F^p\cdot Y\ge\int_{B(x,\varepsilon)}[Y]\wedge
(dd^c\varphi)^p\ge s^p\nu([Y],\psi)=s^p\nu(Y,x)
\ge s^p.\eqno\square$$
\endproof

\begstat{(11.17) Theorem} Let $s\ge 1$ and $m\ge 2$ be arbitrary integers
and let $L$ be an ample line bundle. If $L$ satisfies the numerical condition
$$(m-1)\,\mu(L,X)+s\ge 6(n+s)^n,$$
then $2K_X+mL$ is $s$-jet ample. If $s=1$, the result still holds with
$6(n+s)^n$ replaced by $12n^n$, in particular $2K_X+12n^nL$ is always very
ample.
\endstat

\begproof{} We denote here simply $\mu(L,X)=\mu(L)$, for $\Xi=X$ everywhere
in the proof. As $L$ is ample, there exists an integer
$q$ (possibly very large) such that
$$\cases{
K_X+qL\phantom{\big)}&is ample,\cr
T_X\otimes\cO\big((2n+1)(K_X+qL)\big)&is nef,\cr
\mu(K_X+qL)\ge 3(n+s)^n.&\cr}\leqno(11.18)$$
By lemma 11.15 applied to $L'=K_X+qL$, we find that $F=K_X+L'=
2K_X+qL$ is very ample and generates $s$-jets. In particular
$K_X+{q\over 2}L$ is an ample $\bbbq$-divisor, and for any $p$-dimensional
subvariety $Y\subset X$ we have
$$\eqalign{
\big(K_X+(q-1)L\big)^p\cdot Y
&=\Big({1\over 2}F+(q/2-1)L\Big)^p\cdot Y\cr
&=\sum_{0\le k\le p}{p\choose k}2^{k-p}(q/2-1)^kF^{p-k}\cdot L^k\cdot Y.\cr}$$
By the convexity inequality 10.2~b) and Lemma 11.16 we get
$$F^{p-k}\cdot L^k\cdot Y\ge (F^p\cdot Y)^{1-k/p}(L^p\cdot Y)^{k/p}\ge
s^{p-k}\big(\mu(L)\big)^k.$$
Hence $\big(K_X+(q-1)L\big)^p\cdot Y\ge\big((q/2-1)\,\mu(L)+s/2\big)^p$
and
$$\mu\big(K_X+(q-1)L,X\big)\ge{1\over 2}\big((q-2)\,\mu(L)+s\big).$$
Moreover, Lemma~11.13 applied to $F$ shows that
$$T_X\otimes\cO(K_X+nF)=T_X\otimes\cO\big((2n+1)K_X+nqL\big)$$
is nef. As $nq/(2n+1)\le q/2\le q-1$ for $q\ge 2$, we find that all
properties (11.18) except perhaps the last one remain valid with $q-1$
in place of $q\,$:
$$\cases{
K_X+(q-1)L\phantom{\big)}&is ample,\cr
T_X\otimes\cO\big((2n+1)(K_X+(q-1)L)\big)&is nef,\cr
\mu\big(K_X+(q-1)L,X\big)\ge{1\over 2}\big((q-2)\,\mu(L)+s\big).&\cr}
\leqno(11.19)$$
By induction we conclude that (11.19) is still true for the smallest
integer $q-1=m$ such that
$${1\over 2}\big((q-2)\,\mu(L)+s\big)={1\over2}\big((m-1)\,\mu(L)+s\big)
\ge 3(n+s)^n.$$
For this value of $m$, Lemma~11.15 implies that $2K_X+mL$
generates $s$-jets.\qed
\endproof

\begstat{(11.20) Remark} \rm The condition $(m-1)\,\mu(L,X)+s\ge 6(n+s)^n$
is never satisfied for $m=1$. However, Lemma~8.6 applied with $F=K_X+L$
and $\mu=2$ allows us to obtain also a sufficient condition
in order that $2K_X+L$ generates $s$-jets. It is sufficient that
$2(K_X+L)$ generates jets of order $s'_j=2(n+s_j)+1$ at any of the points
$x_j$ whenever $\sum(s_j+1)=s+1$. For $n\ge 2$ we get
$$\sum(n+s'_j)^n=\sum(2s_j+3n+1)^n\le(3n+3+2s)^n$$
after a short computation. The proof of Theorem~11.17 then yields the
sufficient condition $\mu(L,X)\ge 6(3n+3+2s)^n$.\qed
\endstat

\begstat{(11.21) Remark} \rm If $G$ is a nef line bundle, the Main
Theorem 11.6 is still valid for the line bundle $K_X+L+G$, with the
same lower bounds in the numerical conditions for~$L$; indeed, the
proof rests on the existence of suitable singular hermitian metrics
with positive definite curvature on~$L$, and adding $G$ preserves all
properties of these metrics. It follows that Theorem 11.17 and Remark
11.20 can be applied as well to the line bundle $2K_X+mL+G$, under the
same numerical conditions.\qed
\endstat

\titleb{12.}{Holomorphic Morse Inequalities}
Let $X$ be a compact K\"ahler manifold, $E$ a holomorphic vector bundle of
rank $r$ and $L$ a line bundle over~$X$. If $L$ is equipped with a smooth
metric of curvature form $\Theta(L)$, we define the $q$-index set of $L$
to be the open subset
$$X(q,L)=\left\{x\in X~;~\ii\Theta(L)_x \hbox{~has~}
{\displaystyle q\atop\displaystyle n-q}~
{\hbox{negative eigenvalues}\atop\hbox{positive eigenvalues}}\right\}
\leqno(12.1)$$
for $0\le q\le n$. Hence $X$ admits a partition $X=\Delta\amalg
\coprod_q X(q,L)$ where $\Delta=\{\hbox{$x\in X$}\,;\,\det(\Theta(L)_x)=0\}$
is the degeneracy set. We also introduce
$$X(\le q,L)=\bigcup_{0\le j\le q}X(j,L).\leqno(12.1')$$
It is shown in [Dem85b] that the cohomology groups
$H^q\big(X,E\otimes\cO(kL)\big)$ satisfy the following asymptotic
{\it weak Morse inequalities} as $k\to+\infty$
$$h^q\big(X,E\otimes\cO(kL)\big)\le
r{k^n\over n!}\int_{X(q,L)} (-1)^q\Big({\ii\over2\pi}\Theta(L)\Big)^n+o(k^n).
\leqno(12.2)$$
A sharper form is given by the {\it strong Morse inequalities}
$$\leqalignno{
\sum_{0\le j\le q}(-1)^{q-j}h^j\big(X,E&{}\otimes\cO(kL)\big)\cr
&\le r{k^n\over n!}\int_{X(\le q,L)}(-1)^q\Big({\ii\over2\pi}\Theta(L)\Big)^n
+o(k^n).&(12.2')\cr}$$
These inequalities are a useful complement to the Riemann-Roch formula
when information is needed about individual cohomology groups, and not
just about the Euler-Poincar\'e characteristic.

One difficulty in the application of these inequalities is that the
curvature integral is in general quite uneasy to compute, since it is
neither a topological nor an algebraic invariant. However, the Morse
inequalities can be reformulated in a more algebraic setting in which
only algebraic invariants are involved. We give here two such
reformulations.

\begstat{(12.3) Theorem} Let $L=F-G$ be a holomorphic line bundle over a
compact K\"ahler manifold~$X$, where $F$ and $G$ are numerically effective
line bundles. Then for every $q=0,1,\ldots,n=\dim X$, there is an asymptotic
strong Morse inequality
$$\sum_{0\leq j\leq q}(-1)^{q-j}h^j(X,kL)\leq{k^n\over n!}
\sum_{0\leq j\leq q}(-1)^{q-j}{n\choose j}F^{n-j}\cdot G^j+o(k^n).$$
\endstat

\begproof{} By adding $\varepsilon$ times a K\"ahler metric $\omega$ to the
curvature forms of $F$ and $G$, $\varepsilon>0$ one can write ${\ii\over
2\pi}\Theta(L)=\theta_\varepsilon(F)-\theta_\varepsilon(G)$ where
$\theta_\varepsilon(F)={\ii\over 2\pi}\Theta(F)+\varepsilon\omega$ and
$\theta_\varepsilon(G)={\ii\over 2\pi}\Theta(G)+\varepsilon\omega$ are
positive definite. Let $\lambda_1\geq\ldots\geq\lambda_n>0$ be the
eigenvalues of $\theta_\varepsilon(G)$ with respect to $\theta_\varepsilon
(F)$. Then the eigenvalues of ${\ii\over 2\pi}\Theta(L)$ with respect
to $\theta_\varepsilon(F)$ are the real numbers $1-\lambda_j$ and the
set $X(\le q,L)$ is the set $\{\lambda_{q+1}<1\}$ of points $x\in X$ such
that $\lambda_{q+1}(x)<1$. The strong Morse inequalities yield
$$\sum_{0\leq j\leq q}(-1)^{q-j}h^j(X,kL)\leq{k^n\over n!}
\int_{\{\lambda_{q+1}<1\}}\!\!\!(-1)^q\!\!
\prod_{1\leq j\leq n}(1-\lambda_j)\theta_\varepsilon(F)^n+o(k^n).$$
On the other hand we have
$${n\choose j}\theta_\varepsilon(F)^{n-j}\wedge\theta_\varepsilon(G)^j=
\sigma_n^j(\lambda)\,\theta_\varepsilon(F)^n,$$
where $\sigma_n^j(\lambda)$ is the $j$-th elementary symmetric function in
$\lambda_1,\ldots,\lambda_n\,$, hence
$$\sum_{0\leq j\leq q}(-1)^{q-j}{n\choose j}F^{n-j}\cdot G^j=
\lim_{\varepsilon\to 0}\int_X\sum_{0\leq j\leq q}(-1)^{q-j}
\sigma_n^j(\lambda)\,\theta_\varepsilon(F)^n.$$
Thus, to prove the Lemma, we only have to check that
$$\sum_{0\leq j\leq n}(-1)^{q-j}\sigma_n^j(\lambda)-
\bbbone_{\{\lambda_{q+1}<1\}}(-1)^q\prod_{1\leq j\leq n}(1-\lambda_j)\ge 0$$
for all $\lambda_1\geq\ldots\geq\lambda_n\geq 0$, where $\bbbone_
{\{\ldots\}}$ denotes the characteristic function of a set.
This is easily done by induction on~$n$ (just split apart the parameter
$\lambda_n$ and write $\sigma_n^j(\lambda)=
\sigma_{n-1}^j(\lambda)+\sigma_{n-1}^{j-1}(\lambda)\,\lambda_n$).\qed
\endproof

In the case $q=1$, we get an especially interesting lower bound (this
bound has been observed and used by S.~Trapani [Tra92] in a similar
context).

\begstat{(12.4) Consequence}
$h^0(X,kL)-h^1(X,kL)\geq{k^n\over n!}(F^n-nF^{n-1}\cdot G)-o(k^n)$.\newline
Therefore some multiple $kL$ has a section as soon as $F^n-nF^{n-1}\cdot
G>0$.
\endstat

\begstat{(12.5) Remark} \rm The weaker inequality
$$h^0(X,kL)\geq{k^n\over n!}(F^n-nF^{n-1}\cdot G)-o(k^n)$$
is easy to prove if $X$ is projective algebraic. Indeed, by adding a
small ample $\bbbq$-divisor to $F$ and $G$, we may assume that
$F$, $G$ are ample. Let $m_0G$ be very ample and let $k'$ be the smallest
integer $\geq k/m_0$. Then $h^0(X,kL)\geq h^0(X,kF-k'm_0G)$. We select
$k'$ smooth members $G_j$, $1\leq j\leq k'$ in the linear system $|m_0G|$
and use the exact sequence
$$0\to H^0(X,kF-\sum G_j)\to H^0(X,kF)\to\bigoplus H^0(G_j,kF_{|G_j}).$$
Kodaira's vanishing theorem yields $H^q(X,kF)=0$ and $H^q(G_j,kF_{|G_j})=0$
for $q\geq 1$ and $k\geq k_0$. By the exact sequence combined with
Riemann-Roch, we get
$$\eqalign{
h^0(X,kL)&\geq h^0(X,kF-\sum G_j)\cr
&\geq{k^n\over n!}F^n-O(k^{n-1})-
\sum\Big({k^{n-1}\over(n-1)!}F^{n-1}\cdot G_j-O(k^{n-2})\Big)\cr
&\geq{k^n\over n!}\Big(F^n-n{k'm_0\over k}F^{n-1}\cdot G\Big)-O(k^{n-1})\cr
&\geq{k^n\over n!}\Big(F^n-n\,F^{n-1}\cdot G\Big)-O(k^{n-1}).\cr}$$
(This simple proof is due to F.~Catanese.)\qed
\endstat

\begstat{(12.6) Corollary} Suppose that $F$ and $G$ are nef and that $F$ is
big. Some multiple of $mF-G$ has a section as soon as
$$m>n\,{F^{n-1}\cdot G\over F^n}.$$
\endstat

In the last condition, the factor $n$ is sharp: this is easily seen by taking
$X=\bbbp_1^n$ and $F=\cO(a,\ldots,a)$ and $G=\cO(b_1,\ldots,b_n)$ over
$\bbbp_1^n\,$; the condition of the Corollary is then
$m>\sum b_j/a$, whereas $k(mF-G)$ has a section if and only if
$m\ge\sup b_j/a$; this shows that we cannot replace $n$ by
$n(1-\varepsilon)$.\qed
\medskip

We now discuss another application of Morse inequalities in the case where
$c_1(L)\in N_\psef$. Then the regularization theorem 9.1 allows
us to measure the distance of $L$ to the nef cone $N_\nef$. In that
case, a use of singular metrics combined with 9.1 produces smooth metrics
on $L$ for which an explicit bound of the negative part of the curvature
is known. It follows that (12.2) gives an explicit upper bound of the
cohomology groups of $E\otimes\cO(kL)$ in terms of a polynomial in the first
Chern class~$c_1(L)$ (related techniques have already been used
in [Sug87] in a slightly different context).

\begstat{(12.7) Theorem} Suppose that there is a nef cohomology class
$\{u\}$ in $H^{1,1}(X)$ such that $c_1\big(\cO_{T_X}(1)\big)+\pi^\star\{u\}$
is nef over the hyperplane bundle $P(T^\star_X)$. Suppose
moreover that $L$ is equipped with a singular metric such that
$T={\ii\over 2\pi}\Theta(L)\ge 0$. For $p=1,2\ld n,n+1$ set
$$b_p=\inf\{c>0\,;\,\codim E_c(T)\ge p\},$$
with $b_{n+1}=\max_{x\in X}\nu(T,x)$. Then for any holomorphic
vector bundle $E$ of rank $r$ over~$X$ we have
$$h^q\big(X,E\otimes\cO(kL)\big)\le A_qr\,k^n+o(k^n)$$
where $A_q$ is the cup product
$$A_q={1\over q!\,(n-q)!}\big(b_{n-q+1}\{u\}\big)^q\cdot
\big(c_1(L)+b_{n-q+1}\{u\}\big)^{n-q}$$
in $H^{2n}(X,\bbbr)$, identified with a positive number.
\endstat

\begstat{(12.8) Remark} \rm When $X$ is projective algebraic and $\kappa(L)=n$,
the proof of 6.6~f) shows  that $mL\simeq\cO(A+D)$ with $A$ ample and $D$
effective, for some $m\ge 1$. Then we can choose a singular metric on $L$
such that $T={\ii\over2\pi}\Theta(L)=\omega+m^{-1}[D]$,
where $\omega=m^{-1}{\ii\over2\pi}\Theta(A)$ is a K\"ahler metric. As
$\nu(T,x)=m^{-1}\nu(D,x)$
at each point, the constants $b_j$ of theorem~12.7 are obtained
by counting the multiplicities of the singular points of $D\,$;
for example, if $D$ only has isolated singularities, then $b_1=0$,
$b_2=\ldots=b_n=1/m$. Observe moreover that the nefness
assumption on $\cO_{T_X}(1)$ is satisfied with $\{u\}=c_1(G)$ if $G$ is a nef
$\bbbq$-divisor such that $\cO(T_X)\otimes\cO(G)$ is nef, e.g.\ if
$\cO(S^mT_X)\otimes\cO(mG)$ is spanned by sections for some $m\ge 1$.\qed
\endstat

\begproof{of theorem 12.7.} By definition, we have $0=b_1\le b_2\le\ldots\le
b_n\le b_{n+1}$, and for $c\in{}]b_p,b_{p+1}]$,
$E_c(T)$ has codimension $\ge p$ with some component(s) of codimension~$p$
exactly. Let $\omega$ be a fixed K\"ahler metric on $X$. By adding
$\varepsilon\omega$ to $u$ if necessary, we may assume that $u\ge 0$
and that $\cO_{T_X}(1)$ has a smooth hermitian metric such that
$c\big(\cO_{T_X}(1)\big)+\pi^\star u\ge 0$.

Under this assumption, the approximation theorem 9.1 shows that the
metric of $L$ can be approximated by a sequence of smooth metrics such
that the associated curvature forms $T_j$ satisfy the uniform lower bound
$$T_j\ge -\lambda_j(x)\,u(x)-\varepsilon_j\,\omega(x)\leqno(12.9)$$
where ${\lim\!\downarrow\,}_{j\to+\infty}\varepsilon_j=0$ and
$(\lambda_j)_{j>0}$ is a decreasing sequence of continuous functions on
$X$ such that $\lim_{j\to+\infty}\lambda_j(x)=\nu(T,x)$ at each point.

The estimate (12.2) cannot be used directly with $T={\ii\over2\pi}\Theta(L)$
because wedge products of currents do not make sense in general. Therefore,
we replace ${\ii\over2\pi}\Theta(L)$ by its approximations $T_j$ and try
to find an upper bound for the limit.
\endproof

\begstat{(12.10) Lemma} Let $U_j=X(q,T_j)$ be the $q$-index set
associated to $T_j$ and let $c$ be a positive number. On the open
set $\Omega_{c,j}=\{x\in X\,;\,\lambda_j(x)<c\}$
we have
$$(-1)^q\bbbone_{U_j}T_j^n\le{n!\over q!\,(n-q)!}
\big(c\,u+\varepsilon_j\,\omega\big)^q\wedge\big(T_j+c\,u
+\varepsilon_j\,\omega\big)^{n-q}.$$
\endstat

\begproof{} Write $v=c\,u+\varepsilon_j\,\omega>0$ and let
$\alpha_{1,j}\le\ldots
\le\alpha_{n,j}$ be the eigenvalues of $T_j$ with
respect to $v$ at each point. Then $T_j^n=
\alpha_{1,j}\ldots\alpha_{n,j}\,v^n$ and
$$v^q\wedge(T_j+v)^{n-q}={q!\,(n-q)!\over n!}\!\!
\sum_{1\le i_1<\ldots<i_{n-q}\le n}(1+\alpha_{i_1,j})\ldots
(1+\alpha_{i_{n-q},j})\,v^n.$$
On $\Omega_{c,j}$ we get $T_j\ge -v$
by inequality (12.9), thus $\alpha_{i,j}\ge -1$; moreover, we
have $\alpha_1\le\ldots\le\alpha_q<0$ and $0<\alpha_{q+1}\le\ldots\le
\alpha_n$ on $U_j$. On $\Omega_{c,j}$ we thus find
$$0\le(-1)^q\bbbone_{U_j}\alpha_{1,j}\ldots
\alpha_{n,j}\le
\bbbone_{U_j}\alpha_{q+1,j}\ldots\alpha_{n,j}\le
(1+\alpha_{q+1,j})\ldots(1+\alpha_{n,j}),$$
therefore $(-1)^q\bbbone_{U_j}T_j^n\le
\big(n!/q!\,(n-q)!\big)\,v^q\wedge(T_j+v)^{n-q}$.\qed
\endproof

\begproof{of theorem 12.7 (end).} Set $\Lambda=
\max_X\lambda_1(x)$. By Lemma~12.10 applied with an arbitrary
$c>\Lambda$ we have
$$(-1)^q\bbbone_{U_j}T_j^n\le{n!\over q!(n-q)!}
(\Lambda\,u+\varepsilon_1\omega)^q\wedge(T_j+\Lambda\,u+\varepsilon_1
\omega)^{n-q}~~~~\hbox{\rm on}~~X.$$
Then estimate (12.2) and Lemma 12.10 again imply
$$\leqalignno{
h^q\big(X,E\otimes\cO(kL)\big)
&\le r{k^n\over n!}\int_X(-1)^q\bbbone_{U_j}T_j^n+o(k^n)\cr
&\le{r\,k^n\over q!\,(n-q)!}\Big\{\int_{\Omega_{c,j}}
(c\,u+\varepsilon_j\,\omega)^q\wedge(T_j+c\,u+\varepsilon_j\,\omega)^{n-q}\cr
&{}+\int_{X\ssm\Omega_{c,j}}
(\Lambda\,u+\varepsilon_1\omega)^q\wedge(T_j+\Lambda\,u+\varepsilon_1
\omega)^{n-q}\Big\}+o(k^n).&(12.11)\cr}$$
Since $\lambda_j(x)$ decreases to $\nu(T,x)$ as
$j\to+\infty$, the set $X\ssm\Omega_{c,j}$ decreases
to $E_c(T)$. Now, $T_j+\Lambda\,u+\varepsilon_1\omega$ is a closed
positive $(1,1)$-form belonging to a fixed cohomology class, so the mass of
any wedge power $(T_j+\Lambda\,u+\varepsilon_1\omega)^p$ with respect
to $\omega$ is constant. By weak compactness, there is a subsequence
$(j_\nu)$ such that
$(T_{j_\nu}+\Lambda\,u+\varepsilon_1\omega)^p$ converges
weakly to a closed positive current $\Theta_p$ of bidegree $(p,p)$,
for each $p=1\ld n$.
For $c>b_{p+1}$, we have $\codim E_c(T)\ge p+1$, hence
$\bbbone_{E_c(T)}\Theta_p=0$. It follows that the integral
over $X\ssm\Omega_{c,j}$ in (12.11) converges to $0$ when
$c>b_{n-q+1}$. For the same reason the integral over
$\Omega_{c,j}$ converges to the same limit as
its value over $X$: observe that
$(T_j+c\,u+\varepsilon_j\,\omega)^{n-q}$ can be
expressed in terms of powers of $u,\omega$ and of the positive forms
$(T_j+\Lambda\,u+\varepsilon_1\omega)^p$ with $p\le n-q\,$; thus
the limit is a linear combination with smooth coefficients of the
currents $\Theta_p$, which carry no mass on $E_c(T)$. In the limit,
we obtain
$$h^q\big(X,E\otimes\cO(kL)\big)\le
{r\,k^n\over q!\,(n-q)!}(c\{u\})^q\cdot\big(c_1(L)+c\{u\}\big)^{n-q}+o(k^n),$$
and since this is true for every $c>b_{n-q+1}$, Theorem~12.7
follows.\qed
\endproof

\titleb{13.}{Effective Version of Matsusaka's Big Theorem}
An important problem of algebraic geometry is to find effective bounds
$m_0$ such that multiples $mL$ of an ample line bundle become very
ample for $m\ge m_0$. From a theoretical point of view, this problem
has been solved by Matsusaka [Mat72] and Koll\'ar-Matsusaka [KoM83].
Their result is that there is a bound $m_0=m_0(n,L^n,L^{n-1}\cdot K_X)$
depending only on the dimension and on the first two coefficients $L^n$
and $L^{n-1}\cdot K_X$ in the Hilbert polynomial of~$L$. Unfortunately,
the original proof does not tell much on the actual dependence of $m_0$
in terms of these coefficients. The goal of this section is to find
effective bounds for such an integer~$m_0$, along the lines of [Siu93].
However, one of the technical lemmas used in [Siu93] to deal with
dualizing sheaves can be sharpened. Using this sharpening of the lemma,
Siu's bound will be here substantially improved. We first start
with the simpler problem of obtaining merely a nontrivial section
in $mL$. The idea, more generally, is to obtain a criterion for the
ampleness of $mL-B$ when $B$ is nef. In this way, one is able to subtract
from $mL$ any multiple of $K_X$ which happens to get added by the
application of Nadel's vanishing theorem (for this, replace $B$ by $B$
plus a multiple of $K_X+(n+1)L$).

\begstat{(13.1) Proposition} Let $L$ be an ample line bundle over a
projective $n$-fold $X$ and let $B$ be a nef line bundle over~$X$. Then
$K_X+mL-B$ has a nonzero section for some integer $m$ such that
$$m\le n\,{L^{n-1}\cdot B\over L^n}+n+1.$$
\endstat

\begproof{} Let $m_0$ be the smallest integer ${}>n\,{L^{n-1}\cdot B\over
L^n}$. Then $m_0L-B$ can be equipped with a singular hermitian metric
of positive definite curvature. Let $\varphi$ be the weight of this metric.
By Nadel's vanishing theorem, we have
$$H^q(X,\cO(K_X+mL-B)\otimes\cI(\varphi))=0\qquad\hbox{for $q\ge 1$,}$$
thus $P(m)=h^0(X,\cO(K_X+mL-B)\otimes\cI(\varphi))$ is a polynomial for
$m\ge m_0$. Since $P$ is a polynomial of degree $n$ and is not identically
zero, there must be an integer $m\in[m_0,m_0+n]$ which is not a root.
Hence there is a nontrivial section in
$$H^0(X,\cO(K_X+mL-B))\supset H^0(X,\cO(K_X+mL-B)\otimes\cI(\varphi))$$
for some $m\in[m_0,m_0+n]$, as desired.\qed
\endproof

\begstat{(13.2) Corollary} If $L$ is ample and $B$ is nef, then $mL-B$ has
a nonzero section for some integer
$$m\le n\Big({L^{n-1}\cdot B+L^{n-1}\cdot K_X\over L^n}+n+1\Big).$$
\endstat

\begproof{} By Fujita's result 8.3~a), $K_X+(n+1)L$ is nef. We can thus
replace $B$ by $B+K_X+(n+1)L$ in the result of Prop.~13.1.
Corollary 13.2 follows.\qed
\endproof

\begstat{(13.3) Remark} \rm We do not know if the above Corollary is
sharp, but it is certainly not far from being so. Indeed, for $B=0$, the
initial constant $n$ cannot be replaced by anything smaller than
$n/2\,$: take $X$ to be a product of curves $C_j$ of large genus $g_j$ and
$B=0$; our bound for $L=\cO(a_1[p_1])\otimes\ldots\otimes\cO(a_n[p_n])$
to have $|mL|\neq\emptyset$ becomes $m\le\sum(2g_j-2)/a_j+n(n+1)$, which
fails to be sharp only by a factor $2$ when $a_1=\ldots=a_n=1$ and
$g_1\gg g_2\gg\ldots\gg g_n\to+\infty$. On the other hand, the additive
constant $n+1$ is already best possible when $B=0$ and $X=\bbbp^n$.\qed
\endstat

So far, the method is not really sensitive to singularities (the Morse
inequalities are indeed still true in the singular case as is easily seen
by using desingularizations of the ambient variety). The same is true
with Nadel's vanishing theorem, provided that $K_X$ is replaced by the
$L^2$ dualizing sheaf $\omega_X$ (according to the notation introduced
in Remark 5.17, $\omega_X=K_X(0)$ is the sheaf of holomorphic $n$-forms
$u$ on $X_\reg$ such that $\smash{\ii^{n^2}}u\wedge\ovl u$ is integrable
in a neighborhood of the singular set). Then Prop.~13.1 can be
generalized as
\medskip

\begstat{(13.4) Proposition} Let $L$ be an ample line bundle over a
projective $n$-fold $X$ and let $B$ be a nef line bundle over~$X$. For every
$p$-dimensional $($reduced$)$ algebraic subvariety $Y$ of~$X$,
there is an integer
$$m\le p{L^{p-1}\cdot B\cdot Y\over L^p\cdot Y}+p+1$$
such that the sheaf $\omega_Y\otimes\cO_Y(mL-B)$ has a nonzero section.\qed
\endstat

To proceed further, we need the following useful ``upper estimate''
about $L^2$ dualizing sheaves (this is one of the crucial steps in
Siu's approach; unfortunately, it has the effect of producing rather
large final bounds when the dimension increases).

\begstat{(13.5) Proposition} Let $H$ be a very ample line bundle on a
projective algebraic mani\-fold~$X$, and let $Y\subset X$ be a
$p$-dimensional irreducible algebraic subvariety. If
\hbox{$\delta=H^p\cdot Y$} is the degree of $Y$ with respect to $H$,
the sheaf
$$\cHom\big(\omega_Y,\cO_Y((\delta-p-2)H)\big)$$
has a nontrivial section.
\endstat

Observe that if $Y$ is a smooth hypersurface of degree $\delta$ in
$(X,H)=(\bbbp^{p+1},\cO(1))$, then $\omega_Y=\cO_Y(\delta-p-2)$ and the
estimate is optimal. On the other hand, if $Y$ is a smooth complete
intersection of multidegree $(\delta_1\ld\delta_r)$ in $\bbbp^{p+r}$,
then $\delta=\delta_1\ldots\delta_r$ whilst
$\omega_Y=\cO_Y(\delta_1+\ldots+\delta_r-p-r-1)\,$; in this case,
Prop.~(13.5) is thus very far from being sharp.

\begproof{} Let $X\subset\bbbp^N$ be the embedding given by $H$, so
that $H=\cO_X(1)$. There is a linear projection
$\bbbp^n\dasharrow\bbbp^{p+1}$ whose restriction $\pi:Y\to\bbbp^{p+1}$
to $Y$ is a finite and regular birational map of $Y$ onto an algebraic
hypersurface $Y'$ of degree $\delta$ in $\bbbp^{p+1}$.  Let $s\in
H^0(\bbbp^{p+1},\cO(\delta))$ be the polynomial of degree $\delta$
defining $Y'$. We claim that for any small Stein open set
$W\subset\bbbp^{p+1}$ and any $L^2$ holomorphic $p$-form $u$ on $Y'\cap
W$, there is a $L^2$ holomorphic $(p+1)$-form $\wt u$ on $W$ with
values in $\cO(\delta)$ such that $\wt u_{\restriction Y'\cap
W}=u\wedge ds$. In fact, this is precisely the conclusion of the
Ohsawa-Takegoshi extension theorem [OT87], [Ohs88] (see also [Man93]
for a more general version); one can also invoke more standard local
algebra arguments (see Hartshorne [Har77], Th.~III-7.11). As
$K_{\bbbp^{p+1}}=\cO(-p-2)$, the form $\wt u$ can be seen as a section
of $\cO(\delta-p-2)$ on~$W$, thus the sheaf morphism $u\mapsto u\wedge
ds$ extends into a global section of
$\cHom\big(\omega_{Y'},\cO_{Y'}(\delta-p-2)\big)$. The pull-back by
$\pi^\star$ yields a section of $\cHom\big(\pi^\star\omega_{Y'},
\cO_Y((\delta-p-2)H)\big)$. Since $\pi$ is finite and generically
$1:1$, it is easy to see that $\pi^\star\omega_{Y'}=\omega_Y$. The
Proposition follows.\qed
\endproof

By an appropriate induction process based on the above results, we can
now improve Siu's effective version of the Big Matsusaka Theorem
[Siu93]. Our version depends on a constant $\lambda_n$ such that
$m(K_X+(n+2)L)+G$ is very ample for $m\ge\lambda_n$ and every nef line
bundle~$G$. Corollary~(8.5) shows that $\lambda_n\le{3n+1\choose n}-2n$,
and a similar argument involving the recent results of Angehrn-Siu
[AS94] implies $\lambda_n\le n^3-n^2-n-1$ for $n\ge 2$.
Of course, it is expected that $\lambda_n=1$ in view of the Fujita
conjecture.

\begstat{(13.6) Effective version of the Big Matsusaka Theorem} Let $L$ and
$B$ be nef line bundles on a projective $n$-fold~$X$. Assume that $L$
is ample and let $H$ be the very ample line bundle $H=\lambda_n(K_X+(n+2)L)$.
Then $mL-B$ is very ample for
$$m\ge(2n)^{(3^{n-1}-1)/2}{(L^{n-1}\cdot(B+H))^{(3^{n-1}+1)/2}
(L^{n-1}\cdot H)^{3^{n-2}(n/2-3/4)-1/4}\over(L^n)^{3^{n-2}(n/2-1/4)+1/4}}.$$
In particular $mL$ is very ample for
$$m\ge C_n~(L^n)^{3^{n-2}}\left(n+2+{L^{n-1}\cdot K_X\over L^n}
\right)^{3^{n-2}(n/2+3/4)+1/4}$$
with $C_n=(2n)^{(3^{n-1}-1)/2}(\lambda_n)^{3^{n-2}(n/2+3/4)+1/4}$.
\endstat

\begproof{} We use Prop.~(13.4) and Prop.~(13.5) to construct inductively a
sequence of (non necessarily irreducible) algebraic subvarieties
$X=Y_n\supset Y_{n-1}\supset\ldots\supset Y_2\supset Y_1$ such that
$Y_p=\bigcup_jY_{p,j}$ is $p$-dimensional, and $Y_{p-1}$ is obtained
for each $p\ge 2$ as the union of zero sets of sections
$$\sigma_{p,j}\in H^0(Y_{p,j},\cO_{Y_{p,j}}(m_{p,j}L-B))$$
with suitable integers $m_{p,j}\ge 1$. We proceed by induction on
decreasing values of the dimension~$p$, and find inductively upper
bounds $m_p$ for the integers $m_{p,j}$.

By Cor.~(13.2), an integer $m_n$ for $m_nL-B$ to have a section
$\sigma_n$ can be found with
$$m_n\le n\,{L^{n-1}\cdot(B+K_X+(n+1)L)\over L^n}
\le n\,{L^{n-1}\cdot(B+H)\over L^n}.$$
Now suppose that the sections $\sigma_n$, $\ldots$, $\sigma_{p+1,j}$ have
been constructed. Then we get inductively a $p$-cycle $\wt
Y_p=\sum\mu_{p,j}Y_{p,j}$ defined by $\wt Y_p={}$ sum of zero divisors
of sections $\sigma_{p+1,j}$ in $\wt Y_{p+1,j}$, where the mutiplicity
$\mu_{p,j}$ on $Y_{p,j}\subset Y_{p+1,k}$ is obtained by multiplying
the corresponding multiplicity $\mu_{p+1,k}$ with the vanishing order
of $\sigma_{p+1,k}$ along $Y_{p,j}$. As cohomology classes, we find
$$\wt Y_p\equiv\sum(m_{p+1,k}L-B)\cdot (\mu_{p+1,k}Y_{p+1,k})\le
m_{p+1} L\cdot \wt Y_{p+1}.$$
Inductively, we thus have the numerical inequality
$$\wt Y_p\le m_{p+1}\ldots m_n L^{n-p}.$$
Now, for each component $Y_{p,j}$, Prop.~(13.4) shows that there exists
a section of $\omega_{Y_{p,j}}\otimes\cO_{Y_{p,j}}(m_{p,j}L-B)$ for some
integer
$$m_{p,j}\le p{L^{p-1}\cdot B\cdot Y_{p,j}\over L^p\cdot Y_{p,j}}+p+1
\le pm_{p+1}\ldots m_n\,L^{n-1}\cdot B+p+1.$$
Here, we have used the obvious lower bound $L^{p-1}\cdot Y_{p,j}\ge 1$
(this is of course a rather weak point in the argument). The degree of
$Y_{p,j}$ with respect to $H$ admits the upper bound
$$\delta_{p,j}:=H^p\cdot Y_{p,j}\le
m_{p+1}\ldots m_n H^p\cdot L^{n-p}.$$
We use the Hovanski-Teissier concavity inequality (10.2$\,$b)
$$(L^{n-p}\cdot H^p)^{1\over p}(L^n)^{1-{1\over p}}\leq L^{n-1}\cdot H$$
to express our bounds in terms of the intersection numbers $L^n$ and
$L^{n-1}\cdot H$ only. We then get
$$\delta_{p,j}\le m_{p+1}\ldots m_n{(L^{n-1}\cdot H)^p\over (L^n)^{p-1}}.$$
By Prop.~(13.5), there is a nontrivial section in
$$\cHom\big(\omega_{Y_{p,j}},\cO_{Y_{p,j}}((\delta_{p,j}-p-2)H)\big).$$
Combining this section with the section in $\omega_{Y_{p,j}}\otimes
\cO_{Y_{p,j}}(m_{p,j}L-B)$ already cons\-truc\-ted, we get a section of
$\cO_{Y_{p,j}}(m_{p,j}L-B+(\delta_{p,j}-p-2)H)$ on~$Y_{p,j}$. Since we do
not want $H$ to appear at this point, we replace $B$ with
$B+(\delta_{p,j}-p-2)H$ and thus get a section $\sigma_{p,j}$ of
$\cO_{Y_{p,j}}(m_{p,j}L-B)$ with some integer $m_{p,j}$ such that
$$\eqalign{m_{p,j}
&\le pm_{p+1}\ldots m_n\,L^{n-1}\cdot(B+(\delta_{p,j}-p-2)H)+p+1\cr
&\le p\,m_{p+1}\ldots m_n\,\delta_{p,j}\,L^{n-1}\cdot(B+H)\cr
&\le p\,(m_{p+1}\ldots m_n)^2{(L^{n-1}\cdot H)^p\over(L^n)^{p-1}}\,
L^{n-1}\cdot(B+H).\cr}$$
Therefore, by putting $M=n\,L^{n-1}\cdot(B+H)$, we get the recursion relation
$$m_p\le M\,{(L^{n-1}\cdot H)^p\over (L^n)^{p-1}}\,(m_{p+1}\ldots m_n)^2\qquad
\hbox{for $2\le p\le n-1$},$$
with initial value $m_n\le M/L^n$. If we let $(\ovl m_p)$ be the sequence
obtained by the same recursion formula with equalities instead of
inequalities, we get $m_p\le\ovl m_p$ with $\ovl m_{n-1}=
M^3(L^{n-1}\cdot H)^{n-1}/(L^n)^n$ and
$$\ovl m_p={L^n\over L^{n-1}\cdot H}\,\ovl m_{p+1}^2\ovl m_{p+1}$$
for $2\le p\le n-2$. We then find inductively
$$m_p\le\ovl m_p=M^{3^{n-p}}{(L^{n-1}\cdot H)^{3^{n-p-1}(n-3/2)+1/2}\over
(L^n)^{3^{n-p-1}(n-1/2)+1/2}}.$$
We next show that $m_0L-B$ is nef for
$$m_0=\max\big(m_2\,,\,m_3\ld m_n\,,\,m_2\ldots m_n\,L^{n-1}\cdot B\big).$$
In fact, let $C\subset X$ be an arbitrary irreducible curve. Either
$C=Y_{1,j}$ for some $j$ or there exists an integer $p=2\ld n$ such that
$C$ is contained in $Y_p$ but not in~$Y_{p-1}$. If $C\subset
Y_{p,j}\ssm Y_{p-1}$, then $\sigma_{p,j}$ does not vanish identically
on~$C$. Hence $(m_{p,j}L-B)_{\restriction C}$ has nonnegative degree and
$$(m_0L-B)\cdot C\ge(m_{p,j}L-B)\cdot C\geq 0.$$
On the other hand, if $C=Y_{1,j}\,$, then
$$(m_0L-B)\cdot C\ge m_0-B\cdot\wt Y_1\ge m_0-m_2\ldots m_n\,L^{n-1}\cdot
B\ge 0.$$
By the definition of $\lambda_n$ (and by Cor.~(8.5) showing that such a
constant exists), $H+G$ is very ample for every nef line bundle~$G$, in
particular $H+m_0L-B$ is very ample. We thus replace again $B$ with $B+H$.
This has the effect of replacing $M$ with $M=n\,\big(L^{n-1}\cdot(B+2H)\big)$
and $m_0$ with
$$m_0=\max\big(m_n\,,\,m_{n-1}\ld m_2\,,\,m_2\ldots m_n\,L^{n-1}\cdot(B+H)
\big).$$
The last term is the largest one, and from the estimate on $\ovl m_p\,$, we
get
$$\eqalignno{m_0
&\le M^{(3^{n-1}-1)/2}{(L^{n-1}\cdot H)^{(3^{n-2}-1)(n-3/2)/2+(n-2)/2}
(L^{n-1}\cdot(B+H))\over(L^n)^{(3^{n-2}-1)(n-1/2)/2+(n-2)/2+1}}\cr
&\le(2n)^{(3^{n-1}-1)/2}{(L^{n-1}\cdot(B+H))^{(3^{n-1}+1)/2}
(L^{n-1}\cdot H)^{3^{n-2}(n/2-3/4)-1/4}\over(L^n)^{3^{n-2}(n/2-1/4)+1/4}}
\cr}$$
\qed
\endproof

\begstat{(13.7) Remark} \rm In the surface case $n=2$, one can take
$\lambda_n=1$ and our bound yields $mL$ very ample for
$$m\ge 4\,{(L\cdot(K_X+4L))^2\over L^2}.$$
If one looks more carefully at the proof, the initial constant $4$ can be
replaced by~$2$. In fact, it has been shown recently by Fern\'andez del
Busto that $mL$ is very ample for
$$m>{1\over 2}\left[{(L\cdot(K_X+4L)+1)^2\over L^2}+3\right],$$
and an example of G.~Xiao shows that this bound is essentially optimal
(see [FdB94]).
\endstat

\titlec{}{\phantom{$~$}\newline\tbfontt References}
\eightpoint
\letter KoMM87|

\article AN54|Y.\ Akizuki {\rm and} S.\ Nakano|Note on Kodaira-Spencer's
proof of Lefschetz theorems|Proc.\ Jap.\ Acad.|30|1954|266-272|

\miscell AS94|U.\ Angehrn {\rm and} Y.T.\ Siu|Effective freeness and point
separation for adjoint bundles|Preprint December {\oldstyle 1994}|

\article Aub78|Aubin, T|Equations du type Monge-Amp\`ere sur les vari\'et\'es
k\"ahl\'eriennes compactes|C.R.\ Acad.\ Sci.\ Paris Ser.\ A {\bf
283}$\,$({\oldstyle 1976}), 119-121$\,$; Bull.\ Sci.\ Math|102|1978|63-95|

\article AV65|A.\ Andreotti {\rm and} E.\ Vesentini|Carleman estimates
for~the Laplace-Beltrami equation in complex manifolds|Publ.\ Math.\
I.H.E.S.|25|1965|81-130|

\article BBS89|M.\ Beltrametti, A.\ Biancofiore {\rm and} A.J.\
Sommese|Pro\-jective $n$-folds of log-general type, I|Trans.\ Amer.\
Math.\ Soc.|314|1989|825-849|

\miscell BeS93|M.\ Beltrametti {\rm and} A.J.\ Sommese|On $k$-jet
ampleness|Complex Analysis and Geometry, Univ.\ Series in
Math., edited by V.~Ancona and A.~Silva, Plenum Press, New-York,
{\oldstyle 1993}, 355-376|

\article BFS89|M.\ Beltrametti, P.\ Francia {\rm and} A.J.\ Sommese|On
Reider's method and higher order embeddings|Duke Math.\ J.|58|1989|425-439|

\article Boc48|S.\ Bochner|Curvature and Betti numbers (I) and (II)|Ann.\
of Math.|49|1948|379-390$\,$; {\bf 50} ({\oldstyle 1949}), 77-93|

\article Bom70|E.\ Bombieri|Algebraic values of meromorphic maps|Invent.\
Math.|10|1970|267-287 and {\it Addendum}, Invent.\ Math.\ {\bf 11}
({\oldstyle 1970}), 163-166|

\article Bom73|E.\ Bombieri|Canonical models of surfaces of general
type|Publ.\ Math.\ IHES|42|1973|171-219|

\miscell Bon93|L.\ Bonavero|In\'egalit\'es de Morse holomorphes
singuli\`eres|Pr\'epublication Universit\'e de Grenoble~I, Institut
Fourier, n${}^\circ\,259$, {\oldstyle 1993}|

\article BSk74|J.\ Brian\c con {\rm et} H.\ Skoda|Sur la cl\^oture
int\'egrale d'un id\'eal de germes de fonctions holomorphes en un point
de $\bbbc^n$|C.\ R.\ Acad.\ Sc.\ Paris, s\'er.\ A|278|1974|949-951|

\article BT76|E.\ Bedford {\rm and} B.A.\ Taylor|The Dirichlet problem for the
complex Monge-Amp\`ere equation|Invent.\ Math.|37|1976|1-44|

\article BT82|E.\ Bedford {\rm and} B.A.\ Taylor|A new capacity for
plurisubharmonic functions|Acta Math.|149|1982|1-41|

\article Cam92|F.\ Campana|Connexit\'e rationnelle des vari\'et\'es de
Fano|Ann.\ Sci.\ Ecole Norm.\ Sup.|25|1992|539-545|

\miscell Cat88|F.\ Catanese|Footnotes to a theorem of I.\ Reider|Proc.\
Intern.\ Conf.\ on Algebraic Geometry (L'Aquila, June {\oldstyle 1988}),
Lecture Notes in Math., Vol.~1417, Springer-Verlag, Berlin,
{\oldstyle 1990}, 67-74|

\bookref CLN69|S.S.\ Chern, H.I.\ Levine {\rm and} L.\ Nirenberg|Intrinsic
norms on a complex manifold|Global Analysis (papers in honor of K.Kodaira),
p.~119-139, Univ.\ of Tokyo Press, Tokyo|1969|

\article DeI87|P.\ Deligne {\rm and} L.\ Illusie|Rel\`evements modulo
$p^2$ et d\'e\-com\-po\-sition du complexe de De Rham|Invent.\
Math.|89|1987|247-270|

\article Dem82a|J.-P.\ Demailly|Sur les nombres de Lelong associ\'es \`a
l'image directe d'un courant positif ferm\'e|Ann.\ Inst.\ Fourier
(Grenoble)|32|1982|37-66|

\article Dem82b|J.-P.\ Demailly|Estimations $L^2$ pour l'op\'erateur $\ovlp$
d'un fibr\'e vectoriel holomorphe semi-positif au dessus d'une vari\'et\'e
k\"ahl\'e\-rienne compl\`ete|Ann.\ Sci.\ Ec.\ Norm.\ Sup.|15|1982|457-511|

\article Dem85a|J.-P.\ Demailly|Mesures de Monge-Amp\`ere et caract\'erisation
g\'eom\'etrique des vari\'et\'es alg\'ebriques affines|M\'em.\ Soc.\ Math.\
France (N.S.)|19|1985|1-124|

\article Dem85b|J.-P.\ Demailly|Champs magn\'etiques et in\'egalit\'es de
Morse pour la $d''$-coho\-mo\-logie|Ann.\ Inst.\ Fourier
(Grenoble)|35|1985|189-229|

\article Dem87|J.-P.\ Demailly|Nombres de Lelong g\'en\'eralis\'es,
th\'eor\`emes d'int\'egralit\'e et
d'analyticit\'e|Acta Math.|159|1987|153-169|

\article Dem89|J.-P.\ Demailly|Transcendental proof of a generalized
Ka\-wa\-mata-Viehweg vanishing theorem|C.\ R.\ Acad.\ Sci.\ Paris S\'er.\
I Math.|309|1989|123-126~~~and:\hfill\break
Proceedings of the Conference ``Geometrical and algebraical aspects
in several complex variables" held at Cetraro, Univ.\ della Calabria,
June {\oldstyle 1989}|

\miscell Dem90|J.-P.\ Demailly|Singular hermitian metrics on positive
line bundles|Proc.\ Conf.\ Complex algebraic varieties (Bayreuth,
April~2-6, 1990), edited by K.~Hulek, T.~Peternell, M.~Schneider,
F.~Schreyer, Lecture Notes in Math., Vol.~1507, Springer-Verlag, Berlin,
{\oldstyle 1992}|

\article Dem92|J.-P.\ Demailly|Regularization of closed positive currents
and Intersection Theory|J.\ Alg.\ Geom.|1|1992|361-409|

\miscell Dem93a|J.-P.\ Demailly|Monge-Amp\`ere operators, Lelong numbers
and intersection theory|Complex Analysis and Geometry, Univ.\ Series in
Math., edited by V.~Ancona and A.~Silva, Plenum Press, New-York,
{\oldstyle 1993}, 115-193|

\article Dem93b|J.-P.\ Demailly|A numerical criterion for very ample line
bundles|J.\ Differential Geom.|37|1993|323-374|

\miscell Dem95|J.-P.\ Demailly|Effective bounds for very ample line
bundles|Preprint Institut Fourier (January {\oldstyle 1995})|

\article DPS94|J.-P.\ Demailly, Th.\ Peternell, M.\ Schneider|Compact
complex manifolds with numerically effective tangent bundles|J.\
Algebraic Geometry|3|1994|295-345|

\miscell Ein94|L.\ Ein|Note on higher dimensional adjoint linear
systems|preprint {\oldstyle 1994}, personal communication to the
author|

\miscell EKL94|L.\ Ein, O.\ K\"uchle {\rm and} R.\ Lazarsfeld|Local
positivity of ample line bundles|preprint August {\oldstyle 1994},
to appear|

\miscell EL92|L.\ Ein and R.\ Lazarsfeld|Seshadri constants on smooth
surfaces|Journ\'ees de G\'eom\'etrie Alg\'ebrique d'Orsay, July
{\oldstyle 1992}, Ast\'erisque {\bf 282}, {\oldstyle 1993}, 177-186|

\article EL93|L.\ Ein {\rm and} R.\ Lazarsfeld|Global generation of
pluricanonical and adjoint linear series on smooth projective
threefolds|Jour.\ of Am.\ Math.\ Soc.|6|1993|875-903|

\article EM84|H.\ El Mir|Sur le prolongement des courants positifs
ferm\'es|Acta Math.|153|1984|1-45|

\miscell Eno93|I.\ Enoki|Kawamata-Viehweg vanishing theorem for compact
K\"ahler manifolds|Einstein metrics and Yang-Mills connections
(ed.\ T.~Mabuchi, S.~Mukai), Marcel Dekker, {\oldstyle 1993}, 59-68|

\article EV86|H.\ Esnault {\rm and} E.\ Viehweg|Logarithmic De Rham
complexes and vanishing theorems|Invent.\ Math.|86|1986|161-194|

\miscell EV92|H.\ Esnault {\rm and} E.\ Viehweg|Lectures on vanishing
theorems|DMV Seminar, Band {\bf 20}, Birkh\"auser Verlag, {\oldstyle 1992}|

\bookref Fed69|H.\ Federer|Geometric measure theory|Springer Verlag,
Vol. 153, Berlin, Heidelberg, New-York|1969|

\miscell FdB93|G.\ Fern\'andez del Busto|Bogomolov instability and
Kawamata-Viehweg vanishing|preprint UCLA, {\oldstyle 1993}, to appear|

\miscell FdB94|G.\ Fern\'andez del Busto|UCLA Thesis|{\oldstyle 1994}|

\article Fuj83|T.\ Fujita|Semipositive line bundles|J.\ Fac.\ Sci.\ Univ.\
of Tokyo|30|1983|353-378|

\miscell Fuj87|T.\ Fujita|On polarized manifolds whose adjoint bundles
are not semipositive|Algebraic Geometry, Sendai, {\oldstyle 1985},
Adv.\ Stud.\ in Pure Math., Vol.~10, North Holland, T.~Oda (ed.),
{\oldstyle 1987}, 167-178|

\miscell Fuj88|T.\ Fujita|Problem list|Conference held at the
Taniguchi Foundation, Katata, Japan, August {\oldstyle 1988}|

\miscell Fuj93|T.\ Fujita|Remarks on Ein-Lazarsfeld criterion of
spannedness of adjoint bundles of polarized threefolds|preprint
{\oldstyle 1993}, to appear|

\article Gir76|J.\ Girbau|Sur le th\'eor\`eme de Le Potier d'annulation de
la cohomologie|C.~R.\ Acad.\ Sci.\ Paris, s\'erie A|283|1976|355-358|

\miscell Gri69|P.A.\ Griffiths|Hermitian differential geometry, Chern
clas\-ses and positive vector bundles|Global Analysis, papers in honor of
K.~Kodaira, Princeton Univ.\ Press, Princeton, {\oldstyle 1969}, 181-251|

\bookref GH78|P.A.\ Griffiths, J.\ Harris|Principles of algebraic
geometry|Wiley, New York|1978|

\miscell Har70|R.\ Hartshorne|Ample subvarieties of algebraic
varieties|\rm Lecture Notes in Math., Vol.~156, Springer-Verlag,
Berlin ({\oldstyle 1970})|

\article Hir64|H.\ Hironaka|Resolution of singularities of an algebraic
variety over a field of characteristic zero|Ann.\ of Math.|79|1964|109-326|

\article H\"or65|L.\ H\"ormander|$L^2$ estimates and existence theorems for
the $\ovl\partial$ operator|Acta Math.|113|1965|89-152|

\bookref H\"or66|L.\ H\"ormander|An introduction to Complex Analysis in
several variables|{\oldstyle 1966}, 3rd edition, North-Holland Math.\
Libr., vol.7, Amsterdam, London|1990|

\article Hov79|A.G.\ Hovanski|Geometry of convex bodies and algebraic
geometry|Uspehi Mat.\ Nauk|34 {\rm(4)} |1979|160-161|

\article Kaw82|Y.\ Kawamata|A generalization of Kodaira-Ramanujam's
va\-nishing theorem|Math.\ Ann.|261|1982|43-46|

\article Kaw84|Y.\ Kawamata|The cone of curves of algebraic varieties|Ann.\
of Math.|119|1984|603-633|

\article Kis78|C.O.\ Kiselman|The partial Legendre transformation for
plu\-ri\-sub\-harmonic functions|Invent.\ Math.|49|1978|137-148|

\article Kis79|C.O.\ Kiselman|Densit\'e des fonctions
plurisousharmoniques|Bull.\ Soc.\ Math.\ France|107|1979|295-304|

\article Kis84|C.O.\ Kiselman|Sur la d\'efinition de l'op\'erateur de
Monge-Amp\`ere complexe|Analyse Complexe, Proceedings of the Jour\-n\'ees
Fermat (SMF), Toulouse, {\oldstyle 1983}, Lecture Notes in Math., Vol. 1094,
Springer-Verlag|\rm Berlin|1984|139-150|

\miscell KMM87|Y.\ Kawamata, K.\ Matsuda {\rm and} K.\ Matsuki|Introduction
to the minimal problem|Algebraic Geometry, Sendai, {\oldstyle 1985},
Adv.\ Stud.\ in Pure Math., Vol.~10, T.~Oda (ed.), North Holland, Amsterdam,
{\oldstyle 1987}, 283-360|

\article KobO73|S.\ Kobayashi {\rm and} T.\ Ochiai|Characterisations of
complex projective spaces and hyperquadrics|J.\ Math.\ Kyoto
Univ.|13|1973|31-47|

\article Kod53|K.\ Kodaira|On a differential geometric method in the theory
of analytic stacks|Proc.\ Nat.\ Acad.\ Sci.\ USA|39|1953|1268-1273|

\article Kod54|K.\ Kodaira|On K\"ahler varieties of restricted type|Ann.\
of Math.|60|1954|28-48|

\miscell Kol85|J.\ Koll\'ar|Vanishing theorems for cohomology groups|
Algebraic Geometry, Bowdoin {\oldstyle 1985}, Proc.\ Symp.\ Pure Math.,
Vol.~46, {\oldstyle 1987}, 233-243|

\article Kol92|J.\ Koll\'ar|Effective basepoint freeness|Math.\
Ann.|296|1993|595-605|

\article KoM83|J.\ Koll\'ar {\rm and} T.\ Matsusaka|Riemann-Roch type
inequalities|Amer.\ J.\ of Math.|105|1983|229-252|

\article KoMM92|J.\ Koll\'ar, Y.\ Miyaoka {\rm and} S.\ Mori|Rational
connectedness and boundedness of Fano manifolds|J.\ Differential
Geom.|36|1992|765-779|

\miscell Laz93|R.\ Lazarsfeld, {\rm with the assistance of} G.\ Fern\'andez
del Busto|Lec\-tu\-res on linear series|Park City, IAS Mathematics Series,
Vol.~3, {\oldstyle 1993}|

\article Lel57|P.\ Lelong|Int\'egration sur un ensemble analytique
complexe|Bull.\ Soc.\ Math.\ France|85|1957|239-262|

\miscell Lel67|P.\ Lelong|Fonctionnelles analytiques et fonctions enti\`eres
($n$ variables)|S\'em.\ de Math.\ Sup\'erieures, 6e session, \'et\'e
{\oldstyle 1967}, Presses Univ.\ Montr\'eal {\oldstyle 1968}|

\bookref Lel69|P.\ Lelong|Plurisubharmonic functions and positive differential
forms|Gordon and Breach, New-York, and Dunod, Paris|1969|

\article LP75|J.\ Le Potier|Annulation de la cohomologie \`a valeurs dans un
fibr\'e vectoriel holomorphe de rang quelconque|Math.\ Ann.|218|1975|35-53|

\article Man93|L.\ Manivel|Un th\'eor\`eme de prolongement $L^2$ de
sections holomorphes d'un fibr\'e vectoriel|Math.\
Zeitschrift|212|1993|107-122|

\article Mat72|T.\ Matsusaka|Polarized varieties with a given Hilbert
polynomial|Amer.\ J.\ of Math.|94|1972|1027-1077|

\article Mor82|S.\ Mori|Threefolds whose canonical bundles are not
numerically effective|Ann.\ of Math.|116|1982|133-176|

\article Nad89|A.M.\ Nadel|Multiplier ideal sheaves and
K\"ahler-Einstein metrics of positive scalar curvature|Proc.\ Nat.\
Acad.\ Sci.\ U.S.A.|86|1989|7299-7300~~~and~~
Annals of Math., {\bf 132} ({\oldstyle 1990}), 549-596|

\article Nak55|S.\ Nakano|On complex analytic vector bundles|J.~Math.\
Soc.\ Japan|7|1955|1-12|

\article Ohs88|T.\ Ohsawa|On the extension of $L^2$ holomorphic
functions, II|Publ.\ RIMS, Kyoto Univ.|24|1988|265-275|

\article OT87|T.\ Ohsawa {\rm and} K.\ Takegoshi|On the extension of $L^2$
holomorphic functions|Math.\ Zeitschrift|195|1987|197-204|

\article Rei88|I.\ Reider|Vector bundles of rank 2 and linear systems on
algebraic surfaces|Ann.\ of Math.|127|1988|309-316|

\miscell Sak88|F.\ Sakai|Reider-Serrano's method on normal surfaces|Proc.\
Intern.\ Conf.\ on Algebraic Geometry (L'Aquila, June {\oldstyle 1988}),
Lecture Notes in Math., Vol.~1417, Springer-Verlag, Berlin, {\oldstyle 1990},
301-319|

\article Sch74|M.\ Schneider|Ein einfacher Beweis des Verschwindungs\-satzes
f\"ur positive holomorphe Vektorraumb\"undel|Manuscripta Math.|11|1974|95-101|

\miscell Ser54|J.-P.\ Serre|Fonctions automorphes: quelques majorations
dans le cas o\`u $X/G$ est compact|S\'em.\ Cartan (1953-54) 2-1 \`a 2-9|

\article Ser55|J.-P.\ Serre|Un th\'eor\`eme de dualit\'e|Comment.\
Math.|29|1955|9-26|

\article Ser56|J.-P.\ Serre|G\'eom\'etrie alg\'ebrique et g\'eom\'etrie
analytique|Ann.\ Inst.\ Fourier (Grenoble)|6|1956|1-42|

\bookref ShSo85|B.\ Shiffman, A.J.\ Sommese|Vanishing theorems on
complex manifolds|Progress in Math.\ no~{\bf 56}, Birkh\"auser|1985|

\article Sib85|N.\ Sibony|Quelques probl\`emes de prolongement de courants
en analyse complexe|Duke Math.\ J.|52|1985|157-197|

\article Sie55|C.L.\ Siegel|Meromorphic Funktionen auf kompakten
Mannigfaltigkeiten~|Nachrichten der Akademie der Wissenschaften in
G\"ottingen, Math.-Phys.\ Klasse|4|1955|71-77|

\article Siu74|Y.T.\ Siu|Analyticity of sets associated to Lelong numbers
and the extension of closed positive currents|Invent.\ Math.|27|1974|53-156|

\article Siu93|Y.T.\ Siu|An effective Matsusaka big theorem|Ann.\ Inst.\
Fourier|43|1993|1387-1405|

\miscell Siu94a|Y.T.\ Siu|Effective Very Ampleness|Preprint
{\oldstyle 1994}, to appear in Inventiones Math.|

\miscell Siu94b|Y.T.\ Siu|Very ampleness criterion of double adjoint of
ample line bundles|Preprint {\oldstyle 1994}, to appear in Annals of
Math.\ Studies, volume in honor of Gunning and Kohn, Princeton Univ.\
Press|

\article Sko72a|H.\ Skoda|Sous-ensembles analytiques d'ordre fini ou infini
dans $\bbbc^n$|Bull.\ Soc.\ Math.\ France|100|1972|353-408|

\article Sko72b|H.\ Skoda|Applications des techniques $L^2$ \`a la th\'eorie
des id\'eaux d'une alg\`ebre de fonctions holomorphes avec poids|Ann.\
Scient.\ Ec.\ Norm.\ Sup.\ 4e S\'erie|5|1972|545-579|

\article Sko75|H.\ Skoda|Estimations $L^2$ pour l'op\'erateur $\ovl\partial$
et applications arithm\'etiques|S\'eminaire P.\ Lelong (Analyse), ann\'ee
1975/76, Lecture Notes in Math., Vol.~538,
Springer-Verlag|\rm Berlin|1977|314-323|

\article Sko82|H.\ Skoda|Prolongement des courants positifs ferm\'es de
masse finie|Invent.\ Math.|66|1982|361-376|

\article Som78|A.J.\ Sommese|Submanifolds of abelian varieties|Math.\
Ann.|233|1978|229-256|

\miscell Sug87|K-I.\ Sugiyama|A geometry of K\"ahler cones|preprint
University of Tokyo in Hongo, September {\oldstyle 1987}|

\article Tei79|B.\ Teissier|Du th\'eor\`eme de l'index de Hodge aux
in\'egalit\'es isop\'erim\'etriques|C.\ R.\ Acad.\ Sc.\ Paris,
s\'er.\ A|288|29\ {\rm Janvier}\ 1979|287-289|

\miscell Tei82|B.\ Teissier|Bonnesen-type inequalities in algebraic
geometry|Sem.\ on Diff.\ Geom.\ edited by S.T.\ Yau, {\oldstyle 1982},
Princeton Univ.\ Press, 85-105|

\article Thi67|P.\ Thie|The Lelong number of a point of a complex analytic
set|Math.\ Annalen|172|1967|269-312|

\miscell Tra91|S.\ Trapani|Numerical criteria for the positivity of the
difference of ample divisors|preprint Universit
{\oldstyle 1992}, to appear|

\article Vie82|E.\ Viehweg|Vanishing theorems|J.\ Reine Angew.\
Math.|335|1982|1-8|

\article Yau78|S.T.\ Yau|On the Ricci curvature of a complex K\"ahler
manifold and the complex Monge-Amp\`ere equation I|Comm.\ Pure and Appl.\
Math.|31|1978|339-411|

\vskip30pt
\noindent(July 24, 1994, revised on February 27, 1995)
\vskip10pt
\noindent
Jean-Pierre Demailly\hfill\break
Universit\'e de Grenoble I,\hfill\break
Institut Fourier, URA 188 du CNRS, BP 74,\hfill\break
F-38402 Saint-Martin d'H\`eres, France

\noindent{\it e-mail address:}\/ demailly@fourier.grenet.fr
\end